\documentclass{article}

\usepackage{graphicx}
\usepackage{amssymb}

\usepackage{siunitx}

\usepackage[bbgreekl]{mathbbol}
\DeclareSymbolFontAlphabet{\mathbbm}{bbold}
\DeclareSymbolFontAlphabet{\mathbb}{AMSb}%
\usepackage[pdftex]{hyperref}
\usepackage{xcolor}
\hypersetup{
	colorlinks,
	linkcolor={red!50!black},
	citecolor={blue!90!black},
	urlcolor={blue!80!black}
}
\usepackage{cite}
\usepackage[hmargin={2.5cm,2.5cm}, vmargin={3cm,3cm}, dvips]{geometry}

\usepackage{fancyhdr}
\usepackage{amsmath} 
\usepackage[bbgreekl]{mathbbol}
\DeclareSymbolFontAlphabet{\mathbbm}{bbold}
\DeclareSymbolFontAlphabet{\mathbb}{AMSb}%
\usepackage{amssymb} 
\usepackage{amsthm} 
\usepackage{array}
\usepackage{placeins}
\usepackage{enumerate}
\usepackage{esint}
\usepackage{graphics}
\usepackage{graphicx}
\usepackage{tikz}
\usetikzlibrary{spy}
\usetikzlibrary{patterns}
\usetikzlibrary{arrows}
\usetikzlibrary{calc}

\usetikzlibrary{positioning,shadows,trees,bending}
\usetikzlibrary{arrows.meta,chains}

\usepackage{pgfplots}
\pgfplotsset{compat=1.15} 
\usepackage{caption}
\usepackage{subcaption}
\usepackage{listings}
\usepackage[shortlabels]{enumitem}
\makeatletter
\def\namedlabel#1#2{\begingroup
	#2%
	\def\@currentlabel{#2}%
	\phantomsection\label{#1}\endgroup
}
\usepackage{multirow}
\usepackage{multicol}
\usepackage{wrapfig}

\usepackage{afterpage}
\usepackage{mathrsfs}
\usepackage{mathtools}

\newtheoremstyle{droit}
{}
{}
{\upshape}
{}
{\bfseries}
{}
{ }
{}
\newtheoremstyle{italique}
{}
{}
{\itshape}
{}
{\bfseries}
{}
{ }
{}
\theoremstyle{italique}
\newtheorem{theorem}{Theorem}[section]

\theoremstyle{droit}

\newtheorem{remark}[theorem]{Remark}

\newtheorem{definition}[theorem]{Definition}
\newtheorem{assumption}[theorem]{Assumption}

\newenvironment{repeatassumption}[1]{%
	\renewcommand{\thetheorem}{\ref{#1}}%
	\assumption}
{\endtheorem}

\usepackage{algorithm,algcompatible}
\usepackage{algpseudocode}
\algnewcommand{\lst}{\texttt{lst}}
\algnewcommand{\slst}{\texttt{slst}}
\algnewcommand{\SEND}{\textbf{send}}
\newsavebox{\algleft}
\newsavebox{\algright}

\usepackage{geometry}
\usetikzlibrary{shapes,backgrounds}

\definecolor{darkgreen}{rgb}{0,0.4,0} 
\definecolor{darkbrown}{rgb}{0.5, 0.396, 0.09}

\definecolor{c1}{rgb}{0.0, 0.4196078431372549, 0.6431372549019608}
\definecolor{c2}{rgb}{1.0, 0.5019607843137255, 0.054901960784313725}
\definecolor{c3}{rgb}{0.6705882352941176, 0.6705882352941176,
	0.6705882352941176} \definecolor{c}{rgb}{0.34901960784313724, 0.34901960784313724, 0.34901960784313724}
\definecolor{c4}{rgb}{0.37254901960784315, 0.6196078431372549,
	0.8196078431372549} 
\definecolor{c5}{rgb}{0.5372549019607843, 0.5372549019607843,
	0.5372549019607843} 
\definecolor{c6}{rgb}{1.0, 0.7372549019607844, 0.4745098039215686}
\definecolor{c7}{rgb}{0.8117647058823529, 0.8117647058823529,
	0.8117647058823529}

\newcommand{\setminussign}{{\mathrm r}}
\newcommand{\intersign}{{\mathrm s}}

\newcommand{\vertiii}[1]{{\left\vert\kern-0.25ex\left\vert\kern-0.25ex\left\vert #1 
		\right\vert\kern-0.25ex\right\vert\kern-0.25ex\right\vert}}
\newcommand{\boundarypiece}{{\underline{\smash{\gamma}}}}
\newcommand{\tensorsigma}{{\boldsymbol{\sigma}}}
\newcommand{\tensorepsilon}{{\boldsymbol{\varepsilon}}}

\newcommand{\picsDir}{images}
\newcommand*{\dataPath}{data}
\newcommand{\changes}[1]{{#1}}
\newcommand{\changesbis}[1]{{#1}}
\newcommand{\review}[1]{{#1}}

\begin{document}



\title{\textbf{Analysis-aware defeaturing of complex geometries \changes{with Neumann features}}} 


\author{\vspace{-0.1cm}P. Antol\'in$^{1}$ and O. Chanon$^{2}$\\\vspace{-0.1cm}\footnotesize{pablo.antolin@epfl.ch, ondine.chanon@asc.tuwien.ac.at}\\\vspace{-0.1cm}
	\footnotesize{$^1$ Institute of Mathematics, \'Ecole Polytechnique F\'ed\'erale de Lausanne, Switzerland}\\\vspace{-0.2cm}
	\footnotesize{$^2$ Institute of Analysis and Scientific Computing, TU Wien, Austria}\\
}

\maketitle
\vspace{-0.8cm}
\noindent\rule{\linewidth}{0.4pt}
\thispagestyle{fancy}
\begin{abstract}
  Local modifications of a computational domain are often performed in order to simplify the meshing process and to reduce computational costs and memory requirements. However, removing geometrical features of a domain often introduces a non-negligible error in the solution of a differential problem in which it is defined. In this work, 
  we \changes{extend the results from \cite{paper1defeaturing}} by studying the case of domains containing an arbitrary number of distinct \changes{Neumann} features, and by performing an analysis on Poisson's, linear elasticity, and Stokes' equations. We introduce a simple, computationally cheap, \changes{reliable, and efficient} \textit{a posteriori} estimator of the geometrical defeaturing error. \changes{Moreover, we also} introduce a geometric refinement strategy that accounts for the defeaturing error: Starting from a fully defeatured geometry, the algorithm determines at each iteration step which features need to be added to the geometrical model to reduce the defeaturing error. These important features are then added to the (partially) defeatured geometrical model at the next iteration, until the solution attains a prescribed accuracy. A wide range of \changes{two- and three-dimensional} numerical experiments are finally reported to illustrate this work.
\end{abstract}

{\small \noindent\textit{Keywords}: Geometric defeaturing, geometric refinement, \textit{a posteriori} error estimation, adaptivity, mesh generation.

}

\section{Introduction} \label{s:intro}
With the advance of engineering knowledge, simulations are performed on objects of increasing geometric complexity, nowadays mainly described by three-dimensional Computer-Aided Design (CAD) models. These models often contain a large number of geometric details of different scales, also called geometric features. Unfortunately, the construction of a finite element mesh on such complex domains may fail, or if it does not, 
\review{the mesh generation may be very difficult; see for example~\cite{qian2012automatic} dealing with the complexity arising from an automatic all-hexahedral mesh generation for complex B-Reps (boundary representations). Moreover, the resulting mesh} may require a very large number of elements, therefore leading to simulations which are too costly or even unfeasible. For instance, it has been shown in \cite{white2003meshing,lee2005small} that the cost of the underlying simulation may be increased by up to a factor $10$ in the presence of a single geometric feature of relatively small size. 

\begin{figure}
	\begin{subfigure}{0.49\textwidth}
		\begin{tikzpicture}[spy using outlines={circle,magnification=2,size=4cm, connect spies}] 
		\spy [c2, magnification=3, size=1.4cm] on (2.535,0.05) in node at (2.5,-2.7);
		\spy [c2, magnification=3, size=1.4cm] on (4,-0.9) in node at (4.3,-2.7);
		\node at (0,0) {\includegraphics[width=0.49\textwidth]{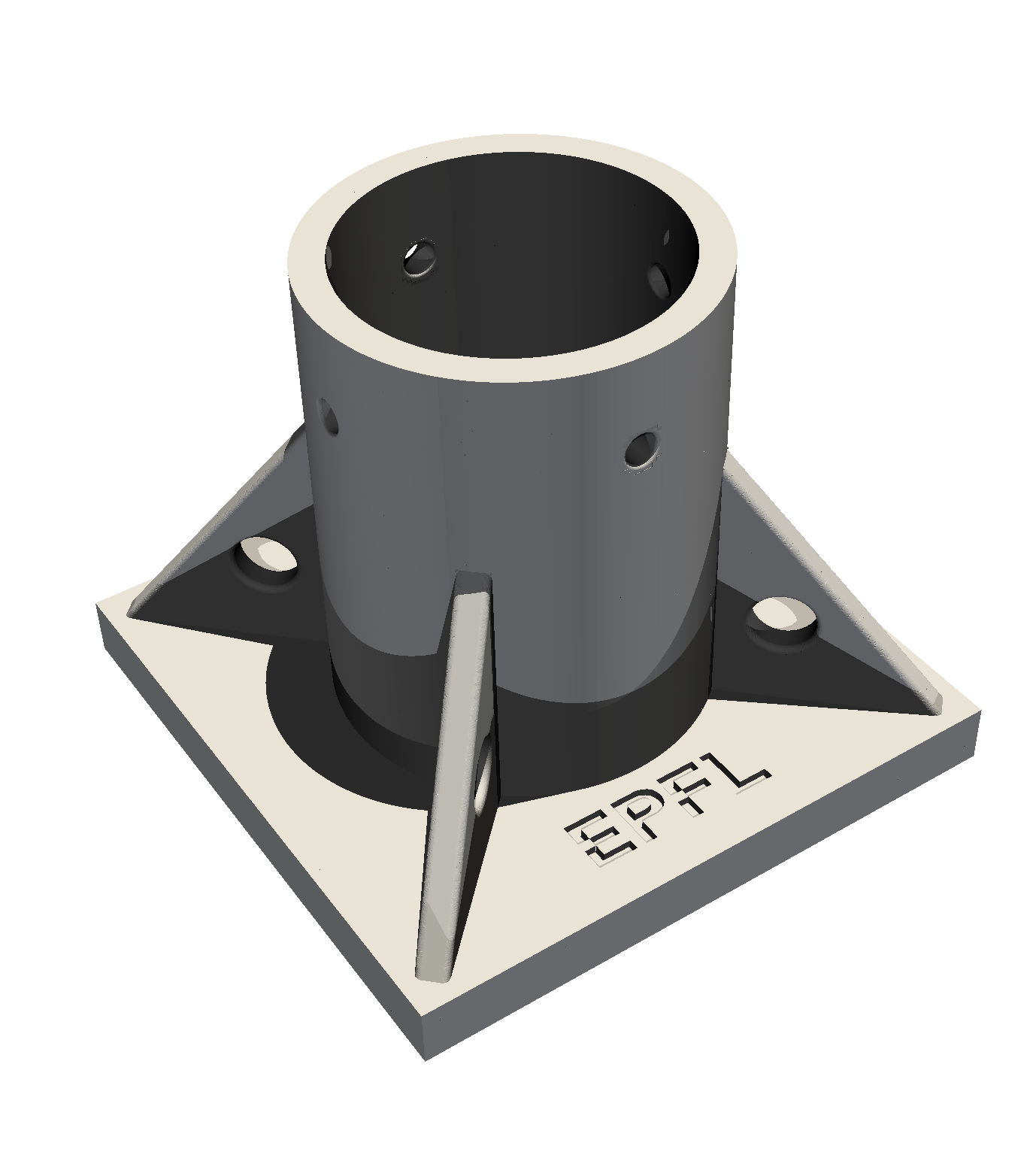}};
		\node at (3.5,0.3) {\includegraphics[width=0.49\textwidth]{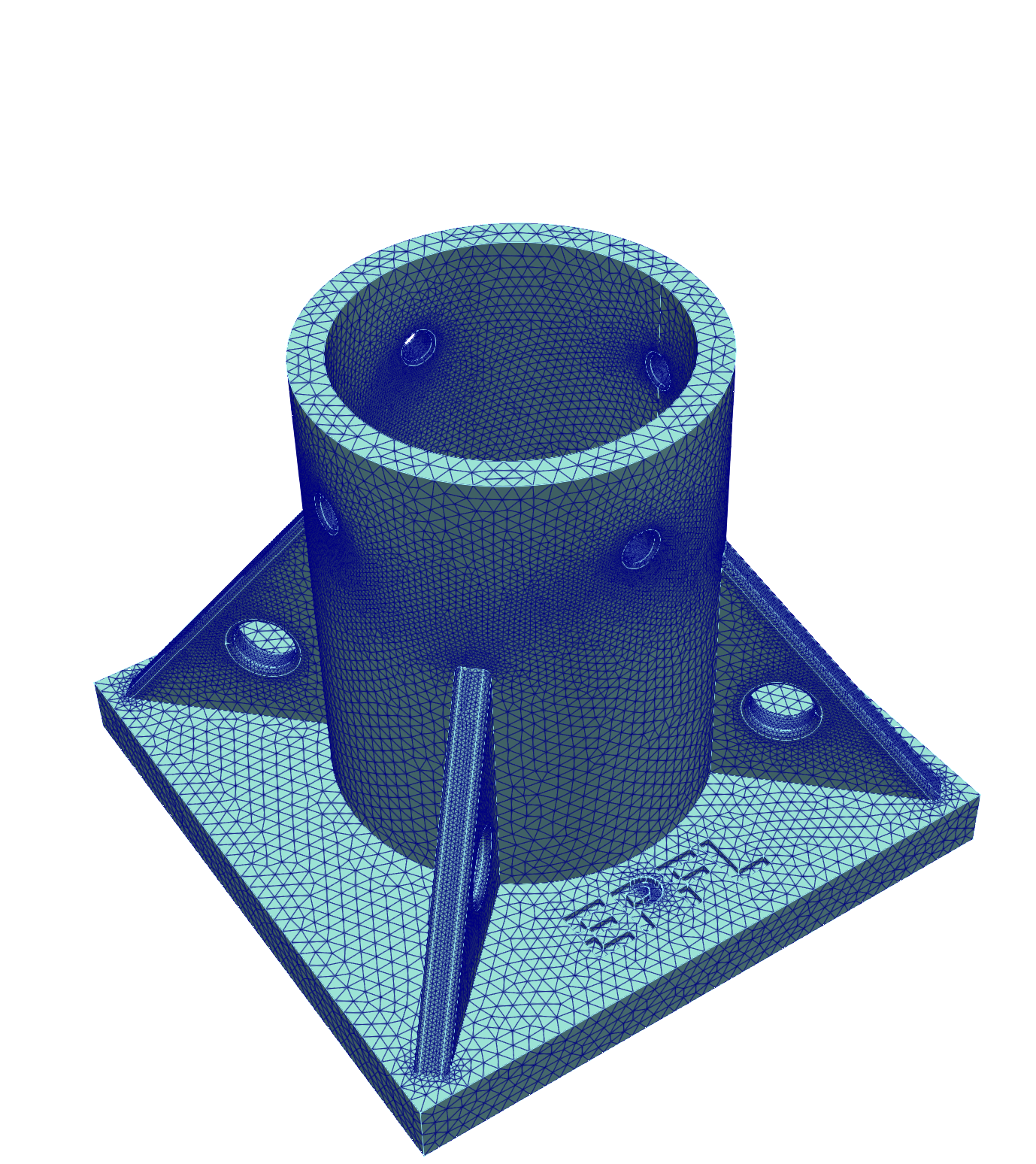}};
		\end{tikzpicture}
		\caption{Original design and mesh (75575 nodes).}
		\label{fig:complexmesh_org}
	\end{subfigure}~
	\begin{subfigure}{0.49\textwidth}
		\begin{tikzpicture}[spy using outlines={circle,magnification=2,size=4cm, connect spies}] 
		\spy [c2, magnification=3, size=1.4cm] on (2.535,0.05) in node at (2.5,-2.7);
		\spy [c2, magnification=3, size=1.4cm] on (4,-0.9) in node at (4.3,-2.7);
		\node at (0,0) {\includegraphics[width=0.49\textwidth]{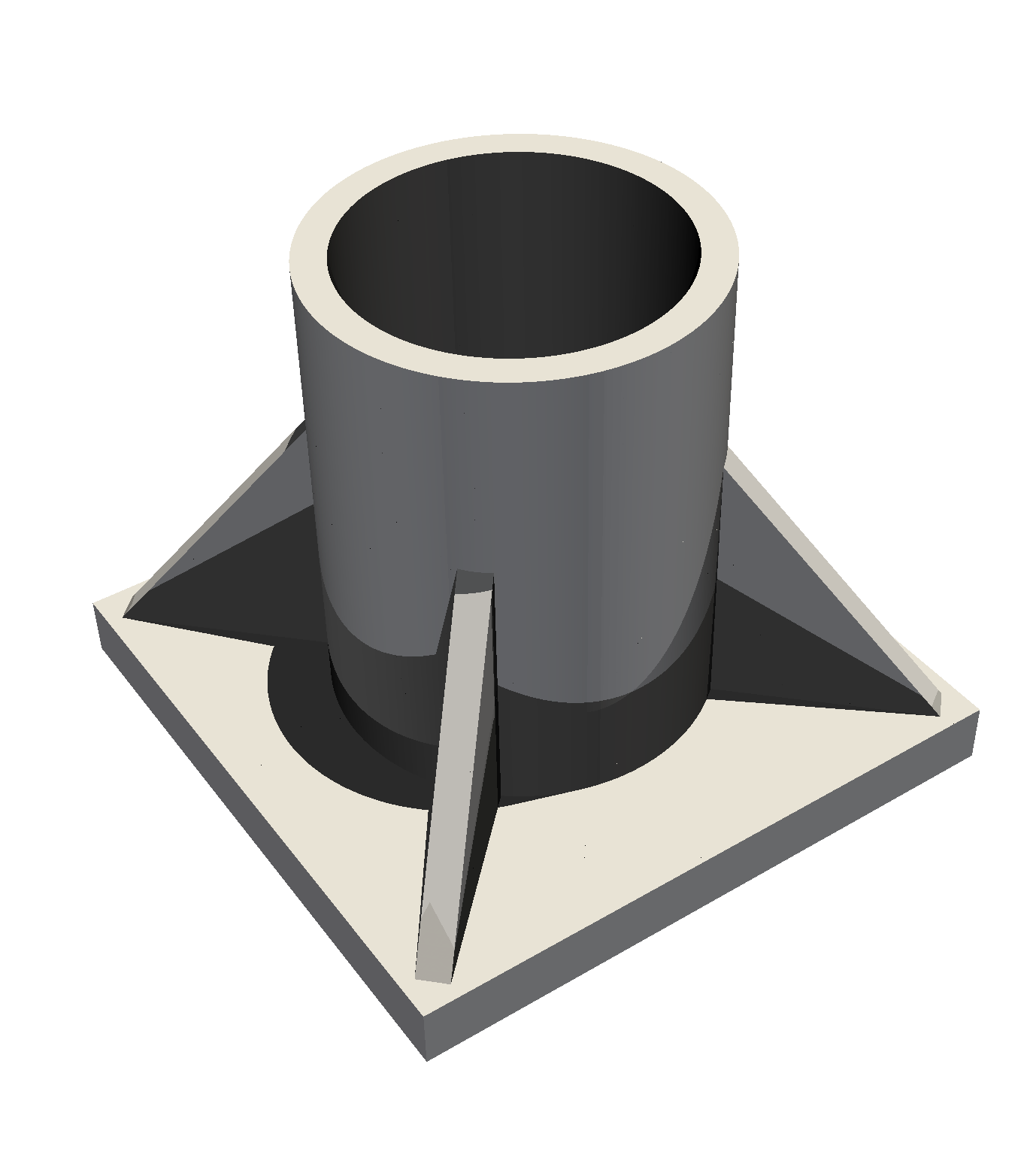}};
		\node at (3.5,0.3) {\includegraphics[width=0.49\textwidth]{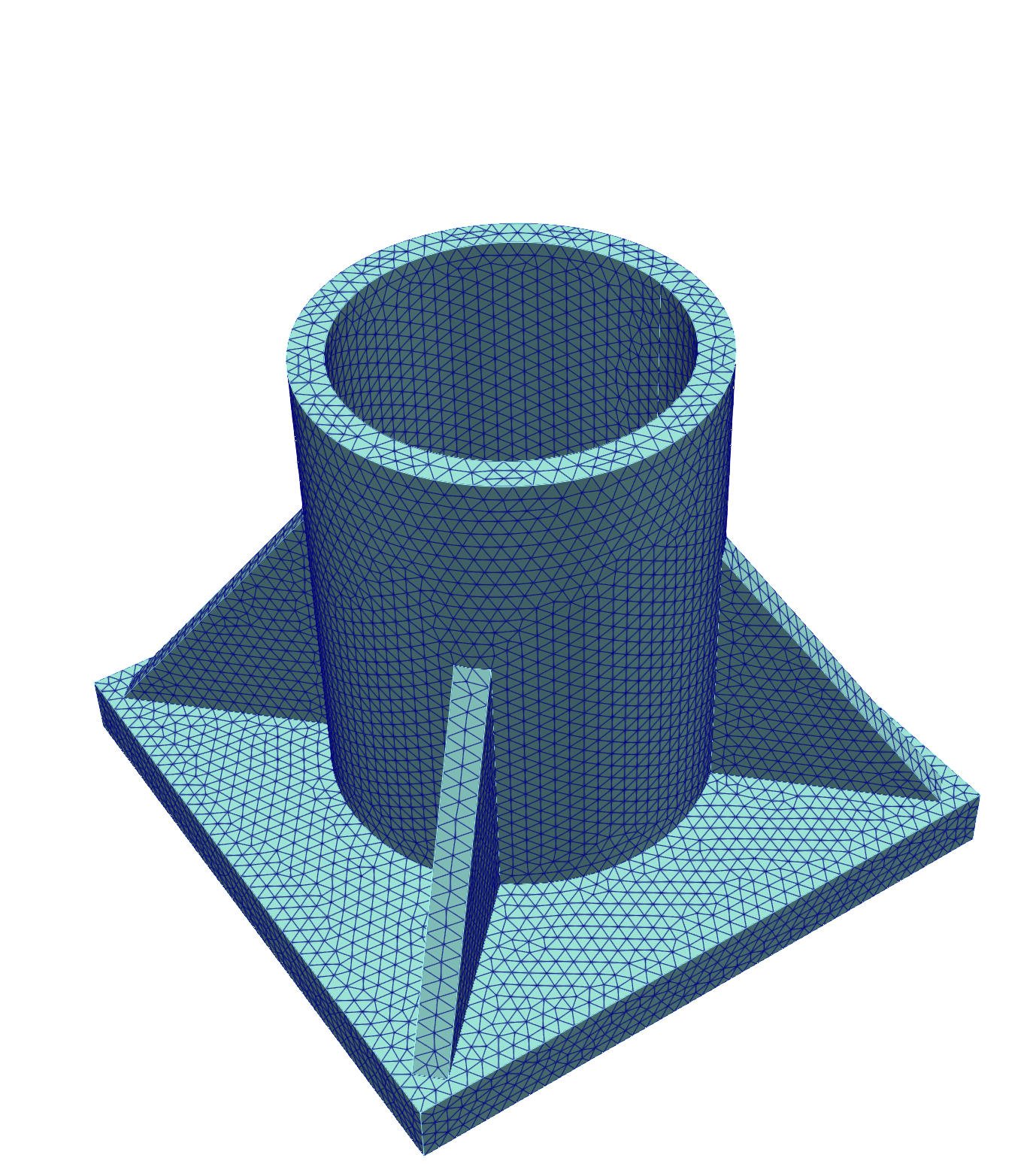}};
		\end{tikzpicture}
		\caption{Defeatured design and mesh (15211 nodes).}
		\label{fig:complexmesh_def}
	\end{subfigure}
	\caption{Example of an original (left) and defeatured CAD designs (right).
		Both finite element meshes were generated with Gmsh~\cite{geuzaine2009gmsh} using the same mesh algorithm and parameters.}
	\label{fig:complexmesh}
\end{figure}

However, depending on the problem at hand, such high model complexity may be unnecessary. That is, the geometric description of the object may require a high number of degrees of freedom, but not all of them are needed to perform an accurate analysis, and taking all of them into account is potentially too costly. To deal with complex geometries and to accelerate the process of analysis-aware geometric design, it is therefore essential to be able to simplify the geometric model, process also known as defeaturing. This is a very common practice among finite element analysts. See, as matter of example, the case illustrated in Figure~\ref{fig:complexmesh}.
There, a CAD design with numerous features as holes, rounds, and a carved logo (Figure~\ref{fig:complexmesh_org}) is defeatured to create a simpler model (Figure~\ref{fig:complexmesh_def}) that is easier to mesh.
Each finite element mesh in Figure~\ref{fig:complexmesh} was generated using Gmsh~\cite{geuzaine2009gmsh} with the same mesh algorithm and parameters.
Nevertheless, the mesh of the original design has 5 times more nodes than the one of the defeatured model.
Likely, most of those extra degrees of freedom, required for correctly representing the geometrical details, will improve little the accuracy of the solution obtained with such mesh, but will increase the computational cost and memory requirements significantly.

Nevertheless, it is important to consider how
such geometrical simplifications will impact the analysis solution, i.e., to control the error introduced by defeaturing, in order to provide an accurate solution of the problem at hand. The literature on the subject is still relatively scarce, and a lot remains to be done.
To estimate the defeaturing error on the solution of a partial differential equation (PDE), some \textit{a posteriori} criteria have been developed: The one introduced in \cite{ferrandesgiannini2009} uses an approximation of the error in energy norm; in \cite{tsa1,tsa2}, an estimator is found using topological sensitivity analysis; adjoint theory is used in the series of works \cite{gopalakrishnan2007formal,gopalakrishnansuresh2008,turevsky2009efficient} to describe the first order defeaturing error on a quantity of interest; an estimator is introduced in \cite{TANG2013413} based on the reciprocal theorem stating the conservation of solution flux in the features; and in the series of works \cite{ligaomartin2011,ligao2011,ligaozhang2012,ligaomartin2013}, defeaturing error is expressed as a modeling error directly on the differential problem, both for negative features (holes) and positive ones (protrusions). In these latter papers, the modeling error is then estimated with the dual weighted residual method \cite{odenprudhomme}. Nevertheless, very few of those works come with a sound mathematical theory, and most of them rely on some heuristics.

In the recent paper \cite{paper1defeaturing}, the authors have tackled this issue: A precise mathematical framework is defined for geometries for which a single feature of very generic shape is removed, and an efficient and reliable \textit{a posteriori} estimator of the defeaturing error is derived in the context of Poisson's equation in the energy norm. 
\changes{Worth mentioning related research papers include \cite{odenvemaganti2000,carstensensauter2004,vegamanti2004} that study heterogeneous and perforated materials, \cite{repinsautersmolianski} that are interested in the error introduced by the approximation of boundary conditions, and \cite{repinsauterbook} that more generally studies modeling errors coming from dimension reduction, homogenization and model simplification.}

\changes{We generalize here the work \cite{paper1defeaturing}. The original contributions of the presented research are the following: 
	\begin{itemize}
		\item We consider geometries for which an \textit{arbitrary number} $N$ of distinct Neumann geometrical features are removed from the computational domain, and we introduce a defeaturing error estimator whose effectivity index is independent from $N$. Note in particular that the features do not need to be small. 
		\item We do not only consider Poisson's equation, but also \textit{linear elasticity and Stokes' problems}. While the considered problems are exclusively linear, the proposed estimator could be extended to encompass nonlinear problems as well.
		\item Based on the proposed defeaturing error estimator, an \textit{adaptive strategy for geometric refinement} is introduced: The algorithm begins with a completely simplified geometry and, at each iteration, identifies the necessary features to be added to the geometric model in order to reduce the defeaturing error.
		\item Finally, more complex \textit{numerical experiments} with engineering interest are presented, in particular the simulation of a three-dimensional CAD design.
	\end{itemize}}
\review{The error estimator formulation depends on the problem at hand. Hence, in this work, different expressions will be provided for the Poisson, Stokes, and linear elasticity problems. In particular, the differential problems and their corresponding defeaturing error estimator depend on the material properties and on the boundary conditions. E.g., some features will play more or less predominant roles depending on the chosen boundary conditions.}

In general, and in practice, differential problems cannot be solved exactly, and thus a numerical method is used to solve them approximately in a finite dimensional space. However, in this article, we neglect the contribution of the error coming from the numerical approximation of the solution, as we concentrate our efforts on the defeaturing error. Deriving an estimator that also includes the numerical error contribution is the subject of a subsequent work, see \cite{adaptivedefeaturing,phdthesis}. \changes{An early method for defeaturing and coarsening, called composite finite elements and developed by Hackbusch and Sauter in \cite{hackbusch1997composite,hackbusch1997composite2}, deserves mentioning. It assumes that defeaturing arises solely from the limitation of the mesh in resolving geometric details within the computational domain. In contrast, our approach aims to treat defeaturing independently of any discretization, allowing us to separate the strict defeaturing error from the numerical approximation error.}

Therefore, in Section~\ref{s:aad}, we first introduce the problem of defeaturing and the corresponding notation that will be used throughout the article, and we precisely define the defeaturing error that we aim at estimating. Then, we state the main results that are obtained, namely the reliability and the efficiency of an \textit{a posteriori} estimator of the defeaturing error whose effectivity index is independent from the number of the features and their size. Subsequently, in \changes{Sections~\ref{s:poisson} and~\ref{s:linelaststokes},} we precisely define the exact and defeatured problems when the differential problem at hand is, respectively, Poisson's \changes{and the linear elasticity or Stokes' equations}, and we propose in each case a defeaturing error estimator. We then introduce in Section~\ref{s:adaptivity} an adaptive geometric refinement strategy driven by the defeaturing error estimators previously defined. We perform in Section~\ref{s:numexp} a validation of the presented theory, of the aforementioned estimators' properties, and of the proposed adaptive strategy, thanks to an extensive set of \changes{two- and three-dimensional} numerical experiments. To perform these tests, we use isogeometric analysis (IGA) \cite{igabook,igabasis} on very fine meshes, in order to have a negligible numerical error with respect to the defeaturing error. Since IGA and defeaturing both pursue the scope of reducing the gap between the design and the analysis phases, IGA is a natural method of choice.
Nevertheless, it is important to remark that the proposed techniques are completely discretization agnostic and any other numerical method could be used. Some conclusions are finally drawn in Section~\ref{s:ccl}.  \changesbis{In the Appendix~\ref{app:proofs}, the reliability and efficiency of the proposed defeaturing error estimator are demonstrated in the framework of Stokes' equations, as the proofs for Poisson's and the linear elasticity equations are very similar, easier, and can be found in details in the corresponding thesis \cite{phdthesis}.} 

In the following, the operator $\lesssim$ is used to mean any inequality which neither depends on the number of features nor on their size, but which can depend on their shape, \review{on their type (i.e., whether they are holes, also called negative features, or protrusions, also called positive features), and on the space dimension $n=2$ or $n=3$}. Moreover, for all $D\subset \mathbb R^n$, and for all $\Lambda \subset \partial D$, we denote by $|D|$ the $n$-dimensional Lebesgue measure of $D$, by $|\Lambda|$ the $(n-1)$-dimensional Hausdorff measure of $\Lambda$, by $\overline{D}$ and $\overline{\Gamma}$ the closure of $D$ and of $\Gamma$, respectively, and by $\text{int}(D)$ and $\text{int}(\Lambda)$ the interior of $D$ and of $\Gamma$, respectively.

\section{Analysis-aware defeaturing} \label{s:aad}
\subsection{Presentation of the problem} \label{ss:pbpres}
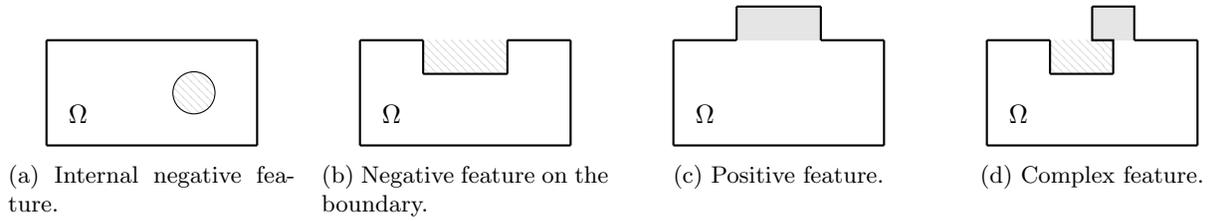
\begin{figure}
	\centering
	\begin{subfigure}[t]{0.23\textwidth}
		\begin{center}
			\begin{tikzpicture}[scale=2.8]
				\draw[thick] (2,0.5) -- (3,0.5) ;
				\draw[thick] (2,0.5) -- (2,1) ;
				\draw[thick] (3,0.5) -- (3,1) ;
				\draw[thick] (2,1) -- (3,1) ;
				\fill[pattern=north west lines, pattern color=gray, opacity=0.5] (2.7,0.75) circle (0.1);
				\draw (2.7,0.75) circle (0.1);
				\draw (2.15,0.65) node{$\Omega$} ;
			\end{tikzpicture}
			\caption{Internal negative feature.}
			\label{fig:internalnegfeat}
		\end{center}
	\end{subfigure}
	~
	\begin{subfigure}[t]{0.23\textwidth}
		\begin{center}
			\begin{tikzpicture}[scale=2.8]
				\draw[thick] (2,0.5) -- (3,0.5) ;
				\fill[pattern=north west lines, pattern color=gray, opacity=0.5] (2.3,1) rectangle (2.7,0.84);
				\draw[thick] (2,0.5) -- (2,1) ;
				\draw[thick] (3,0.5) -- (3,1) ;
				\draw[thick] (2,1) -- (2.3,1) ;
				\draw[thick] (2.7,1) -- (3,1) ;
				\draw[thick] (2.3,1) -- (2.3,0.84) ;
				\draw[thick] (2.7,1) -- (2.7,0.84) ;
				\draw[thick] (2.3,0.84) -- (2.7,0.84) ;
				\draw (2.15,0.65) node{$\Omega$} ;
			\end{tikzpicture}
			\caption{Negative feature on the boundary.}
			\label{fig:bdnegfeat}
		\end{center}
	\end{subfigure}
	~
	\begin{subfigure}[t]{0.23\textwidth}
		\begin{center}
			\begin{tikzpicture}[scale=2.8]
				\fill[gray, fill, opacity=0.2] (0.3,1) rectangle (0.7,1.16);
				\draw[thick] (0,0.5) -- (1,0.5) ;
				\draw[thick] (0,0.5) -- (0,1) ;
				\draw[thick] (1,0.5) -- (1,1) ;
				\draw[thick] (0,1) -- (0.3,1) ;
				\draw[thick] (0.7,1) -- (1,1) ;
				\draw[thick] (0.3,1) -- (0.3,1.16) ;
				\draw[thick] (0.7,1) -- (0.7,1.16) ;
				\draw[thick] (0.3,1.16) -- (0.7,1.16) ;
				\draw (0.15,0.65) node{$\Omega$} ;
			\end{tikzpicture}
			\caption{Positive feature.}
			\label{fig:posfeat}
		\end{center}
	\end{subfigure}
	~
	\begin{subfigure}[t]{0.23\textwidth}
		\begin{center}
			\begin{tikzpicture}[scale=2.8]
				\fill[gray, fill, opacity=0.2] (0.5,1) rectangle (0.7,1.16);
				\fill[pattern=north west lines, pattern color=gray, opacity=0.5] (0.3,1) rectangle (0.6,0.84);
				\draw[thick] (0,0.5) -- (1,0.5) ;
				\draw[thick] (0,0.5) -- (0,1) ;
				\draw[thick] (1,0.5) -- (1,1) ;
				\draw[thick] (0,1) -- (0.3,1) ;
				\draw[thick] (0.7,1) -- (1,1) ;
				\draw[thick] (0.3,1) -- (0.3,0.84) ;
				\draw[thick] (0.3,0.84) -- (0.6,0.84) ;
				\draw[thick] (0.6,0.84) -- (0.6,1) -- (0.5,1) -- (0.5,1.16);
				\draw[thick] (0.5,1.16) -- (0.7,1.16) -- (0.7,1);
				\draw (0.15,0.65) node{$\Omega$} ;
			\end{tikzpicture}
			\caption{Complex feature.}
			\label{fig:complexfeat}
		\end{center}
	\end{subfigure}
	\caption{Domains with different types of geometrical features $F$. In each case, the negative component of $F$ is dashed while its positive component is filled in gray.} \label{fig:featurekinds}
\end{figure}

\changes{In this section, we present the considered problem, following and extending the setting in \cite{paper1defeaturing}.} Let us consider a potentially complicated, \review{connected open Lipschitz} domain $\Omega\subset \mathbb{R}^n$, $n=2$ or $n=3$, on which we want to solve a differential problem $\mathcal P(\Omega)$ which contains some \review{boundary conditions}. More precisely, let us assume that $\Omega$ contains geometrical details of smaller scale also called features, \review{and assume that the feature information of $\Omega$ is known \textit{a priori}}. As illustrated in Figure~\ref{fig:featurekinds}, \changes{we say that a feature $F\subset \mathbb{R}^n$ is \textit{negative} if $\left(\overline{F}\cap\overline{\Omega}\right) \subset \partial \Omega$, \textit{positive} if $F\subset \Omega$, or \textit{complex} if it has both negative and positive components.} 
A negative feature corresponds to a part where some material has been removed (e.g., a hole), a positive feature corresponds to the addition of some material (e.g., a protrusion), and a feature is complex in the most general situation \changes{in which there has been both the addition and the removal of some material}. Note that an internal feature (e.g., an internal hole) is a special case of negative feature. \review{We assume that the positive and the negative part of each feature is a connected open Lipschitz domain of $\mathbb{R}^n$.}

However, solving the given differential problem $\mathcal P(\Omega)$ can be very complicated due to the complexity of $\Omega$, coming from the presence of the features.
Therefore, we solve instead a similar problem but in a \emph{defeatured domain} $\Omega_0$, where the features of $\Omega$ are removed: Holes are filled with some material and protrusions are cut out of the computational domain. This differential problem on $\Omega_0$ is denoted by $\mathcal P(\Omega_0)$ and it is called \emph{defeatured} (or simplified) \emph{problem}.

The defeaturing of the computational domain introduces an error in the problem's solution, and we are interested in controlling the energy norm of this so-called \emph{defeaturing error}. In other words, if $u$ is the solution of the exact problem $\mathcal P(\Omega)$ and $u_0$ is the solution of the defeatured problem $\mathcal P(\Omega_0)$, then we are interested in controlling the error ``$u-u_0$'' in the exact energy norm, which is the energy norm defined by problem $\mathcal P(\Omega)$ and denoted as $\vertiii{\cdot}_\Omega$. Since $u$ is defined in $\Omega$ and $u_0$ is defined in $\Omega_0$, the defeaturing error needs to be more accurately defined. To do so, we need to introduce some further geometric notation. \\

Generalizing the work from \cite{paper1defeaturing}, let $N_f\geq 1$, $N_f\in \mathbb N$, denote the total number of (possibly complex) features of $\Omega$, gathered into the set $\mathfrak{F} := \big\{F^k\big\}_{k=1}^{N_f}$.
Let us make the following assumption on the features.
\begin{assumption} \label{as:separatednew}
	The features in $\mathfrak{F}$ are separated, that is,
	\begin{enumerate}[(a)]
		\item $\overline{F^k}\cap \overline{F^\ell} = \emptyset$ for every $k,\ell=1,\ldots,N_f$, $k\neq\ell$,\label{it:conditionaseparated}
		\item \review{one cannot have an increasingly large number of features that are arbitrarily close to one another.} \label{it:conditionbseparated}
	\end{enumerate}
\end{assumption}

\begin{remark}
	In the currently considered setting, features are discrete objects. \review{Note that:}
	\begin{itemize}
		\item It is always possible to satisfy \review{condition~\ref{it:conditionaseparated}} of Assumption~\ref{as:separatednew} by changing the numbering of the features. Indeed, if there are $k,\ell = 1,\ldots,N_f$, $k\neq\ell$, such that $\overline{F^k}\cap\overline{F^\ell}\neq \emptyset$, then \hbox{$F^{k,\ell}:= \textrm{int}( \overline{F^k} \cup \overline{F^\ell})$} can be considered as a single feature that replaces the two features \hbox{$F^k$ and $F^\ell$}. However, note that the treatment of a geometry in which the boundary is \changes{complicated} everywhere (for instance a fractal-like domain) is not considered here\changes{, see instead \cite{hiptmairshapeapprox,HEYDAROV2022102157,shapederpaper} for first results in this different framework.}
		\item \review{The precise mathematical definition of condition~\ref{it:conditionbseparated} will be given at the end of Section~\ref{ss:notationbdext}, as it requires technical geometric definitions that would hinder here the readability of the main results of this manuscript.}
	\end{itemize}
\end{remark}

\begin{figure}
	\centering
	\begin{subfigure}[t]{0.23\textwidth}
		\begin{center}
			\begin{tikzpicture}[scale=2.8]
				\fill[pattern=north west lines, pattern color=c4, opacity=0.3] (2.3,0.84) rectangle (2.7,1.00);
				\draw[domain=-180:0] plot ({2.5+0.5*cos(\x)}, {1+0.5*sin(\x)}); 
				\draw (2,1) -- (2.3,1) ;
				\draw (2.7,1) -- (3,1) ;
				\draw[c1,thick] (2.3,1) -- (2.3,0.84) ;
				\draw[c1,thick] (2.7,1) -- (2.7,0.84) ;
				\draw[c1,thick] (2.3,0.84) -- (2.7,0.84) ;
				\draw (2.2,0.75) node{$\Omega$} ;
				\draw[c4] (2.5,0.92)node{\small $F=F_\mathrm n$};
				\draw[c1,thick] (2.5,0.84) node[below]{$\gamma=\gamma_\mathrm n$} ;
			\end{tikzpicture}
			\caption{Negative feature in a domain $\Omega$.}
			\label{fig:ex1acomplex}
		\end{center}
	\end{subfigure}
	~
	\begin{subfigure}[t]{0.23\textwidth}
		\begin{center}
			\begin{tikzpicture}[scale=2.8]
				\fill[gray, fill, opacity=0.3] (0.3,1) rectangle (0.7,1.16);
				\draw[domain=-180:0] plot ({0.5+0.5*cos(\x)}, {1+0.5*sin(\x)});
				\draw (0,1) -- (0.3,1) ;
				\draw (0.7,1) -- (1,1) ;
				\draw[c1,thick,dash dot] (0.3,1) -- (0.3,1.16) ;
				\draw[c1,thick,dash dot] (0.7,1) -- (0.7,1.16) ;
				\draw[c1,thick] (0.5,1.15) node[above]{$\gamma=\gamma_\mathrm p$} ;
				\draw[c1,thick,dash dot] (0.3,1.16) -- (0.7,1.16) ;
				\draw (0.2,0.75) node{$\Omega$} ;
				\draw[gray] (0.5,1.08)node{\small $F=F_\mathrm p$};
			\end{tikzpicture}
			\caption{Positive feature in a domain $\Omega$.}
			\label{fig:ex1ccomplex}
		\end{center}
	\end{subfigure}
	~
	\begin{subfigure}[t]{0.23\textwidth}
		\begin{center}
			\begin{tikzpicture}[scale=2.8]
				\fill[pattern=north west lines, pattern color=c4, opacity=0.3] (0.3,0.84) rectangle (0.5,1);
				\fill[gray, fill, opacity=0.3] (0.5,1) rectangle (0.7,1.16);
				\draw[domain=-180:0] plot ({0.5+0.5*cos(\x)}, {1+0.5*sin(\x)});
				\draw (0,1) -- (0.3,1) ;
				\draw (0.7,1) -- (1,1) ;
				\draw[c1,thick] (0.3,1) -- (0.3,0.84) ;
				\draw[c1,thick,dash dot] (0.5,1.16) -- (0.5,1) ;
				\draw[c1,thick] (0.5,1)--(0.5,0.84) ;
				\draw[c1,thick,dash dot] (0.7,1) -- (0.7,1.16) ;
				\draw[c1,thick] (0.4,0.84) node[below]{$\gamma_\mathrm n$} ;
				\draw[c1,thick] (0.6,1.16) node[above]{$\gamma_\mathrm p$} ;
				\draw[c1,thick,dash dot] (0.5,1.16) -- (0.7,1.16) ;
				\draw[c1,thick] (0.3,0.84) -- (0.5,0.84) ;
				\draw (0.2,0.75) node{$\Omega$} ;
				\draw[gray] (0.6,1.08)node{\small $F_\mathrm p$};
				\draw[c4] (0.4,0.92)node{\small $F_\mathrm n$};
			\end{tikzpicture}
			\caption{General complex feature in a domain $\Omega$.}
			\label{fig:ex1dcomplex}
		\end{center}
	\end{subfigure}
	~
	\begin{subfigure}[t]{0.23\textwidth}
		\begin{center}
			\begin{tikzpicture}[scale=2.8]
				\draw[domain=-180:0] plot ({0.5+0.5*cos(\x)}, {1+0.5*sin(\x)});
				\draw (0,1) -- (0.3,1) ;
				\draw (0.7,1) -- (1,1) ;
				\draw[c2,thick] (0.3,1) -- (0.7,1) ;
				\draw[c2,thick] (0.3,0.98) -- (0.3, 1.02);
				\draw[c2,thick] (0.7,0.98) -- (0.7, 1.02);
				\draw[c2,thick] (0.5,1) node[below]{$\gamma_0$} ;
				\draw (0.2,0.75) node{$\Omega_0$} ;
			\end{tikzpicture}
			\caption{Simplified domain for cases in (a)--(c), (e)--(f).}
			\label{fig:ex1bcomplex}
		\end{center}
	\end{subfigure}
	~
	\begin{subfigure}[t]{0.23\textwidth}
		\begin{center}
			\begin{tikzpicture}[scale=2.8]
				\fill[pattern=north west lines, pattern color=c4, opacity=0.3] (0.3,0.84) rectangle (0.6,1);
				\fill[gray, fill, opacity=0.3] (0.5,1) rectangle (0.7,1.2);
				\draw[domain=-180:0] plot ({0.5+0.5*cos(\x)}, {1+0.5*sin(\x)});
				\draw (0,1) -- (0.3,1) ;
				\draw (0.7,1) -- (1,1) ;
				\draw[c1,thick] (0.3,1) -- (0.3,0.84) ;
				\draw[c1,thick,dash dot] (0.5,1.2) -- (0.5,1);
				\draw[c1,thick,dash dot] (0.5,1) -- (0.6,1);
				\draw[c1,thick](0.6,1) -- (0.6,0.84) ;
				\draw[c1,thick,dash dot] (0.7,1) -- (0.7,1.2) ;
				\draw[c1,thick] (0.45,0.84) node[below]{$\gamma_\mathrm n$} ;
				\draw[c1,thick] (0.6,1.18) node[above]{$\gamma_\mathrm p$} ;
				\draw[c1,thick,dash dot] (0.5,1.2) -- (0.7,1.2) ;
				\draw[c1,thick] (0.3,0.84) -- (0.6,0.84) ;
				\draw (0.2,0.75) node{$\Omega$} ;
				\draw[gray] (0.6,1.1)node{\small $F_\mathrm p$};
				\draw[c4] (0.45,0.92)node{\small $F_\mathrm n$};
			\end{tikzpicture}
			\caption{General complex feature in a domain $\Omega$.}
			\label{fig:ex1fcomplex}
		\end{center}
	\end{subfigure}
	~
	\begin{subfigure}[t]{0.23\textwidth}
		\begin{center}
			\begin{tikzpicture}[scale=2.8]
				\fill[gray, fill, opacity=0.3] (0.4,0.87) rectangle (0.5,1.2);
				\fill[gray, fill, opacity=0.3] (0.5,1) rectangle (0.7,1.2);
				\fill[pattern=north west lines, pattern color=c4, opacity=0.5](0.3,0.84) rectangle (0.6,1);
				\draw[domain=-180:0] plot ({0.5+0.5*cos(\x)}, {1+0.5*sin(\x)});
				\draw (0,1) -- (0.3,1) ;
				\draw (0.7,1) -- (1,1) ;
				\draw[c1,thick] (0.3,1) -- (0.3,0.84) ;
				\draw[c1,thick,dash dot] (0.4,1.2) -- (0.4,0.87);
				\draw[c1,thick,dash dot] (0.5,1) -- (0.5,0.87);
				\draw[c1,thick,dash dot] (0.4,0.87) -- (0.5,0.87);
				\draw[c1,thick,dash dot] (0.5,1) -- (0.6,1);
				\draw[c1,thick](0.6,1) -- (0.6,0.84) ;
				\draw[c1,thick,dash dot] (0.7,1) -- (0.7,1.2) ;
				\draw[c1,thick] (0.45,0.84) node[below]{$\gamma_\mathrm n$} ;
				\draw[c1,thick] (0.55,1.18) node[above]{$\gamma_\mathrm p$} ;
				\draw[c1,thick,dash dot] (0.4,1.2) -- (0.7,1.2) ;
				\draw[c1,thick] (0.3,0.84) -- (0.6,0.84) ;
				\draw (0.2,0.75) node{$\Omega$} ;
				\draw[gray] (0.5,1.1)node{\small$F_\mathrm p$};
				\draw[c4] (0.15,0.92)node{\small $F_\mathrm n$};
				\draw[c4] (0.35,0.92) -- (0.23,0.92);
			\end{tikzpicture}
			\caption{General complex feature in a domain $\Omega$.}
			\label{fig:ex1ecomplex}
		\end{center}
	\end{subfigure}
	~
	\begin{subfigure}[t]{0.23\textwidth}
		\begin{center}
			\begin{tikzpicture}[scale=2.8]
				\draw[domain=-180:0] plot ({0.5+0.5*cos(\x)}, {1+0.5*sin(\x)});
				\draw (0,1) -- (0.3,1) ;
				\draw (0.7,1) -- (1,1) ;
				\draw[c2,thick] (0.3,1) -- (0.6, 1);
				\draw[c2,thick] (0.6,0.98) -- (0.6, 1.02);
				\draw[c2,thick] (0.3,0.98) -- (0.3, 1.02);
				\draw[c2,thick] (0.7,0.98) -- (0.7, 1.02);
				\draw[c2,thick,dash dot] (0.6,1) -- (0.7,1) ;
				\draw[c2,thick] (0.45,1) node[below]{$\gamma_{0,\mathrm n}$} ;
				\draw[c2,thick] (0.65,1) node[above]{$\gamma_{0,\mathrm p}$} ;
				\draw (0.2,0.75) node{$\Omega_0$} ;
			\end{tikzpicture}
			\caption{Simplified cases for domains in (e) and (f).}
			\label{fig:ex1gcomplex}
		\end{center}
	\end{subfigure}
	~
	\begin{subfigure}[t]{0.23\textwidth}
		\begin{center}
			\begin{tikzpicture}[scale=2.8]
				\fill[gray, fill, opacity=0.3] (0.3,1.05) rectangle (0.7,1.2);
				\fill[gray, fill, opacity=0.3] (0.3,1) rectangle (0.4,1.05);
				\fill[gray, fill, opacity=0.3] (0.6,1) rectangle (0.7,1.05);
				\draw[domain=-180:0] plot ({0.5+0.5*cos(\x)}, {1+0.5*sin(\x)});
				\draw (0,1) -- (0.3,1) ;
				\draw (0.7,1) -- (1,1) ;
				\draw (0.4,1) -- (0.6,1);
				\draw[c1,thick,dash dot] (0.3,1) -- (0.3, 1.2);
				\draw[c1,thick,dash dot] (0.3,1.2) -- (0.7, 1.2);
				\draw[c1,thick,dash dot] (0.7,1) -- (0.7, 1.2);
				\draw[c1,thick,dash dot] (0.4,1) -- (0.4, 1.05);
				\draw[c1,thick,dash dot] (0.4,1.05) -- (0.6, 1.05);
				\draw[c1,thick,dash dot] (0.6,1) -- (0.6, 1.05);
				\draw[c1,thick] (0.5,1.18) node[above]{$\gamma = \gamma_{\mathrm p}$} ;
				\draw (0.2,0.75) node{$\Omega_0$} ;
				\draw[gray] (0.5,1.1)node{\small $F=F_\mathrm p$};
				\draw[c2,thick,dash dot] (0.6,1) -- (0.7,1) ;
				\draw[c2,thick,dash dot] (0.3,1) -- (0.4,1) ;
				\draw[c2] (0.49,0.9) -- (0.35,0.98);
				\draw[c2] (0.51,0.9) -- (0.65,0.98);
				\draw[c2] (0.545,0.91) node[below]{$\gamma_0 = \gamma_{0,\mathrm p}$} ;
			\end{tikzpicture}
			\caption{Positive feature in domain $\Omega := \text{int}\left(\overline{\Omega_0} \cup \overline{F}\right)$.}
			\label{fig:ex1hcomplex}
		\end{center}
	\end{subfigure}
	\caption{Geometries containing different types of features.} \label{fig:exfeatcomplex}
\end{figure}
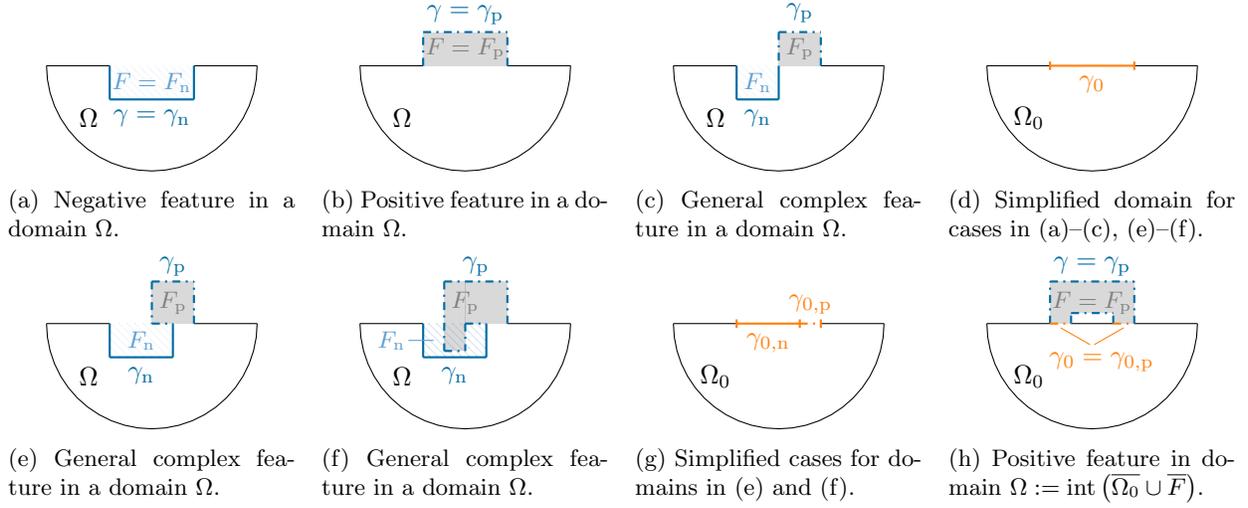

Since the features in $\mathfrak F$ are assumed to be generally complex, this means that, for all $k=1,\ldots,N_f$, the feature $F^k\in\mathfrak F$ is an open domain which is composed of a (not necessarily connected) negative component $F_\mathrm n^k$, a hole, and a (not necessarily connected) positive component $F^k_\mathrm p$, a protrusion, that can have a non-empty intersection, see Figure~\ref{fig:exfeatcomplex}. More precisely,
$$F^k = \text{int}\left(\overline{F_\mathrm p^k} \cup \overline{F_\mathrm n^k}\right),$$
where $F_\mathrm n^k$ and $F_\mathrm p^k$ are open domains such that, if we define
\begin{equation}\label{eq:defOmegastar}
	F_\mathrm p := \text{int}\left(\bigcup_{k=1}^{N_f} \overline{F_\mathrm p^k}\right), \quad F_\mathrm n := \text{int}\left(\bigcup_{k=1}^{N_f} \overline{F_\mathrm n^k}\right), \quad \Omega_\star := \Omega \setminus \overline{F_\mathrm p},
\end{equation}
then
\begin{align*}
	 & F_\mathrm p \subset \Omega\ \textrm{and}\ \left(\overline{F_\mathrm n} \cap \overline{\Omega_\star} \right) \subset \partial \Omega_\star.
\end{align*}
In this setting, we define the defeatured geometry by
\begin{equation} \label{eq:defomega0multi}
	\Omega_0 := \text{int}\left( \overline{\Omega_\star} \cup \overline{F_\mathrm n} \right) \subset \mathbb{R}^n.
\end{equation}
Note that in the case in which the exact domain $\Omega$ contains a single feature $F$, i.e. $N_f = 1$, then
\begin{itemize}
	\item if $F$ is negative, i.e., $F_\mathrm p = \emptyset$ and $F_\mathrm n = F$, then ${\Omega_0} := \mathrm{int}\left(\overline\Omega \cup \overline F\right)$,
	\item if $F$ is positive, i.e., $F_\mathrm p = F$ and $F_\mathrm n = \emptyset$,  then $\Omega_0 := \Omega \setminus \overline F$,
	\item if $F$ is complex, definitions \eqref{eq:defOmegastar} and \eqref{eq:defomega0multi} and apply,
\end{itemize}
therefore generalizing the definition of $\Omega_0$ given in \cite{paper1defeaturing}.

\begin{figure}
	\centering
	\begin{subfigure}[t]{0.3\textwidth}
		\begin{center}
			\begin{tikzpicture}[scale=2.2]
				\draw[thick] (1,0) -- (2,0) ;
				\draw[thick] (1,0) -- (1,1) ;
				\draw[thick] (2,0) -- (2,0.5) ;
				\draw[thick] (1,1) -- (1.5,1) ;
				\draw[thick] (1.9,0.5) -- (2,0.5);
				\draw[thick] (1.5,0.9) -- (1.5,1);
				\draw[thick] (1.5,0.9) -- (1.7,0.9);
				\draw [thick](1.7,0.9) -- (1.9,0.5);
				\draw (1.45,0.45) node{$\Omega$} ;
			\end{tikzpicture}
			\caption{Exact domain $\Omega$.}
			\label{fig:Omegaexactnoteasy}
		\end{center}
	\end{subfigure}
	~
	\begin{subfigure}[t]{0.3\textwidth}
		\begin{center}
			\begin{tikzpicture}[scale=2.2]
				\fill[pattern=north east lines, opacity=0.4] (1.5,1) -- (1.5,0.9) -- (1.7,0.9) -- (1.9,0.5) -- (2,0.5) -- (2,1) -- cycle;
				\draw[thick] (1,0) -- (2,0) ;
				\draw[thick] (1,0) -- (1,1) ;
				\draw[thick] (2,0) -- (2,0.5) ;
				\draw[thick] (1,1) -- (1.5,1) ;
				\draw[thick] (1.9,0.5) -- (2,0.5);
				\draw[thick] (1.5,0.9) -- (1.5,1);
				\draw[thick] (1.5,0.9) -- (1.7,0.9);
				\draw[thick] (1.7,0.9) -- (1.9,0.5);
				\draw[thick] (1.5,1) -- (2,1);
				\draw[thick] (2,1) -- (2,0.5);
				\draw (1.875,0.85) node{$F$};
				\draw (1.45,0.45) node{$\Omega$} ;
			\end{tikzpicture}
			\caption{Possible simplified domain $\Omega_0:=\mathrm{int}(\overline{\Omega}\cup\overline{F})$ with a negative feature $F$.}
			\label{fig:negfeatOmeganoteasy}
		\end{center}
	\end{subfigure}
	~
	\begin{subfigure}[t]{0.3\textwidth}
		\begin{center}
			\begin{tikzpicture}[scale=2.2]
				\fill[gray!20] (1.5,0.5) -- (1.9,0.5) -- (1.7,0.9) -- (1.5,0.9) -- cycle;
				\draw[thick] (1,0) -- (2,0) ;
				\draw[thick] (1,0) -- (1,1) ;
				\draw[thick] (2,0) -- (2,0.5) ;
				\draw[thick] (1,1) -- (1.5,1) ;
				\draw[thick] (1.5,0.5) -- (1.9,0.5);
				\draw[thick] (1.9,0.5) -- (2,0.5);
				\draw[thick] (1.5,0.9) -- (1.5,1);
				\draw[thick] (1.5,0.5) -- (1.5,0.9);
				\draw[thick] (1.5,0.9) -- (1.7,0.9);
				\draw[thick] (1.7,0.9) -- (1.9,0.5);
				\draw (1.25,0.2) node{$\Omega_0$} ;
				\draw (1.65,0.67) node{$F$} ;
			\end{tikzpicture}
			\caption{Possible simplified domain $\Omega_0= \Omega \setminus \overline F$ with a positive feature $F$.}
			\label{fig:posfeatOmeganoteasy}
		\end{center}
	\end{subfigure}
	\caption{Exact geometry $\Omega$ and different possible defeatured geometries $\Omega_0$.} \label{fig:exnegposnoteasy}
\end{figure}
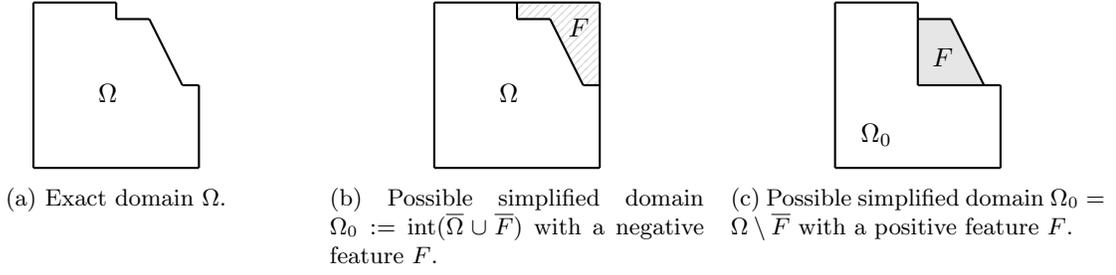

\begin{remark}
	Given a complicated geometry $\Omega$ without any further information, one cannot always easily tell whether the features it contains are negative or positive, see Figure~\ref{fig:exnegposnoteasy}. Therefore, this is often a choice that the user needs to make, based on the available geometric information at hand. If one has access to a simplified geometry, for instance thanks to the history of CAD operations from which the exact geometry $\Omega$ is built, then it is possible to define the features from $\Omega$ and $\Omega_0$, instead of defining $\Omega_0$ from $\Omega$ and the features. The identification of features in a given geometry and the construction of a corresponding simplified geometric model can be complicated tasks, see \cite{surveymodelsimpl} for a review of possible techniques. However, this goes beyond the scope of this work, which supposes at its roots that the feature information is known \textit{a priori}.
\end{remark}

\begin{figure}
	\centering
	\begin{subfigure}[t]{0.23\textwidth}
		\begin{center}
			\begin{tikzpicture}[scale=3.5, every node/.style={scale=0.75}]
				\fill[c3,opacity=0.2,domain=-180:0] (2.7,1) -- (3,1) -- (3,0.9) -- plot ({2.5+0.5*cos(\x)}, {0.9+0.5*sin(\x)}) -- (2,0.9) -- (2,1) -- (2.3,1) -- (2.3,1.3) -- (2.4,1.3) -- (2.6,1.1) -- (2.6,1) -- (2.4,1) -- (2.4,0.8) -- (2.7,0.8) -- cycle;
				\fill[white] (2.4,0.8) -- (2.7,0.8) -- (2.7,1) -- (2.4,1) -- cycle;
				\draw[thick] (2.7,1) -- (3,1) -- (3,0.9);
				\draw[thick] (2,0.9) -- (2,1) -- (2.3,1) ;
				\draw[thick,domain=-180:0] plot ({2.5+0.5*cos(\x)}, {0.9+0.5*sin(\x)});
				\draw[c1,thick,dashed] (2.4,1) -- (2.4,0.8) -- (2.7,0.8) -- (2.7,1);
				\draw[c1,thick] (2.4,1) -- (2.6,1) -- (2.6,1.1) -- (2.4,1.3) -- (2.3,1.3) -- (2.3,1);
				\draw (2.23,0.7) node[scale=1.25]{$\Omega$} ;
				\draw[c1,thick] (2.3,1.15) node[left]{$\gamma_\mathrm p$} ;
				\draw[c1,thick] (2.55,0.8) node[below]{$\gamma_\mathrm n$} ;
			\end{tikzpicture}
			\caption{Domain $\Omega$ with a complex feature.}
		\end{center}
	\end{subfigure}
	~
	\begin{subfigure}[t]{0.23\textwidth}
		\begin{center}
			\begin{tikzpicture}[scale=3.5, every node/.style={scale=0.75}]
				\fill[c3, opacity=0.4] (2.3,1) rectangle (2.6,1.3);
				\fill[c3, opacity=0.1,domain=-180:0] (3,0.9) -- (3,1) -- (2,1) -- (2,0.9) -- plot ({2.5+0.5*cos(\x)}, {0.9+0.5*sin(\x)}); 
				\draw[thick] (2,0.9) -- (2,1) ;
				\draw[thick] (3,0.9) -- (3,1) ;
				\draw[thick,domain=-180:0] plot ({2.5+0.5*cos(\x)}, {0.9+0.5*sin(\x)});
				\draw[thick] (2,1) -- (2.3,1) ;
				\draw[thick] (2.7,1) -- (3,1) ;
				\draw[c2] (2.6,1.2) node[right,scale=1.25]{$\tilde \gamma$};
				\draw (2.46,1.16)node{$\tilde F_\mathrm p$};
				\draw[c2,thick] (2.6,1.1) -- (2.6,1.3) -- (2.4,1.3);
				\draw[c1,dashed] (2.4,1) -- (2.4,0.8) ;
				\draw[c1,dashed] (2.7,1) -- (2.7,0.8) ;
				\draw[c1,dashed] (2.4,0.8) -- (2.7,0.8) ;
				\draw[c1,dashed] (2.6,1.1) -- (2.4,1.3);
				\draw (2.5,0.75) node[scale=1.25]{$\review{\Omega_0}$} ;
				\draw[c2,thick,densely dotted] (2.3,1) -- (2.4,1);
				\draw[c2] (2.33,1) node[below]{$\gamma_{0,\mathrm p}$};
				\draw[c2] (2.55,1) node[below]{$\gamma_{0,\mathrm n}$};
				\draw[c2, thick] (2.4,1) -- (2.6,1);
				\draw[c2, thick] (2.6,1) -- (2.7,1);
				\draw[c4,thick] (2.4,1.0025) -- (2.6,1.0025) -- (2.6,1.1);
				\draw[c4,thick] (2.4,1.3) -- (2.3,1.3) -- (2.3,1);
				\draw[c2,thick] (2.4,0.9975) -- (2.6,0.9975);
				\draw[c2,thick] (2.6,0.9975) -- (2.7,1);
				\draw[c4] (2.65,1.05) node{$\gamma_\intersign$};
			\end{tikzpicture}
			\caption{Defeatured domain $\Omega_0$ and simplified positive component $\tilde F_\mathrm p$.} \label{fig:boundingbox}
		\end{center}
	\end{subfigure}
	~
	\begin{subfigure}[t]{0.23\textwidth}
		\begin{center}
			\begin{tikzpicture}[scale=3.5, every node/.style={scale=0.75}]
				\draw[white] (2,1.3)--(2,1.35);
				\fill[c3, opacity=0.1,domain=-180:0] (2,0.9) -- (2,1) -- (3,1) -- (3,0.9) -- plot ({2.5+0.5*cos(\x)}, {0.9+0.5*sin(\x)}) -- cycle;
				\fill[c3, opacity=0.55] (2.4,0.8) rectangle (2.7,1.00);
				\fill[c3, opacity=0.4] (2.3,1) -- (2.6,1) -- (2.6,1.1) -- (2.4,1.3) -- (2.3,1.3);
				\fill[c3, opacity=0.8] (2.6,1.1) -- (2.6,1.3) -- (2.4,1.3);
				\draw[c2,thick] (2.6,1.1) -- (2.6,1.3) -- (2.4,1.3);
				\draw[c1] (2.5,1.14) node{$\gamma_\setminussign$};
				\draw (2.54,1.24)node{$G_\mathrm p$};
				\draw (2.55,0.9)node{$F_\mathrm n$};
				\draw (2.4,1.1)node{$F_\mathrm p$};
				\draw[c2,thick,densely dotted] (2.3,1) -- (2.4,1);
				\draw[thick] (2,0.9) -- (2,1) ;
				\draw[thick] (3,0.9) -- (3,1) ;
				\draw[thick,domain=-180:0] plot ({2.5+0.5*cos(\x)}, {0.9+0.5*sin(\x)});
				\draw[thick] (2,1) -- (2.3,1) ;
				\draw[thick] (2.7,1) -- (3,1) ;
				\draw[c1,thick,dashed] (2.4,1) -- (2.4,0.8) ;
				\draw[c1,thick,dashed] (2.7,1) -- (2.7,0.8) ;
				\draw[c1,thick,dashed] (2.4,0.8) -- (2.7,0.8) ;
				\draw[c1,thick,dashed] (2.6,1.1) -- (2.4,1.3);
				\draw[c4,thick] (2.4,1.0025) -- (2.6,1.0025) -- (2.6,1.1);
				\draw[c4,thick] (2.4,1.3) -- (2.3,1.3) -- (2.3,1);
				\draw[c4] (2.25,1.15) node{$\gamma_\intersign$};
				\draw (2.23,0.7) node[scale=1.25]{$\Omega_\star$} ;
				\draw[c2,thick] (2.4,0.9975) -- (2.6,0.9975);
				\draw[c2,thick] (2.6,0.9975) -- (2.7,1);
			\end{tikzpicture}
			\caption{Domains $\Omega_\star$, $F_\mathrm n$, $F_\mathrm p$ and $G_\mathrm p$.}\label{fig:extentionnotation_Gp}
		\end{center}
	\end{subfigure}
	~
	\begin{subfigure}[t]{0.23\textwidth}
		\begin{center}
			\begin{tikzpicture}[scale=3.5, every node/.style={scale=0.75}]
				\fill[c3,opacity=0.25,domain=-180:0] (2.7,1) -- (3,1) -- (3,0.9) -- plot ({2.5+0.5*cos(\x)}, {0.9+0.5*sin(\x)}) -- (2,0.9) -- (2,1) -- (2.3,1) -- (2.3,1.3) -- (2.6,1.3) -- (2.6,1) -- (2.4,1) -- (2.4,0.8) -- (2.7,0.8) --cycle;
				\fill[white] (2.4,0.8) -- (2.7,0.8) -- (2.7,1) -- (2.4,1) -- cycle;
				\draw[c3,thick] (2.6,1) -- (2.6,1.3) -- (2.3,1.3) -- (2.3,1.25);
				\draw[c2,densely dotted] (2.3,1) -- (2.4,1);
				\draw[thick] (2,0.9) -- (2,1) ;
				\draw[thick] (3,0.9) -- (3,1) ;
				\draw[thick,domain=-180:0] plot ({2.5+0.5*cos(\x)}, {0.9+0.5*sin(\x)});
				\draw[thick] (2,1) -- (2.3,1) ;
				\draw[thick] (2.7,1) -- (3,1) ;
				\draw[c1,dashed] (2.4,1) -- (2.4,0.8) ;
				\draw[c1,dashed] (2.7,1) -- (2.7,0.8) ;
				\draw[c1,dashed] (2.4,0.8) -- (2.7,0.8) ;
				\draw[c1,dashed] (2.6,1.1) -- (2.4,1.3);
				\draw[c4,thick] (2.4,1.0025) -- (2.6,1.0025) -- (2.6,1.1);
				\draw[c4,thick] (2.4,1.3) -- (2.3,1.3) -- (2.3,1);
				\draw[c2,thick] (2.4,0.9975) -- (2.6,0.9975);
				\draw[c2,thick] (2.6,0.9975) -- (2.7,1);
				\draw (2.23,0.7) node[scale=1.25]{$\tilde\Omega$} ;
				\draw[c2,thick] (2.6,1.1) -- (2.6,1.3) -- (2.4,1.3);
			\end{tikzpicture}
			\caption{Domain $\tilde \Omega$.}
		\end{center}
	\end{subfigure}
	\caption{Domain with one complex feature and illustration of the notation. \review{The different domains are distinguished by different gray intensities.}} \label{fig:extentionnotation}
\end{figure}
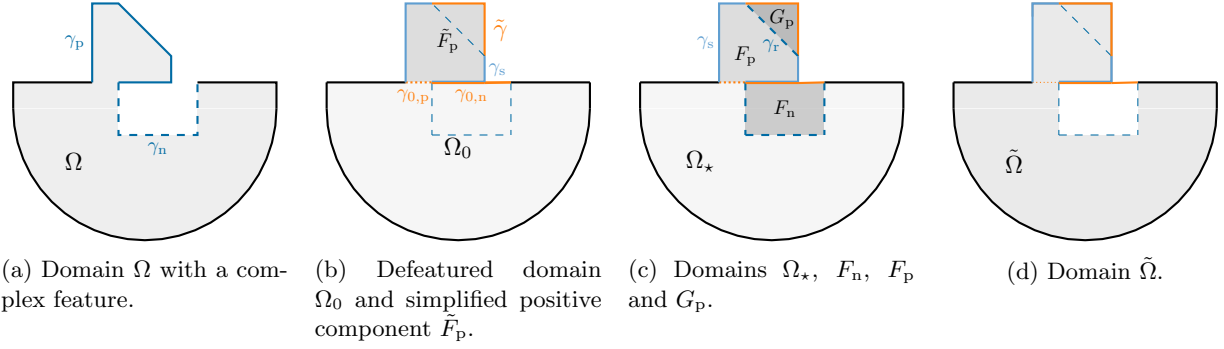

Now, we remark that the solution $u_0$ of $\mathcal P(\Omega_0)$ is not defined in $F_\mathrm p$, i.e., in the positive component of the features of $\Omega$, since $F_\mathrm p \not\subset \Omega_0$ but $F_\mathrm p \subset \Omega$. Thus, to be able to compare $u$ and $u_0$ in $F_\mathrm p$, we need to solve for each feature $F^k\in\mathfrak F$ a local differential problem $\mathcal P\big(F_\mathrm p^k\big)$ in $F_\mathrm p^k$, whose solution suitably extends $u_0$ to $F_\mathrm p^k$.
However, and as already noted in \cite{paper1defeaturing}, the domain $F_\mathrm p^k$ might be complicated or even non-smooth, thus finding the solution of the extension problem $\mathcal P\big(F_\mathrm p^k\big)$ might be cumbersome. To circumvent this difficulty, instead of solving $\mathcal P\big(F_\mathrm p^k\big)$, we solve an extension problem $\mathcal P\big(\tilde F_\mathrm p^k\big)$ in a simpler domain $\tilde F_\mathrm p^k$ containing $F_\mathrm p^k$, see Figure~\ref{fig:extentionnotation}. \review{Let us call $u_k$ the solution of $\mathcal P\big(\tilde F_\mathrm p^k\big)$.} The domain $\tilde F_\mathrm p^k$ can be for instance the bounding box of $F_\mathrm p^k$, if the \review{solution $u_k$ restricted} to $F_\mathrm p^k$ correctly defines an extension of $u_0$ in $F_\mathrm p^k$, see Section~\ref{ss:notationbdext}, Equation \eqref{eq:condextension}. Then we define the \emph{extended defeatured solution} $u_\mathrm d$ in $\Omega$ as
\begin{equation}\label{eqid:defudmulti}
	u_\mathrm d \equiv \begin{cases}
		u_0\vert_{\Omega_\star}  & \text{in } \Omega_\star                                    \\
		u_k\vert_{F_\mathrm p^k} & \text{in } F_\mathrm p^k, \text{ for all } k=1,\ldots,N_f, 
	\end{cases}
\end{equation}
where we recall the definition of $\Omega_\star$ from \eqref{eq:defOmegastar}. We are now finally able to precisely define the defeaturing error in the energy norm by
$\vertiii{u-u_\mathrm d}_\Omega$.

\subsection{Main results}

In this article, we will precisely define the considered differential problems $\mathcal P(\Omega)$, $\mathcal P(\Omega_0)$, and $\mathcal P\big(\tilde F_\mathrm p^k\big)$ for all $k=1,\ldots,N_f$ in the context of Poisson's, linear elasticity, and Stokes' equations. Note that when the linear elasticity equations are considered, then $u$, $u_0$, $u_k$, and $u_\mathrm d$ are vector-valued functions that we will also, respectively, write $\boldsymbol u$, $\boldsymbol u_0$, $\boldsymbol u_k$, and $\boldsymbol u_\mathrm d$. Similarly, for Stokes equations, $u$, $u_0$, $u_k$, and $u_\mathrm d$ also include the pressure, i.e., we will write them $(\boldsymbol u, p)$, $(\boldsymbol u_0, p_0)$, $(\boldsymbol u_k, p_k)$, and $(\boldsymbol u_\mathrm d, p_\mathrm d)$, respectively.

For each case, we will then define an \textit{a posteriori} estimator $\mathscr E(u_\mathrm d)$ of the defeaturing error in the energy norm (Sections \ref{s:poisson} \changes{and \ref{s:linelaststokes}}), in the case in which Neumann boundary conditions are imposed on the boundaries of the features. This assumption is formalized as follows:
\begin{assumption}\label{assu:notrimmingdirichlet}
	Let us decompose $\partial \Omega = \overline{\Gamma_D} \cup \overline{\Gamma_N}$ such that Dirichlet boundary conditions are imposed on $\Gamma_D$ with $|\Gamma_D| > 0$, Neumann boundary conditions are imposed on $\Gamma_N$, and $\Gamma_D \cap \Gamma_N = \emptyset$. Then we assume that for all $k=1,\ldots,N_f$,
	$$\Gamma_D \cap \big( \partial F_\mathrm n^k \cup \partial F_\mathrm p^k \big) = \emptyset.$$
\end{assumption}

\review{Then for} each one of the differential problems $\mathcal P(\Omega)$ that will be considered in this article, the introduced \textit{a posteriori} defeaturing error estimator $\mathscr E(u_\mathrm d)$ verifies the two following Theorems~\ref{thm:upperbound} and~\ref{thm:lowerbound}.
\begin{theorem}[Reliability]\label{thm:upperbound}
	\review{Under \review{Assumptions~\ref{as:separatednew}} and~\ref{assu:notrimmingdirichlet},} the defeaturing error estimator $\mathscr E(u_\mathrm d)$ is reliable, meaning that it is an upper bound for the defeaturing error in the energy norm:
	$$\vertiii{u-u_\mathrm d}_\Omega \lesssim \mathscr{E}(u_\mathrm d).$$
\end{theorem}

\review{Before stating Theorem~\ref{thm:lowerbound}, let us just introduce the following definition.}
\begin{definition} \label{defn:pwsmoothshaperegweak} Let $\Lambda$ be an $(n-1)$-dimensional \review{subset} of $\mathbb R^n$. We say that $\Lambda$ is \emph{regular} if $\Lambda$ is piecewise smooth and shape regular, that is, if \review{there exists $N_\Lambda\in \mathbb{N}\setminus\{0\}$ such that} for all $\ell_1,\ell_2=1,\ldots,N_\Lambda$ with $\ell_1\neq \ell_2$, $\Lambda = \mathrm{int}\left(\displaystyle\bigcup_{\ell=1}^{N_\Lambda} \overline{\Lambda^\ell}\right)$, $\Lambda^{\ell_1} \cap \Lambda^{\ell_2} = \emptyset$, $|\Lambda|\lesssim \left|\Lambda^{\ell_1}\right|$ \review{where the hidden constant may depend on $N_\Lambda$ but not on the measure of each component $\Lambda^{\ell_1}$}, and $\Lambda^{\ell_1}$ is smooth.
\end{definition}

\begin{theorem}[Efficiency] \label{thm:lowerbound}
	\review{Under \review{Assumptions~\ref{as:separatednew}} and~\ref{assu:notrimmingdirichlet},} if the boundaries of the features are shape regular as in Definition~\ref{defn:pwsmoothshaperegweak}, then the defeaturing error estimator $\mathscr E(u_\mathrm d)$ is efficient up to oscillations, meaning that it is a lower bound for the defeaturing error in the energy norm:
	$$\mathscr E(u_\mathrm d) \lesssim \vertiii{u-u_\mathrm d}_\Omega + \mathrm{osc}(u_\mathrm d),$$
	where $\mathrm{osc}(u_\mathrm d)$ is a higher order term with respect to the size of the features.
\end{theorem}
We remark again that the inequalities in Theorems~\ref{thm:upperbound} and~\ref{thm:lowerbound} do not depend neither on the size nor on the number of features.
In addition, as it will be shown in Sections \ref{s:poisson} \changes{and \ref{s:linelaststokes}} for Poisson, Stokes, and linear elasticity problems, the estimator $\mathscr E(u_\mathrm d)$ is simple and computationally cheap, and it is not only driven by geometrical considerations, but also by the PDE at hand.

\changesbis{The explicit expression of the oscillation term $\mathrm{osc}(u_\mathrm d)$ will be provided in Appendix \ref{app:proofs}, together with the proofs of Theorems~\ref{thm:upperbound} and~\ref{thm:lowerbound}, in the framework of Stokes' equations, as the oscillations corresponding to Poisson's and linear elasticity equations are very similar and their details can be found in \cite{phdthesis}.} The key issue in the analysis is to track the dependence of all constants from the sizes of the features and from their number.

\subsection{Some further geometric notation} \label{ss:notationbdext}
To be able to correctly define the defeaturing error estimator $\mathscr{E}(u_\mathrm d)$, we need to introduce some further notation identifying specific pieces of boundaries of the features.
In the following, we use the upper index $k$ to refer to the feature $F^k$, with $k=1,\ldots,N_f$, and the lower indices $\mathrm n$ and $\mathrm p$ to refer to negative and positive components of the features, respectively.

Now in particular, we let $\gamma$ be the union of the pieces of boundaries of $\Omega$ that are removed by defeaturing, and we let $\gamma_0$ be the union of the pieces of boundaries of $\Omega_0$ replacing them, that is, and as illustrated in Figure~\ref{fig:exfeatcomplex},
\begin{align}
	\gamma_0 := \bigcup_{k=1}^{N_f} \gamma_0^k                                                                                              & \quad \text{ and } \quad \gamma := \bigcup_{k=1}^{N_f} \gamma^k, \nonumber                                                                                                                                              \\
	\text{ where }\quad
	\gamma_0^k := \text{int}\left(\overline{\gamma_{0,\mathrm n}^k} \cup \overline{\gamma_{0,\mathrm p}^k}\right) \subset \partial \Omega_0 & \quad \text{with} \quad \gamma_{0,\mathrm n}^k := \partial F_\mathrm n^k \setminus \partial \Omega_\star, \quad \gamma_{0,\mathrm p}^k := \partial F_\mathrm p^k \setminus \partial \Omega, \nonumber                   \\
	\gamma^k :=\text{int}\left(\overline{\gamma_{\mathrm n}^k} \cup \overline{\gamma_{\mathrm p}^k}\right) \subset \partial \Omega          & \quad \text{with} \quad \gamma_\mathrm n^k:= \partial F_\mathrm n^k\setminus \overline{\gamma_{0,\mathrm n}^k}, \quad \gamma_\mathrm p^k := \partial F_\mathrm p^k\setminus \overline{\gamma_{0,\mathrm p}^k},\nonumber
\end{align}
so that $\partial F_\mathrm n^k = \overline{\gamma_\mathrm n^k} \cup \overline{\gamma_{0,\mathrm n}^k}$ with $\gamma_\mathrm n^k\cap \gamma_{0,\mathrm n}^k=\emptyset$, and $\partial F_\mathrm p^k = \overline{\gamma_\mathrm p^k} \cup \overline{\gamma_{0,\mathrm p}^k}$ with $\gamma_\mathrm p^k\cap \gamma_{0,\mathrm p}^k=\emptyset$.
Moreover, let
\begin{align*}
	\gamma_{0,\mathrm n} := \bigcup_{k=1}^{N_f} \gamma_{0,\mathrm n} ^k, \quad \gamma_{0,\mathrm p}  := \bigcup_{k=1}^{N_f} \gamma_{0,\mathrm p}^k,  \nonumber \\
	\gamma_\mathrm n := \bigcup_{k=1}^{N_f} \gamma_\mathrm n^k, \quad \gamma_\mathrm p := \bigcup_{k=1}^{N_f} \gamma_\mathrm p^k. \label{eq:defcupmulti}
\end{align*}

Using this notation and following the discussion in Section~\ref{ss:pbpres}, we can now precisely determine which simple extension $\tilde F_\mathrm p^k$ of the positive component $F_\mathrm p^k$ can be chosen, for all $k=1,\ldots,N_f$ (see Figure~\ref{fig:extentionnotation}). More precisely, for the solution of $\mathcal P\big(\tilde F_\mathrm p^k\big)$ restricted to $F_\mathrm p^k$ to correctly define an extension of the solution of $\mathcal P(\Omega_0)$ in $F_\mathrm p^k$, we need
\begin{equation}\label{eq:condextension}
	\tilde F_\mathrm p^k \supset F_\mathrm p^k, \qquad \gamma_{0,\mathrm p}^k\subset \left(\partial \tilde{F}_\mathrm p^k \cap \partial F_\mathrm p^k\right).
\end{equation}
Note that it is possible to have $\tilde F_\mathrm p^k\cap \Omega_0 \neq \emptyset$.
Let us also define
$$G_\mathrm p^k := \tilde F_\mathrm p^k \setminus \overline{F_\mathrm p^k} \quad \text{ for all } k=1,\ldots,N_f \qquad \text{ and } \qquad G_\mathrm p := \bigcup_{k=1}^{N_f} G_\mathrm p^k.$$
We remark that for all $k=1,\ldots,N_f$, one can look at $\tilde F_\mathrm p^k $ as the defeatured geometry of the positive component $F_\mathrm p^k $, that is, as a geometry simplified from the exact geometry $F_\mathrm p^k $, for which $G_\mathrm p^k $ is a negative feature  (see Figure~\ref{fig:extentionnotation_Gp}).
To simplify the following exposition, and even if this hypothesis could easily be removed, let us make the following assumption regarding these domain extensions.
\begin{assumption} \label{asid:simplassumptionfeat}
	Let us assume that
	\begin{itemize}
		\item $\tilde F_\mathrm p^k \cap \tilde F_\mathrm p^\ell = \emptyset$ for all $k,\ell = 1,\ldots,N_f$ such that $k\neq \ell$,
		\item if we let $\tilde F_\mathrm p := \displaystyle\bigcup_{k=1}^{N_f} \tilde F_\mathrm p^k$, then $\tilde F_\mathrm p\cap \Omega_\star = \emptyset$ for all $k = 1,\ldots,N_f$.
	\end{itemize}
\end{assumption}

In addition, let $\tilde \gamma^k := \partial \tilde F_\mathrm p^k \setminus \partial F_\mathrm p^k$, and let $\gamma_\mathrm p^k$ be decomposed as $\gamma_\mathrm p^k = \text{int}\big(\overline{\gamma_\intersign^k}\cup\overline{\gamma_\setminussign^k}\big)$, being $\gamma_\intersign^k$ and $\gamma_\setminussign^k$ open, $\gamma_\intersign^k$ is the part of $\gamma_\mathrm p^k$, that is shared with $\partial \tilde F_\mathrm p^k$, while $\gamma_\setminussign^k$ is the remaining part of $\gamma_\mathrm p^k$, that is, the part that does not belong to $\partial \tilde F_\mathrm p^k$.
I.e., $\gamma_\intersign^k=\gamma_\mathrm p^k\cap\partial \tilde F_\mathrm p^k$ and $\gamma_\setminussign^k=\gamma_\mathrm p^k\setminus\gamma_\intersign^k$.
This notation is illustrated in Figure~\ref{fig:extentionnotation}.
Then, similarly to the previously introduced notation, let
\begin{equation*}
	\tilde \gamma := \bigcup_{k=1}^{N_f} \tilde \gamma^k, \quad \gamma_\intersign := \bigcup_{k=1}^{N_f} \gamma_\intersign^k, \quad \gamma_\setminussign := \bigcup_{k=1}^{N_f} \gamma_\setminussign^k.
\end{equation*}

In the sequel, we will see that the boundaries $\gamma_\mathrm n$, $\gamma_{0,\mathrm p}$, and $\gamma_\setminussign$ will play an important role in the definition of the defeaturing error estimators $\mathscr{E}(u_\mathrm d)$. Therefore, let us also introduce the following notation:
\begin{align}
	\Gamma^k & := \gamma_\mathrm n^k\cup \gamma_{0,\mathrm p}^k\cup \gamma_\setminussign^k, \quad \text{ for } k=1,\ldots,N_f, \label{eq:gammamulti}
\end{align}
and the sets
\begin{align}
	\Sigma_\mathrm n & := \left\{ \gamma_\mathrm n^k \right\}_{k=1}^{N_f}, \quad \Sigma_{0,\mathrm p} := \left\{\gamma^k_{0,\mathrm p}\right\}_{k=1}^{N_f}, \quad \Sigma_\setminussign := \left\{\gamma^k_\setminussign\right\}_{k=1}^{N_f}, \nonumber \\
	\Sigma^k         & := \left\{\gamma_\mathrm n^k, \gamma_{0,\mathrm p}^k, \gamma_\setminussign^k\right\}, \quad \text{ for } k=1,\ldots,N_f, \nonumber                                                                                              \\
	\Sigma           & := \{\boundarypiece \in \Sigma^k : k=1,\ldots,N_f\}. \label{eq:sigmamulti}
\end{align}

Finally, let $\mathbf n$ and $\mathbf n_0$ respectively denote the unit outward normal vectors to $\Omega$ and $\Omega_0$. Moreover, for all $k=1,\ldots,N_f$, let $\tilde{\mathbf{n}}^k$ be the unitary outward normal to $\tilde F_\mathrm p^k$, and let $\mathbf n^k\equiv \mathbf n_{F^k}$ denote the unitary outward normal vector to $F_\mathrm n^k$ and to $F_\mathrm p^k$. Note that the vectors $\mathbf n^k$ may not be uniquely defined if the outward normal to $F_\mathrm n^k$ is of opposite sign of the outward normal to $F_\mathrm p^k$ \changes{(see Figure~\ref{fig:extentionnotation} for instance)}, but we allow this abuse of notation since the context will always make it clear.\\

\begin{figure}
	\centering
	\begin{center}
		\begin{tikzpicture}[scale=5]
		\draw (-0.3,2) -- (-0.2,2);
		\draw (-0.2,2) -- (-0.2,2.1);
		\draw (-0.2,2.1) -- (-0.10,2.1);
		\draw (-0.1,2.1) -- (-0.1,2);
		\draw (-0.1,2) -- (0.6,2);
		\draw (0.6,2) -- (0.6,2.1);
		\draw (0.6,2.1) -- (0.7,2.1);
		\draw (0.7,2.1) -- (0.7,2);
		\draw (0.7,2) -- (0.65,2);
		\draw (0.65,2) -- (0.65, 1.9);
		\draw (0.65, 1.9) -- (0.8,1.9);
		\draw (0.8,1.9) -- (0.8,2.1);
		\draw (0.8,2.1) -- (0.9,2.1);
		\draw (0.9, 2.1) -- (0.9,2);
		\draw (0.9,2) -- (1,2);
		\draw (1,2) -- (1,1.5);
		\draw (1,1.5) -- (0.4,1.5);
		\draw (0.4,1.5) -- (0.4,1.6);
		\draw (0.4,1.6) -- (0.3,1.6);
		\draw (0.3,1.6) -- (0.3,1.5);
		\draw (0.3,1.5) -- (-0.3,1.5);
		\draw (-0.3,1.5) -- (-0.3,2);
		\fill[c2, fill, opacity=0.3] (-0.2,2) rectangle (-0.1,2.1);
		\draw[c2, dashed, thick] (-0.2,2) rectangle (-0.1,2.1);
		\fill[c2, fill, opacity=0.3] (-0.3,1.6) rectangle (0,2);
		\draw[c2] (-0.15,1.85) node{$\Omega^1$} ;
		\fill[c1, fill, opacity=0.3] (0.2,1.6) rectangle (0.5,2);
		\fill[c1, fill, opacity=0.3] (0.2,1.5) rectangle (0.3,1.6);
		\fill[c1, fill, opacity=0.3] (0.4,1.5) rectangle (0.5,1.6);
		\draw[c1, dashed, thick] (0.3,1.5) rectangle (0.4,1.6);
		\draw[c1] (0.35,1.75) node{$\Omega^2$} ;
		\fill[c3, fill, opacity=0.3] (0.5,1.65) rectangle (1,1.9);
		\fill[c3, fill, opacity=0.3] (0.5,1.9) rectangle (0.65,2);
		\fill[c3, fill, opacity=0.3] (0.6,2) rectangle (0.7,2.1);
		\fill[c3, fill, opacity=0.3] (0.8,1.9) rectangle (1,2);
		\fill[c3, fill, opacity=0.3] (0.8,2) rectangle (0.9,2.1);
		\draw[c3, dashed, thick] (0.6,2)--(0.65, 2) -- (0.65, 1.9) -- (0.8, 1.9) -- (0.8, 2) -- (0.9,2) -- (0.9,2.1) -- (0.8,2.1) -- (0.8,2) -- (0.7,2) -- (0.7,2.1) -- (0.6,2.1) -- (0.6,2);
		\draw[c] (0.75,1.8) node{$\Omega^3$} ;
		\draw (-0.15,2.05) node{$F^1$};
		\draw (0.35,1.55) node{$F^2$};
		\draw (0.725,1.95) node{$F_\mathrm n^3$};
		\draw (0.75,2.15) node{$F_\mathrm p^3$};
		\draw[->] (0.73,2.13) -- (0.65,2.05);
		\draw[->] (0.79,2.13) -- (0.85,2.05);
		\end{tikzpicture}
		\caption{\changesbis{Domain $\Omega$ with three separated features, and a possible choice of subdomains $\Omega^k$, $k=1,2,3$, satisfying Assumption~\ref{as:repeatseparatednew}.}} \label{fig:Omegak}
	\end{center}
\end{figure}
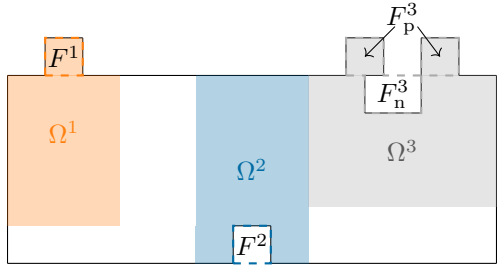

\review{With this technical notation at hand, we can now precisely restate condition~\ref{it:conditionbseparated} of Assumption~\ref{as:separatednew}.
\begin{repeatassumption}{as:separatednew}[b] \label{as:repeatseparatednew}
	There exist sub-domains $\Omega^k\subset \Omega$, $k=1,\ldots,N_f$ such that
	\begin{itemize}
		\item $F_\mathrm p^k \subset \Omega^k,\quad \left(\gamma_\mathrm n^k \cup \gamma_\setminussign^k\right) \subset \partial \Omega^k, \quad \gamma_{0,\mathrm p}^k \subset \partial (\Omega^k \cap \Omega_0)$,
		\item $\left|\Omega^k\right| \simeq |\Omega|$ where the hidden constant is independent of the size of the features, i.e., the measure of $\Omega^k$ is comparable with the measure of $\Omega$, not with the measure of the feature $F^k$,
		\item $N_s := \displaystyle\max_{J\subset \left\{1,\ldots,N_f\right\}} \left(\#J : \bigcap_{k\in J}\Omega^k \neq \emptyset \right) \ll N_f$, that is, the maximum number $N_s$ of superposed sub-domains $\Omega^k$ is limited and notably smaller than the total number of features $N_f$.
	\end{itemize}
\end{repeatassumption}
This condition is illustrated in Figure~\ref{fig:Omegak}. Note that if $N_f = 1$, one can take $\Omega^1 := \Omega$.}

\section{Defeaturing in Poisson's equation}\label{s:poisson}
In this section, we consider the Poisson equation and precisely define the exact problem $\mathcal P(\Omega)$, the defeatured problem $\mathcal P(\Omega_0)$, and the extension problems $\mathcal P\big(\tilde F_\mathrm p^k\big)$ for all $k=1,\ldots,N_f$. Then in this context, we give a precise definition of the proposed reliable and efficient \textit{a posteriori} estimator $\mathscr E(u_\mathrm d)$ of the energy norm of the defeaturing error $\vertiii{u-u_\mathrm d}_\Omega$.

In the following, we denote by $H^s(D)$ the Sobolev space of order $s\in\mathbb R$ in a domain $D\subset \mathbb R^n$, whose classical norm and semi-norm are written $\|\cdot\|_{s,D}$ and $|\cdot|_{s,D}$, respectively. We will also denote by $L^2(D):=H^0(D)$. Moreover, if we let $\varphi \subset \partial D$ and $z\in H^\frac{1}{2}(\varphi)$, and to be able to deal with boundary conditions, we will denote by
$$H_{z,\varphi}^1(D) := \{ y\in H^1(D): \mathrm{tr}_\varphi(y)=z\},$$
where $\mathrm{tr}_\varphi(y)$ denotes the trace of $y$ on $\varphi$.

\subsection{Exact and defeatured problems}
Let us first introduce Poisson's problem $\mathcal P(\Omega)$ in the exact geometry $\Omega$. To do so, let $g_D\in H^{\frac{3}{2}}(\Gamma_D)$, $g\in H^{\frac{1}{2}}(\Gamma_N)$ and $f\in L^2\left(\Omega\right)$. \review{Note that one order of regularity more than usual is required. As we will further see, to define the proposed defeaturing error estimator, we need a solution whose normal derivative along the features’ boundaries belongs to $L^2$.}
Then, the problem reads:
find \changes{$u\in H^1_{g_D,\Gamma_D}(\Omega)$}, the weak solution of
\begin{align}
	\mathcal P(\Omega): & \begin{cases}
		-\Delta u = f                                           & \text{ in } \Omega                  \\
		u = g_D                                                 & \text{ on } \Gamma_D                \\
		\displaystyle\frac{\partial u}{\partial \mathbf{n}} = g & \text{ on } \Gamma_N,\vspace{0.1cm}
	\end{cases}
	\quad \changes{\text{i.e.,}\quad \forall} v\in H^1_{0,\Gamma_D}(\Omega),\quad
	\mathfrak a(u,v) := \int_\Omega \nabla u \cdot \nabla v \,\mathrm dx = \int_\Omega fv \,\mathrm dx + \int_{\Gamma_N} g v \,\mathrm ds.
	\label{eqid:weakoriginalpb}
\end{align}
If $H_{0,\Gamma_D}^1(\Omega)$ is equipped with the $L^2(\Omega)$-norm of the gradient $\|\nabla \cdot\|_{0,\Omega}$, then we know from \review{Riesz representation} theorem that this problem is well-posed. However, note that this problem is usually not solved in practice, as it is assumed to be computationally expensive.

Let us therefore introduce the corresponding Poisson problem $\mathcal P(\Omega_0)$ in the defeatured geometry $\Omega_0$. To do so, we need to consider an $L^2$-extension of the restriction $f\vert_{\Omega_\star}$ in all the negative components $F_\mathrm n^k$, $k=1,\ldots,N_f$, that we still write $f\in L^2(\Omega_0)$ by abuse of notation. Note that such an extension is not needed in the positive components of the features.
Then instead of~\eqref{eqid:weakoriginalpb}, we solve the following defeatured Poisson problem:
after choosing $g_{0}\in H^{\frac{1}{2}}(\gamma_0)$, find \changes{$u_0\in H^1_{g_D,\Gamma_D}(\Omega_0)$}, the weak solution of
\begin{align} 
	\mathcal P(\Omega_0): \begin{cases}
		-\Delta u_0 = f                                                  & \text{ in } \Omega_0                   \\
		u_0 = g_D                                                        & \text{ on } \Gamma_D                   \\ \vspace{1mm}
		\displaystyle\frac{\partial u_0}{\partial \mathbf{n}_0}  = g     & \text{ on } \Gamma_N\setminus {\gamma} \\
		\displaystyle\frac{\partial u_0}{\partial \mathbf{n}_0}  = g_{0} & \text{ on } \gamma_0,\vspace{0.1cm}
	\end{cases}
	\quad \changes{\text{i.e.,}\quad}
	\int_{\Omega_0} \nabla u_0 \cdot \nabla v_0 \,\mathrm dx = \int_{\Omega_0} fv_0 \,\mathrm dx + \int_{\Gamma_N\setminus \gamma} g v_0 \,\mathrm ds + \int_{\gamma_0} g_{0} v_0 \,\mathrm ds\label{eqid:weaksimplpb}
\end{align}
for all $v_0\in H^1_{0,\Gamma_D}(\Omega_0)$. Let us recall that, according to Assumption~\ref{assu:notrimmingdirichlet}, $\Gamma_D\cap\gamma=\emptyset$.
From the \review{Riesz representation} theorem, we know that problem~\eqref{eqid:weaksimplpb} is well-posed. 

\begin{remark}
	The \review{best possible choices for} the defeatured problem data $f$ in $F_\mathrm n^k$ for $k=1,\ldots,N_f$ and $g_0$ on $\gamma_0$ will be guided by Remark~\ref{rmk:compatcondmultipoisson}. \review{However, remark that in applications, the defeatured problem data are rarely chosen, being rather determined by the problem at hand, and being often taken as the natural extension of the exact problem data. For instance, if the right hand side $f$ corresponds to the gravity field, then its only possible physical extension to the negative components of the features is still the gravity field itself. Moreover, one of the aims of defeaturing relies in avoiding the local meshing of the features when the differential problem is solved numerically. Therefore, when the geometry is defeatured, one does not want to have to do a local meshing to capture the behavior of the defeaturing problem data, as this would be equivalent to meshing the feature.} We only anticipate here that the best possible \review{defeatured problem data} allow the conservation of the solution flux in the positive and negative components of the features. That is, \review{the defeaturing error will be smaller if the defeatured problem data verify} for all $k=1,\ldots,N_f$, 
	\begin{align}
	\int_{\gamma_{0,\mathrm p}^{k}} g_{0} \,\mathrm ds = \int_{\gamma_\mathrm p^k} g \,\mathrm{d}s + \int_{F_\mathrm p^k} f\,\mathrm dx, \qquad
	\text{and } \qquad \int_{\gamma_{0,\mathrm n}^{k}} g_{0} \,\mathrm ds = \int_{\gamma_\mathrm n^{k}} g \,\mathrm{d}s - \int_{F_\mathrm n^k} f\,\mathrm dx. \label{eq:compatcondmultipoisson0}
	\end{align}
\end{remark}

Now, we need to extend the solution $u_0$ of~\eqref{eqid:weaksimplpb} to $\tilde F^k_\mathrm p$, for $k=1,\ldots,N_f$, as discussed in Section~\ref{ss:pbpres}, where $\tilde F^k_\mathrm p$ satisfies the properties given in~\eqref{eq:condextension}.
To do so, we need to consider an $L^2$-extension of the restriction $f\vert_{F_\mathrm p^k}$ in each $G_\mathrm p^k := \tilde F_\mathrm p^k\setminus \overline{F_\mathrm p^k}$, that we still write $f$ by abuse of notation. Then, we solve the following extension problem for $k=1,\ldots,N_f$: after choosing $\tilde g^k\in H^\frac{1}{2}(\tilde \gamma^k)$, find \changes{$u_k\in H^1_{u_0,\gamma_{0,\mathrm p}^k}\big(\tilde F_\mathrm p^k\big)$}, the weak solution of
\begin{align} 
	\mathcal P\big(\tilde F_\mathrm p^k\big): \quad\begin{cases}
		-\Delta u_k = f                                                               & \text{ in } \tilde{F}_\mathrm p^k            \\ \vspace{1mm}
		u_k = u_0                                                                     & \text{ on } \gamma_{0,\mathrm p}^k           \\ \vspace{1mm}
		\displaystyle\frac{\partial u_k}{\partial \tilde{\mathbf{n}}^k}  = \tilde g^k & \text{ on } \tilde \gamma^k                  \\
		\displaystyle\frac{\partial u_k}{\partial \tilde{\mathbf{n}}^k}  = g          & \text{ on } \gamma_\intersign^k,\vspace{1mm}
	\end{cases}
	\quad \changes{\text{i.e.,}\quad}
	\int_{\tilde{F}_\mathrm p^k} \nabla u_k \cdot \nabla v^k \,\mathrm dx = \int_{\tilde{F}_\mathrm p^k} fv^k \,\mathrm dx + \int_{\tilde \gamma^k} \tilde{g}^k v^k \,\mathrm ds + \int_{\gamma_\intersign^k} {g} v^k \,\mathrm ds \label{eqid:weakfeaturepbmulti}
\end{align}
for all $v^k\in H^1_{0,\gamma_{0,\mathrm p}^k}\big(\tilde{F}_\mathrm p^k\big)$. As before, from \review{Riesz representation} theorem, we know that problem~\eqref{eqid:weakfeaturepbmulti} is well-posed. Once again, the choice of the defeatured problem data $f$ in $G_\mathrm p^k$ and $\tilde g^k$ on $\tilde \gamma^k$, for $k=1,\ldots,N_f$, will be guided by Remark~\ref{rmk:compatcondmultipoisson}.
Thus, as in~\eqref{eqid:weaksimplpb}, the best possible choices will allow for the conservation of the solution flux in every feature extension $G_\mathrm p^k$. That is, the best possible choices verify for all $k=1,\ldots,N_f$,
\begin{equation}
	\int_{\tilde \gamma^k} \tilde g^k \,\mathrm ds = \int_{\gamma^k_{\setminussign}} g \,\mathrm{d}s - \int_{G_\mathrm p^k} f\,\mathrm dx. \label{eq:compatcondmultipoissontilde}
\end{equation}

Then, the defeatured solution $u_\mathrm d\in H_{g_D,\Gamma_D}(\Omega)$ is defined from $u_0$ and $u_k$ for $k=1,\ldots,N_f$ as in~\eqref{eqid:defudmulti}, and the energy norm of the defeaturing error is defined as
\begin{equation}\label{eq:defeaterrorpoisson}
	\vertiii{u-u_\mathrm d}_\Omega := \big( \mathfrak a(u-u_\mathrm d,u-u_\mathrm d)\big)^\frac{1}{2} = \|\nabla (u-u_\mathrm d)\|_{0,\Omega}.
\end{equation}

\subsection{Defeaturing error estimator}
In this section, we generalize the defeaturing error estimator introduced in \cite{paper1defeaturing}, in the case of a geometry that presents multiple features.
As for the single feature case, the derived estimator is an upper bound and a lower bound (up to oscillations) of the energy norm of the defeaturing error.

Let us recall the definition of the defeatured solution $u_\mathrm d$ from~\eqref{eqid:defudmulti}, and the definitions of $\Sigma$, $\Sigma^k$, $\Sigma_\mathrm n$, $\Sigma_{0,\mathrm p}$, and $\Sigma_\setminussign$ from~\eqref{eq:sigmamulti}. Then for all $\boundarypiece \in \Sigma$, let $k_\boundarypiece \equiv k$ if $\boundarypiece \in \Sigma^k$ for some $k=1,\ldots,N_f$, and let
\begin{equation} \label{eqmf:dsigma}
	d_\boundarypiece \equiv \begin{cases}\vspace{1mm}
		g - \displaystyle\frac{\partial u_\mathrm d}{\partial \mathbf n}                                      & \text{if } \boundarypiece \in \Sigma_\mathrm n \text{ or if } \boundarypiece \in \Sigma_\setminussign, \\\vspace{1mm}
		- \left( g_0 + \displaystyle\frac{\partial u_\mathrm d}{\partial \mathbf n^{k_\boundarypiece}}\right) & \text{if } \boundarypiece \in \Sigma_{0,\mathrm p}.
	\end{cases}
\end{equation}
In other words, $d_\boundarypiece$ represents the Neumann error on the boundaries $\gamma_\mathrm n$ and $\gamma_\setminussign$ due to the defeaturing process, and the jump of the normal derivative of the defeatured solution on the boundaries $\gamma_{0,\mathrm p}$.
Then, if we let $\eta\in\mathbb{R}$ be the unique solution of $\eta = -\log(\eta)$, and if for all $\boundarypiece\in\Sigma$, we let
\begin{align} \label{eqp1:negconstant}
	c_{\boundarypiece} := \begin{cases}
		\max\big(\hspace{-0.05cm}-\log\left(|\boundarypiece|\right), \eta \big)^\frac{1}{2} & \text{ if } n = 2  \\
		1                                                                                   & \text{ if } n = 3,
	\end{cases}
\end{align}
we define the \textit{a posteriori} defeaturing error estimator as:
\begin{equation}\label{eq:multiestimatorpoisson}
	\mathscr{E}( u_\mathrm d) := \left(\sum_{\boundarypiece\in\Sigma} \mathscr{E}_\boundarypiece({u}_\mathrm d)^2 \right)^\frac{1}{2},
\end{equation}
where, for all $\boundarypiece\in\Sigma$,
\begin{equation*}
	\mathscr{E}_\boundarypiece(u_\mathrm d) := \left(|\boundarypiece|^\frac{1}{n-1} \left\|d_\boundarypiece - \overline{d_\boundarypiece}^\boundarypiece\right\|_{0,\boundarypiece}^2 + c_\boundarypiece^2 |\boundarypiece|^\frac{n}{n-1} \left|\overline{ d_\boundarypiece}^\boundarypiece\right|^2\right)^\frac{1}{2},
\end{equation*}
where $\overline{d_\boundarypiece}^\boundarypiece$ denotes the average value of $d_\boundarypiece$ over $\boundarypiece$. 
Note that we can rewrite the estimator feature-wise as follows:
\begin{equation} \label{eq:estforeachfeatpoisson}
	\mathscr{E}(u_\mathrm d) = \left(\sum_{k=1}^{N_f} \sum_{\boundarypiece\in\Sigma^k} \mathscr{E}_{\boundarypiece}({u}_\mathrm d)^2 \right)^\frac{1}{2} = \left(\sum_{k=1}^{N_f} \mathscr{E}^{k}(u_\mathrm d)^2\right)^\frac{1}{2},
\end{equation}
where for all $k=1,\ldots,N_f$, we define $\mathscr{E}^{k}(u_k)$ as the defeaturing error estimator for feature $F^k$, that is,
\begin{equation*}
	\mathscr{E}^k(u_\mathrm d) := \left(\sum_{\boundarypiece\in\Sigma^k} \mathscr{E}_{\boundarypiece}({u}_\mathrm d)^2 \right)^\frac{1}{2}.
\end{equation*}

The proposed estimator indicates that the whole information on the error introduced by defeaturing multiple features, in the energy norm, is encoded in the boundary of the features, and can be accounted by suitably evaluating the error made on the normal derivative of the solution. This result generalizes the one from \cite{paper1defeaturing}.

\begin{remark} \label{rmk:compatcondmultipoisson}
	As analogously noted in \cite[Remarks~4.1 and~5.3]{paper1defeaturing}, the terms involving the average values of $d_\boundarypiece$ in the estimator $\mathscr E(u_\mathrm d)$ only depend on the defeatured problem data since for all $k=1,\ldots,N_f$,
	\begin{align*}
		\overline{d_{\gamma_\mathrm n^k}}^{\gamma_\mathrm n^k}         & = \overline{\left(g-\frac{\partial u_\mathrm d}{\partial \mathbf n}\right)}^{\gamma_\mathrm n^k} = \frac{1}{\left|{\gamma_\mathrm n^k}\right|}\left( \int_{\gamma_\mathrm n^k} g\,\mathrm ds - \int_{\gamma_{0,\mathrm n}^k} g_0\,\mathrm ds - \int_{F_\mathrm n} f\,\mathrm dx\right),                                           \\
		\overline{d_{\gamma_{0,\mathrm p}^k}}^{\gamma_{0,\mathrm p}^k} & = \overline{\left(g_0+\frac{\partial u_\mathrm d}{\partial \mathbf n^k}\right)}^{\gamma_{0,\mathrm p}^k} = \frac{1}{\left|\gamma_{0,\mathrm p}^k\right|}\left( \int_{\gamma_{0,\mathrm p}^k} g_0\,\mathrm ds - \int_{\gamma_\mathrm p^k} g\,\mathrm ds - \int_{F_\mathrm p^k} f\,\mathrm dx\right), \nonumber                     \\
		\overline{d_{\gamma_\setminussign^k}}^{\gamma_\setminussign^k} & = \overline{\left(g-\frac{\partial u_\mathrm d}{\partial \mathbf n}\right)}^{\gamma_\setminussign^k} = \frac{1}{\left|\gamma_\setminussign^k\right|}\left( \int_{\gamma_\setminussign^k} g\,\mathrm ds - \int_{\tilde \gamma^k} \tilde g^k\,\mathrm ds - \int_{\tilde F_\mathrm p^k\setminus F_\mathrm p^k} f\,\mathrm dx\right).
	\end{align*}
	As a consequence, if these terms dominate, this means that the defeatured problem data should be more accurately chosen, namely $g_0$, $\tilde g^k$, and the extension of $f$ to $G_\mathrm p^k$.
	Moreover, under the reasonable flux conservation assumptions~\eqref{eq:compatcondmultipoisson0} and~\eqref{eq:compatcondmultipoissontilde},
	the defeaturing error estimator~\eqref{eq:multiestimatorpoisson} rewrites
	$$\mathscr E(u_\mathrm d) := \left( \displaystyle\sum_{\boundarypiece\in\Sigma} \left|\boundarypiece\right|^{\frac{1}{n-1}} \left\| d_\boundarypiece \right\|^2_{0,\boundarypiece}\right)^\frac{1}{2}.$$
\end{remark}

\begin{remark} \label{rmk:estmultitildepoisson}
	Note that
	\begin{align*}
		\mathscr{E}(u_\mathrm d) \lesssim \left( \displaystyle\sum_{\boundarypiece\in\Sigma} c_\boundarypiece^2 \left|\boundarypiece\right|^{\frac{1}{n-1}} \left\| d_\boundarypiece \right\|^2_{0,\boundarypiece}\right)^\frac{1}{2} =: \tilde{\mathscr{E}}(u_\mathrm d).
	\end{align*}
	However, when $n=2$, and under the flux conservation conditions~\eqref{eq:compatcondmultipoisson0} and~\eqref{eq:compatcondmultipoissontilde}, $\tilde{\mathscr{E}}\left(u_\mathrm d\right)$ is sub-optimal since in this case, $\tilde{\mathscr{E}}(u_\mathrm d) \lesssim \displaystyle\max_{\boundarypiece\in\Sigma}\left( c_{\boundarypiece} \right) {\mathscr{E}}(u_\mathrm d)$. Indeed, no lower bound can be provided for $\tilde{\mathscr{E}}(u_\mathrm d)$.
\end{remark}

%

\section{\changes{Defeaturing in linear elasticity and in Stokes' equations}}\label{s:linelaststokes}
In this section, we consider the linear elasticity \changes{and the Stokes equations}, and, following the same structure as for the Poisson's problem in Section~\ref{s:poisson}, we define precisely the exact $\mathcal P(\Omega)$, defeatured $\mathcal P(\Omega_0)$, and extension problems $\mathcal P\big(\tilde F_\mathrm p^k\big)$ for all $k=1,\ldots,N_f$. Then in this context, we give a precise definition of the proposed \textit{a posteriori} estimator of the energy norm of the defeaturing error, 
\changesbis{together with the proof of its reliability and efficiency (up to oscillations)}. 

In the following, we denote by $\boldsymbol H^s(D):= \left[ H^s(D) \right]^n$ the vector-valued Sobolev space of order $s\in\mathbb R$ in a domain $D\subset \mathbb R^n$, whose classical norm and semi-norm are again written $\|\cdot\|_{s,D}$ and $|\cdot|_{s,D}$, respectively. We will also denote by $\boldsymbol L^2(D):=\boldsymbol H^0(D)$. Moreover, if we let $\varphi \subset \partial D$ and $\boldsymbol z\in \boldsymbol H^\frac{1}{2}(\varphi)$, and to be able to deal with boundary conditions, we will denote by
$$\boldsymbol H_{\boldsymbol z,\varphi}^1(D) := \{ \boldsymbol y\in \boldsymbol H^1(D): \mathrm{tr}_\varphi(\boldsymbol y)=\boldsymbol z\},$$
where $\mathrm{tr}_\varphi(\boldsymbol y)$ denotes the trace of $\boldsymbol y$ on $\varphi$.

\subsection{Exact and defeatured problems}
Let us first introduce the Stokes problem $\mathcal P(\Omega)$ in the exact geometry $\Omega$. To do so, considering a function $\boldsymbol v:\Omega\to\mathbb R^n$, let \review{$\tensorepsilon(\boldsymbol v):=\frac{1}{2}\left( \boldsymbol \nabla \boldsymbol v+ \boldsymbol \nabla \boldsymbol v^T\right)$ be} the linearized strain rate tensor in $\Omega$ and let \review{
\begin{equation}\label{eq:constitutiverelationstokes2}
\boldsymbol \sigma(\boldsymbol v) = 2\mu \boldsymbol{\varepsilon}(v).
\end{equation}} 
be the viscous stress tensor of the considered Newtonian fluid, where the \review{constant $\mu > 0$ 
is the dynamic 
viscosity.} 
Note that $\tensorsigma(\boldsymbol v)$ is the viscous stress tensor and not the total Cauchy stress tensor that would be defined by
$\hat{\tensorsigma}(\boldsymbol v, q) := \tensorsigma(\boldsymbol{v}) - q\,\pmb{\mathbb I}_n$ for some function $q:\Omega\to\mathbb R$ in the space of pressures.
\begin{assumption} \label{asid:lame_constant}
	Hereinafter, for the sake of simplicity in the exposition, we assume that $\lambda$ 
	\review{is} constant everywhere, and \review{it is therefore} naturally extended to the defeatured geometry $\Omega_0$.
\end{assumption}

Now, let $\boldsymbol g_D\in \boldsymbol H^{\frac{3}{2}}(\Gamma_D)$, $\boldsymbol g\in \boldsymbol H^{\frac{1}{2}}(\Gamma_N)$, $f_c\in L^2(\Omega)$, and $\boldsymbol f\in \boldsymbol L^2\left(\Omega\right)$. We are interested in the following Stokes problem defined in the exact geometry $\Omega$: find \changes{$(\boldsymbol u, p)\in\boldsymbol H^1_{\boldsymbol g_D, \Gamma_D}(\Omega)\times L^2(\Omega)$}, the weak solution of
\begin{align}
	\mathcal P(\Omega): \quad
	\begin{cases}
		-\boldsymbol \nabla \cdot \tensorsigma (\boldsymbol u) + \nabla p =\boldsymbol f & \text{ in } \Omega    \\
		\boldsymbol \nabla \cdot \boldsymbol u = f_c                                     & \text{ in } \Omega    \\
		\boldsymbol u = \boldsymbol g_D                                                  & \text{ on } \Gamma_D  \\
		\tensorsigma (\boldsymbol u)\mathbf{n} - p\mathbf n= \boldsymbol g               & \text{ on } \Gamma_N,
	\end{cases}
	\quad \changes{\text{i.e.,}} \quad
	\begin{cases}
		\mathfrak a(\boldsymbol{u}, \boldsymbol{v}) + \mathfrak b(p,\boldsymbol{v}) = \displaystyle\int_\Omega \boldsymbol f \cdot \boldsymbol v \,\mathrm dx + \displaystyle\int_{\Gamma_N} \boldsymbol g \cdot \boldsymbol v \,\mathrm ds \\
		\mathfrak b(q,\boldsymbol{u})                                               = -\displaystyle\int_{\Omega} q f_c\,\mathrm dx\label{eqid:weakoriginalstokespb}
	\end{cases}
\end{align}
for all $(\boldsymbol v,q)\in \boldsymbol H^1_{\boldsymbol 0,\Gamma_D}(\Omega)\times L^2(\Omega)$, where for all $\boldsymbol{v}, \boldsymbol{w}\in \boldsymbol{H}^1(\Omega)$ and all $q\in L^2(\Omega)$,
\begin{align}
	\mathfrak a(\boldsymbol{w}, \boldsymbol{v}) & := \displaystyle\int_\Omega \tensorsigma (\boldsymbol u) : \tensorepsilon(\boldsymbol v) \,\mathrm dx \quad \text{ and } \quad \mathfrak b(q,\boldsymbol{v}) := - \int_\Omega p\boldsymbol \nabla \cdot \boldsymbol{v} \,\mathrm dx. \label{eq:defabstokes}
\end{align}
If we equip $\boldsymbol H^1_{\boldsymbol 0, \Gamma_D}(\Omega)$ with the norm $\| \boldsymbol \nabla \boldsymbol \cdot \|_{0,\Omega}$, it is possible to show by Ladyzhenskaya-Babu\v{s}ka-Brezzi theorem~\cite{boffi2013mixed} that problem~\eqref{eqid:weakoriginalstokespb} is well-posed. \changes{However, note that this problem is never solved in practice, as it is assumed to be computationally too expensive.}

\review{\begin{remark}
	In a more general setting, the viscous stress tensor writes 
	\begin{equation}\label{eq:viscousstressgeneral}
	\tensorsigma(\boldsymbol v) = 2\mu\,\tensorepsilon(\boldsymbol v) + \lambda (\boldsymbol \nabla \cdot \boldsymbol v)\pmb{\mathbb I}_n
	\end{equation}
	instead of (\ref{eq:constitutiverelationstokes2}), where $\mu > 0$  and $\lambda \geq 0$ are the dynamic and bulk viscosities, respectively, and $\pmb{\mathbb I}_n$ is the identity tensor in $\mathbb R^{n\times n}$. However, note that in this case, 
	\begin{itemize}
		\item if $f_c \equiv 0$, then~\eqref{eqid:weakoriginalstokespb} is the system of equations describing a linear elastic problem in the incompressible limit $\lambda \to \infty$, and the constitutive relation~(\ref{eq:viscousstressgeneral}) then simplifies to~(\ref{eq:constitutiverelationstokes2}); 
		\item since
		$\nabla \cdot \big( \lambda (\nabla\cdot \boldsymbol v)\pmb{\mathbb I}_n\big) = \nabla(\lambda \nabla \cdot \boldsymbol v),$
		if we use the change of variables,
		$p' := p-\lambda \nabla u$ and $\boldsymbol{\sigma}'(\boldsymbol{u}) := 2\mu\boldsymbol{\varepsilon}(\boldsymbol{u})$ as in~(\ref{eq:constitutiverelationstokes2}),
		then problem~\eqref{eqid:weakoriginalstokespb} remains identical if we replace $p$ by $p'$ and the general expression $\boldsymbol{\sigma}(\boldsymbol u)$ from~(\ref{eq:viscousstressgeneral}) by $\boldsymbol{\sigma}'(\boldsymbol u)$ from~(\ref{eq:constitutiverelationstokes2}). 
	\end{itemize}
	Hence, instead of~(\ref{eq:viscousstressgeneral}), we choose the simpler expression~(\ref{eq:constitutiverelationstokes2}) for $\boldsymbol{\sigma}$, without loss of generality.
\end{remark}

\begin{remark}
	The linear elasticity equations can be obtained from Stokes' equations~\eqref{eqid:weakoriginalstokespb} by removing the pressure terms and the divergence condition. In this case, $\boldsymbol{\sigma}(\boldsymbol{v})$ is the Cauchy stress tensor in the medium, satisfying Hooke's law, and in the general form~(\ref{eq:viscousstressgeneral}), $\lambda$ and $\mu$ denote the first and second Lam\'e coefficients, respectively, with $\mu>0$ and $\lambda +\frac{2}{3}\mu>0$. Since the corresponding defeatured problem and the derivation of the defeaturing error estimator for the linear elasticity equations can simply be obtained by removing the pressure terms everywhere, then we will not detail this case and the interested reader is referred to \cite{phdthesis}.
\end{remark}}

Now, let us introduce the corresponding Stokes problem $\mathcal P(\Omega_0)$ in the defeatured geometry $\Omega_0$. To do so, we need to choose an $\boldsymbol L^2$-extension of $\boldsymbol f$ and an $L^2$-extension of $f_c$ in the negative components of the features $F_\mathrm n$, that we still write $\boldsymbol f$ and $f_c$ by abuse of notation. Moreover, we assume that the viscous stress tensor $\tensorsigma$ also satisfies~\eqref{eq:constitutiverelationstokes2} on functions defined everywhere in $\Omega_0$.
Then instead of the exact problem~\eqref{eqid:weakoriginalstokespb}, the following defeatured problem is solved: after choosing $\boldsymbol g_0\in \boldsymbol H^{\frac{1}{2}}\left(\gamma_0\right)$, find the weak solution $(\boldsymbol u_0,p_0)\in \boldsymbol H^1(\Omega_0)\times L^2(\Omega_0)$ of
\begin{align} \label{eqid:simplstokespb}
	\mathcal P(\Omega_0): \quad \begin{cases}
		-\boldsymbol \nabla \cdot \tensorsigma (\boldsymbol u_0) + \nabla p_0 =\boldsymbol f & \text{ in } \Omega_0                 \\
		\boldsymbol \nabla \cdot \boldsymbol u_0 = f_c                                       & \text{ in } \Omega_0                 \\
		\boldsymbol u_0 = \boldsymbol g_D                                                    & \text{ on } \Gamma_D                 \\
		\tensorsigma (\boldsymbol u_0)\mathbf{n} - p_0\mathbf n = \boldsymbol g              & \text{ on } \Gamma_N\setminus \gamma \\
		\tensorsigma (\boldsymbol u_0)\mathbf{n}_0 - p_0\mathbf n = \boldsymbol g_0          & \text{ on } \gamma_0,
	\end{cases}
\end{align}
that is, $(\boldsymbol u_0,p_0)\in \boldsymbol H^1_{\boldsymbol g_D,\Gamma_D}(\Omega_0)\times L^2(\Omega_0)$ satisfies for all $(\boldsymbol v_0,q_0)\in \boldsymbol H^1_{\boldsymbol 0,\Gamma_D}(\Omega_0)\times L^2(\Omega_0)$,
\begin{align}
	\int_{\Omega_0} \tensorsigma (\boldsymbol u_0) : \tensorepsilon(\boldsymbol v_0) \,\mathrm dx - \int_\Omega p_0\boldsymbol \nabla \cdot \boldsymbol{v}_0 \,\mathrm dx & = \int_{\Omega_0} \boldsymbol f \cdot \boldsymbol v_0 \,\mathrm dx + \int_{\Gamma_N\setminus \gamma} \boldsymbol g \cdot \boldsymbol v_0 \,\mathrm ds + \int_{\gamma_0} \boldsymbol g_0 \cdot \boldsymbol v_0 \,\mathrm ds,\nonumber \\
	-\displaystyle\int_\Omega q_0\boldsymbol \nabla \cdot \boldsymbol{u}_0 \,\mathrm dx                                                                                   & = -\int_{\Omega} q_0 f_c\,\mathrm dx. \label{eq:weaksimplstokespb}
\end{align}
By Ladyzhenskaya-Babu\v{s}ka-Brezzi theorem, problem~\eqref{eq:weaksimplstokespb} is well-posed. 

\begin{remark} 
	As for the \changes{Poisson problem}, \review{and even though this is rarely a choice in applications, the best possible choices for} the defeatured problem data $\boldsymbol f$ and $f_c$ in $F_\mathrm n^k$ for $k=1,\ldots,N_f$ and $\boldsymbol g_0$ on $\gamma_0$ will be guided by Remark~\ref{rmk:compatcondmultistokes}. \review{We anticipate that the defeaturing error will be smaller if the defeatured problem data} satisfy a conservation assumption of the solution flux in the positive and negative components of the features. \changes{I.e., if they verify for all $k=1,\ldots,N_f$,
	\begin{align}
	\int_{\gamma_{0,\mathrm p}^{k}} \boldsymbol g_{0} \,\mathrm ds = \int_{\gamma_\mathrm p^k} \boldsymbol g \,\mathrm{d}s + \int_{F_\mathrm p^k} \boldsymbol f\,\mathrm dx \qquad
	\text{ and } \qquad \int_{\gamma_{0,\mathrm n}^{k}} \boldsymbol g_{0} \,\mathrm ds = \int_{\gamma_\mathrm n^{k}} \boldsymbol g \,\mathrm{d}s - \int_{F_\mathrm n^k} \boldsymbol f\,\mathrm dx. \label{eq:compatcondmultilinelast}
	\end{align}}
\end{remark}

Now, for all $k=1,\ldots,N_f$, we need to extend the solution $(\boldsymbol u_0, p_0)$ of~\eqref{eq:weaksimplstokespb} to $\tilde F^k_\mathrm p$ as discussed in Section~\ref{ss:pbpres}, where $\tilde F^k_\mathrm p$ satisfies the properties given in~\eqref{eq:condextension}. To do so, let us choose an $\boldsymbol L^2$-extension of the restriction $\boldsymbol f\vert_{F_\mathrm p^{k}}$ and an $L^2$-extension of the restriction $f_c\vert_{F_\mathrm p^{k}}$ in $G_\mathrm p^k := \tilde F_\mathrm p^k \setminus \overline{F_\mathrm p^k}$, that we still write $\boldsymbol f$ and $f_c$ by abuse of notation. Moreover, we assume that the viscous stress tensor $\tensorsigma$ also satisfies~\eqref{eq:constitutiverelationstokes2} on functions defined in $\tilde F_\mathrm p^k$. Then we define for all $k=1,\ldots,N_f$ the following extension of the solution $(\boldsymbol u_0, p_0)$ of \eqref{eqid:simplstokespb} in $\tilde F^k_\mathrm p$:
after choosing $\tilde{\boldsymbol g}^k\in \boldsymbol H^\frac{1}{2}(\tilde \gamma^k)$, find $(\boldsymbol u_k, p_k)\in \boldsymbol H^1\left(\tilde F_\mathrm p^k\right)\times L^2\left(\tilde F_\mathrm p^k\right)$, the weak solution of
\begin{align} \label{eq:featurestokespb}
	\mathcal P\Big(\tilde F_\mathrm p^k\Big): \quad \begin{cases}
		-\boldsymbol \nabla \cdot \tensorsigma (\boldsymbol u_k) + \nabla p_k =\boldsymbol f                    & \text{ in } \tilde{F}_\mathrm p^k  \\
		\boldsymbol \nabla \cdot \boldsymbol u_k = f_c                                                          & \text{ in } \tilde{F}_\mathrm p^k  \\
		\boldsymbol u_k = \boldsymbol u_0                                                                       & \text{ on } \gamma_{0,\mathrm p}^k \\
		\tensorsigma (\boldsymbol u_k)\tilde{\mathbf{n}}^k - p_k \tilde{\mathbf{n}}^k = \tilde{\boldsymbol g}^k & \text{ on } \tilde \gamma^k        \\
		\tensorsigma (\boldsymbol u_k)\tilde{\mathbf{n}}^k - p_k \tilde{\mathbf{n}}^k = \boldsymbol g           & \text{ on } \gamma_\intersign^k,
	\end{cases}
\end{align}
that is, $(\boldsymbol u_k, p_k)\in \boldsymbol H^1_{\boldsymbol u_0,\gamma_{0,\mathrm p}^k}\left(\tilde F_\mathrm p^k\right)\times L^2\left(\tilde F_\mathrm p^k\right)$
satisfies for all $(\boldsymbol v^k,q^k)\in \boldsymbol H^1_{\boldsymbol 0,\gamma_{0,\mathrm p}^k}\left(\tilde{F}_\mathrm p^k\right)\times L^2\left(\tilde F_\mathrm p^k\right)$,
\begin{align}
	\int_{\tilde{F}_\mathrm p^k} \tensorsigma (\boldsymbol u_k) : \tensorepsilon(\boldsymbol v^k) \,\mathrm dx -\displaystyle\int_{\tilde{F}_\mathrm p^k} p_k\boldsymbol \nabla \cdot \boldsymbol{v}^k \,\mathrm dx & = \int_{\tilde{F}_\mathrm p^k} \boldsymbol f \cdot \boldsymbol v^k \,\mathrm dx + \int_{\tilde \gamma^k} \tilde{\boldsymbol g}^k \cdot \boldsymbol v^k \,\mathrm ds + \int_{\gamma_\intersign^k} \boldsymbol g \cdot \boldsymbol v^k \,\mathrm ds, \nonumber \\
	-\displaystyle\int_{\tilde{F}_\mathrm p^k} q^k\boldsymbol \nabla \cdot \boldsymbol{u}_k \,\mathrm dx                                                                                                            & = -\int_{\tilde{F}_\mathrm p^k} q^k f_c\,\mathrm dx. \label{eq:weakfeaturestokespb}
\end{align}
By Ladyzhenskaya-Babu\v{s}ka-Brezzi theorem, problem~\eqref{eq:weakfeaturestokespb} is well-posed. The choice of the defeatured problem data $\boldsymbol f$ and $f_c$ in $G_\mathrm p^k$, and $\tilde{\boldsymbol g}^k$ on $\tilde \gamma^k$ for $k=1,\ldots,N_f$ will be guided by Remark~\ref{rmk:compatcondmultistokes}. We anticipate here that as before, the best possible choices satisfy a conservation assumption of the solution flux in every feature's extension $G_\mathrm p^k$. \changes{I.e., they verify for all $k=1,\ldots,N_f$,
\begin{align}
\int_{\tilde \gamma^k} \tilde{\boldsymbol g}^k \,\mathrm ds = \int_{\gamma^k_{\setminussign}} \boldsymbol g \,\mathrm{d}s - \int_{\tilde F_\mathrm p^k\setminus F_\mathrm p^k} \boldsymbol f\,\mathrm dx. \label{eq:compatcondmultilinelast2}
\end{align}}

Then, the defeatured solution $(\boldsymbol u_\mathrm d,p_\mathrm d) \in \boldsymbol H_{\boldsymbol g_D,\Gamma_D}(\Omega)\times L^2(\Omega)$ is defined from $(\boldsymbol u_0, p_0)$ and $(\boldsymbol u_k, p_k)$ for $k=1,\ldots,N_f$ as in~\eqref{eqid:defudmulti} and we define the defeaturing error as
\begin{equation}\label{eq:defeaterrorstokes}
	\vertiii{ \left(\boldsymbol u - \boldsymbol u_\mathrm d, p-p_\mathrm d\right)}_\Omega = \left(\int_\Omega \tensorsigma(\boldsymbol u-\boldsymbol u_\mathrm d) : \tensorepsilon(\boldsymbol u-\boldsymbol u_\mathrm d) \,\mathrm dx\right)^\frac{1}{2} + \mu^{-\frac{1}{2}} \|p-p_\mathrm d\|_{0,\Omega}.
\end{equation}

\subsection{Defeaturing error estimator}
Let us recall the definition of the defeatured solution $(\boldsymbol u_\mathrm d, p_\mathrm d)$ from~\eqref{eqid:defudmulti}, and the definitions of $\Sigma$, $\Sigma^k$, $\Sigma_\mathrm n$, $\Sigma_{0,\mathrm p}$, and $\Sigma_\setminussign$ from~\eqref{eq:sigmamulti}. Then \changes{for all $\boundarypiece \in \Sigma$}, let $k_\boundarypiece \equiv k$ if $\boundarypiece \in \Sigma^k$ for some $k=1,\ldots,N_f$, and let us redefine $\boldsymbol d_\boundarypiece$ for all $\boundarypiece \in \Sigma$ in the context of Stokes equations. That is, let
\begin{equation} \label{eq:dsigmastokes}
	\boldsymbol d_\boundarypiece \equiv \begin{cases}\vspace{1mm}
		\boldsymbol g - \tensorsigma(\boldsymbol u_\mathrm d)\mathbf{n} + p_\mathrm d\mathbf{n}                                                         & \text{if } \boundarypiece \in \Sigma_\mathrm n \text{ of if } \boundarypiece\in \Sigma_\setminussign, \\
		-\Big(\boldsymbol g_0 + \tensorsigma (\boldsymbol u_\mathrm d){\mathbf{n}}^{k_\boundarypiece} - p_\mathrm d{\mathbf{n}}^{k_\boundarypiece}\Big) & \text{if } \boundarypiece \in \Sigma_{0,\mathrm p}.
	\end{cases}
\end{equation}
In other words, and as for \changes{Poisson's equations}, $\boldsymbol d_\boundarypiece$ represents the Neumann error on the boundaries $\gamma_\mathrm n$ and $\gamma_\setminussign$ due to the defeaturing process, and the jump of the normal derivative of the defeatured solution on the boundaries $\gamma_{0,\mathrm p}$.
Then, we define the \textit{a posteriori} defeaturing error estimator as:
\begin{equation}\label{eq:multiestimatorstokes}
	\mathscr{E}(\boldsymbol u_\mathrm d, p_\mathrm d) := \left(\sum_{\boundarypiece\in\Sigma} \mathscr{E}_\boundarypiece(\boldsymbol{u}_\mathrm d, p_\mathrm d)^2 \right)^\frac{1}{2},
\end{equation}
where for all $\boundarypiece\in\Sigma$,
\begin{equation}
	\mathscr{E}_\boundarypiece(\boldsymbol u_\mathrm d, p_\mathrm d) := \mu^{-\frac{1}{2}}\left(|\boundarypiece|^\frac{1}{n-1} \left\|\boldsymbol d_\boundarypiece - \overline{\boldsymbol d_\boundarypiece}^\boundarypiece\right\|_{0,\boundarypiece}^2 + c_\boundarypiece^2 |\boundarypiece|^\frac{n}{n-1} \left\|\overline{\boldsymbol d_\boundarypiece}^\boundarypiece\right\|_{\ell^2}^2\right)^\frac{1}{2}, \label{eq:estwithweight}
\end{equation}
where $c_\boundarypiece$ is defined in \eqref{eqp1:negconstant}, $\overline{\boldsymbol d_\boundarypiece}^\boundarypiece$ denotes the dimension-wise average value of $\boldsymbol d_\boundarypiece$ over $\boundarypiece$, and $\|\boldsymbol\cdot\|_{\ell^2}$ denotes the discrete $\ell^2$-norm as $\overline{\boldsymbol d_\boundarypiece}^\boundarypiece\in\mathbb R^n$.
Note that, as for \changes{Poisson's problem}, we can rewrite the estimator feature-wise as
\begin{equation}\label{eq:estforeachfeatstokes}
	\mathscr{E}(\boldsymbol u_\mathrm d, p_\mathrm d) = \left(\sum_{k=1}^{N_f} \sum_{\boundarypiece\in\Sigma^k} \mathscr{E}_{\boundarypiece}(\boldsymbol{u}_\mathrm d, p_\mathrm d)^2 \right)^\frac{1}{2} = \left(\sum_{k=1}^{N_f} \mathscr{E}^{k}(\boldsymbol u_\mathrm d, p_\mathrm d)^2\right)^\frac{1}{2},
\end{equation}
where for all $k=1,\ldots,N_f$, we define $\mathscr{E}^{k}(\boldsymbol u_k, p_k)$ as the defeaturing error estimator for feature $F^k$, that is,
\begin{equation*}
	\mathscr{E}^k(\boldsymbol u_\mathrm d, p_\mathrm d) := \left(\sum_{\boundarypiece\in\Sigma^k} \mathscr{E}_{\boundarypiece}(\boldsymbol{u}_\mathrm d, p_\mathrm d)^2 \right)^\frac{1}{2}.
\end{equation*}

The proposed estimator indicates that all the information on the error introduced by defeaturing is encoded in the boundary of the features, and can be accounted by suitably evaluating the error made on the normal viscous stress and pressure of the solution.

\begin{remark} \label{rmk:compatcondmultistokes}
	Similarly to Remark~\ref{rmk:compatcondmultipoisson}, the terms involving the component-wise average values of $\boldsymbol d_\boundarypiece$ in $\mathscr E(\boldsymbol u_\mathrm d, p_\mathrm d)$ only depend on the defeatured problem data. As a consequence, if these terms dominate, this means that the defeatured problem data should be more accurately chosen. Moreover, under the reasonable data compatibility conditions~\eqref{eq:compatcondmultilinelast} and~\eqref{eq:compatcondmultilinelast2} that represent flux conservation assumptions in this context, the defeaturing error estimator~\eqref{eq:multiestimatorstokes} rewrites
	$\mathscr E(\boldsymbol u_\mathrm d, p_\mathrm d) := \mu^{-\frac{1}{2}}\left( \displaystyle\sum_{\boundarypiece\in\Sigma} \left|\boundarypiece\right|^{\frac{1}{n-1}} \left\| \boldsymbol d_\boundarypiece \right\|^2_{0,\boundarypiece}\right)^\frac{1}{2}.$
\end{remark}

\begin{remark} \label{rmk:estmultitildestokes}
	Similarly to Remark~\ref{rmk:estmultitildepoisson}, note that
	\begin{align*}
		\mathscr{E}(\boldsymbol u_\mathrm d, p_\mathrm d) \lesssim \mu^{-\frac{1}{2}}\left( \displaystyle\sum_{\boundarypiece\in\Sigma} c_\boundarypiece^2 \left|\boundarypiece\right|^{\frac{1}{n-1}} \left\| \boldsymbol d_\boundarypiece \right\|^2_{0,\boundarypiece}\right)^\frac{1}{2} =: \tilde{\mathscr{E}}(\boldsymbol u_\mathrm d, p_\mathrm d).
	\end{align*}
	However, when $n=2$ and under the flux conservation conditions~\eqref{eq:compatcondmultilinelast} and~\eqref{eq:compatcondmultilinelast2}, $\tilde{\mathscr{E}}(\boldsymbol u_\mathrm d, p_\mathrm d)$ is sub-optimal since in this case, $\tilde{\mathscr{E}}(\boldsymbol u_\mathrm d, p_\mathrm d) \lesssim \displaystyle\max_{\boundarypiece\in\Sigma}\left( c_{\boundarypiece} \right) {\mathscr{E}}(\boldsymbol u_\mathrm d, p_\mathrm d)$.
\end{remark}
\changes{\begin{remark}
	If the linear elasticity equations are considered instead of Stokes' equations, with $\lambda$ and $\mu$ denoting the Lam\'e coefficients, then the pressure terms in the Neumann error~(\ref{eq:dsigmastokes}) should be removed, and the weight $\mu^{-\frac{1}{2}}$ in equation~(\ref{eq:estwithweight}) should be replaced by $\rho^{-\frac{1}{2}}$, where $\rho=\mu$ if $\lambda \geq 0$, and $\rho = \min\left(\mu, \frac{3}{2}\lambda + \mu\right)$ otherwise. The coefficient $\rho$ corresponds to the coercivity constant of the bilinear form $\mathfrak a(\boldsymbol{\cdot},\boldsymbol{\cdot})$, \changesbis{see Appendix~\ref{app:proofs}}. Then equations~(\ref{eq:multiestimatorstokes}) and~(\ref{eq:estforeachfeatstokes}) defining the estimator remain the same for the linear elasticity equations.
\end{remark}}

\section{An adaptive geometric refinement strategy} \label{s:adaptivity}
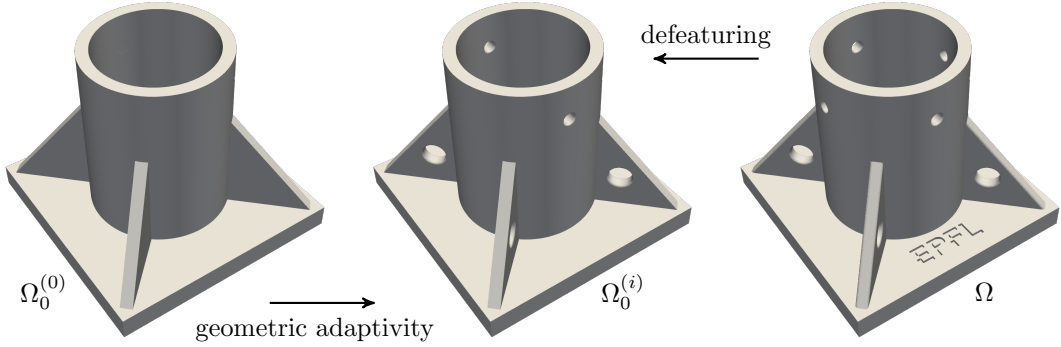
\begin{figure}
	\begin{center}
		\begin{tikzpicture}
		\draw node{\includegraphics[scale=0.115,trim=450 0 400 20, clip]{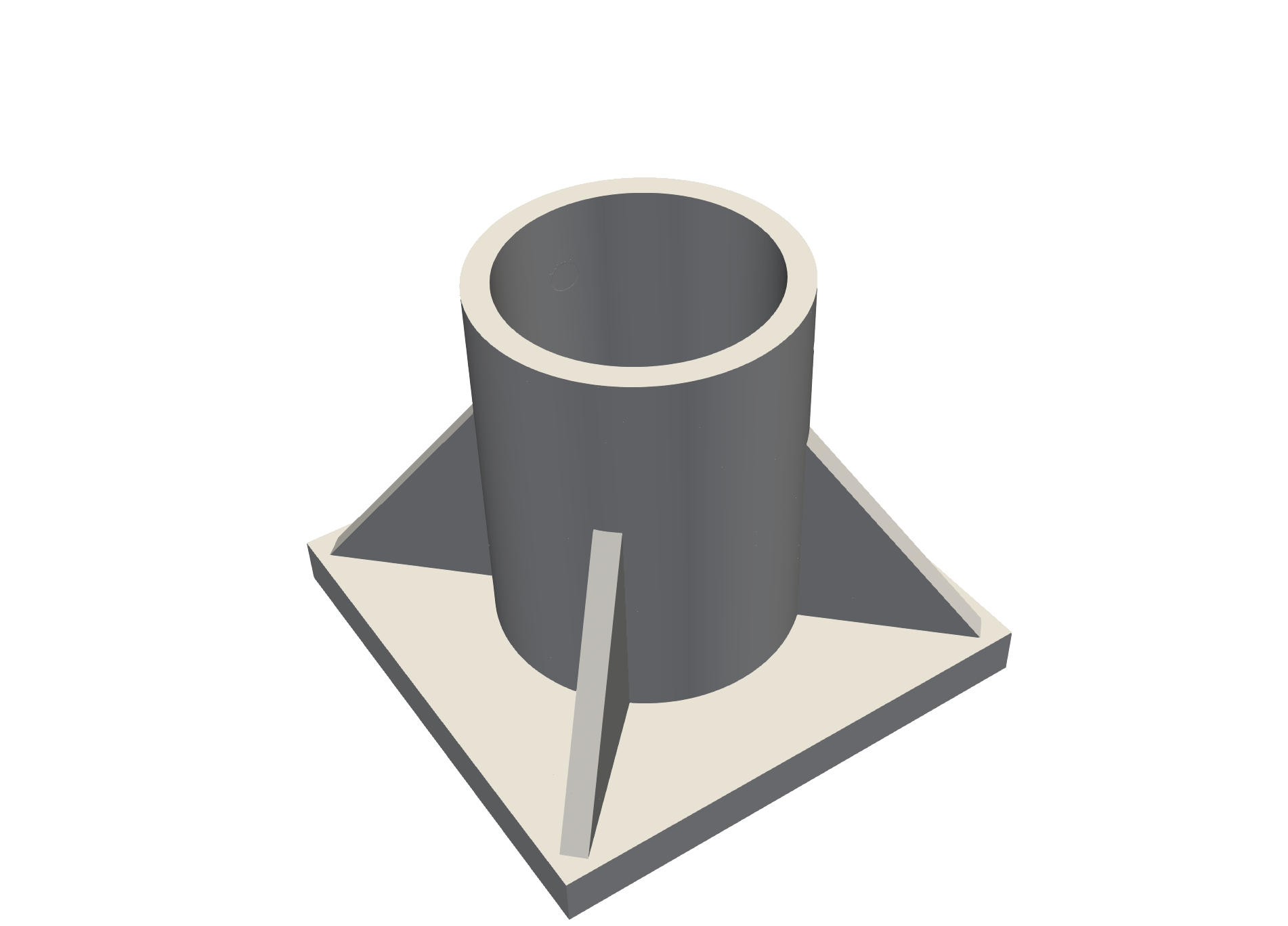}\hspace{5mm} 
			\includegraphics[scale=0.115,trim=450 0 400 20, clip]{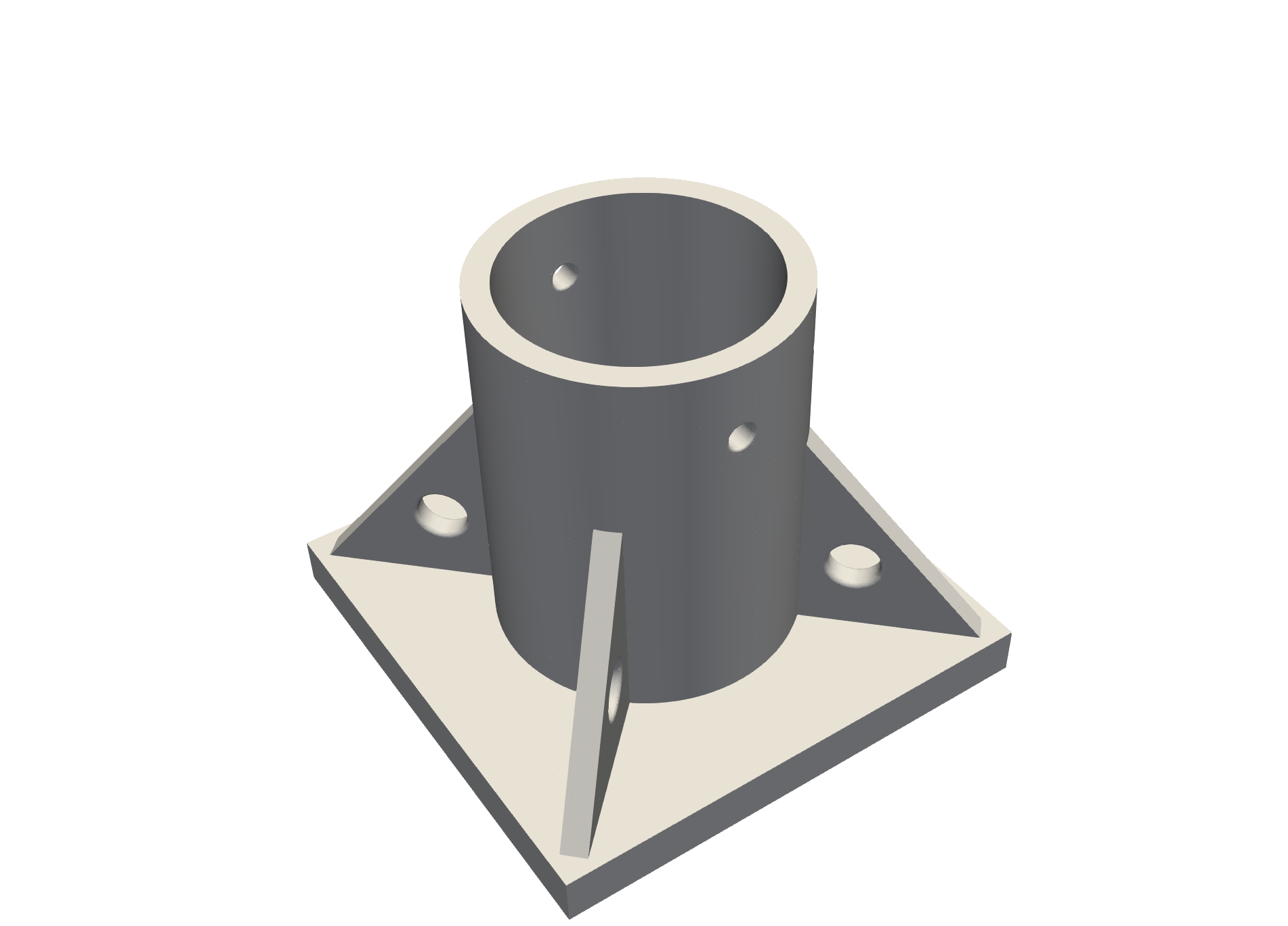}\hspace{5mm} 
			\includegraphics[scale=0.115,trim=450 0 400 20, clip]{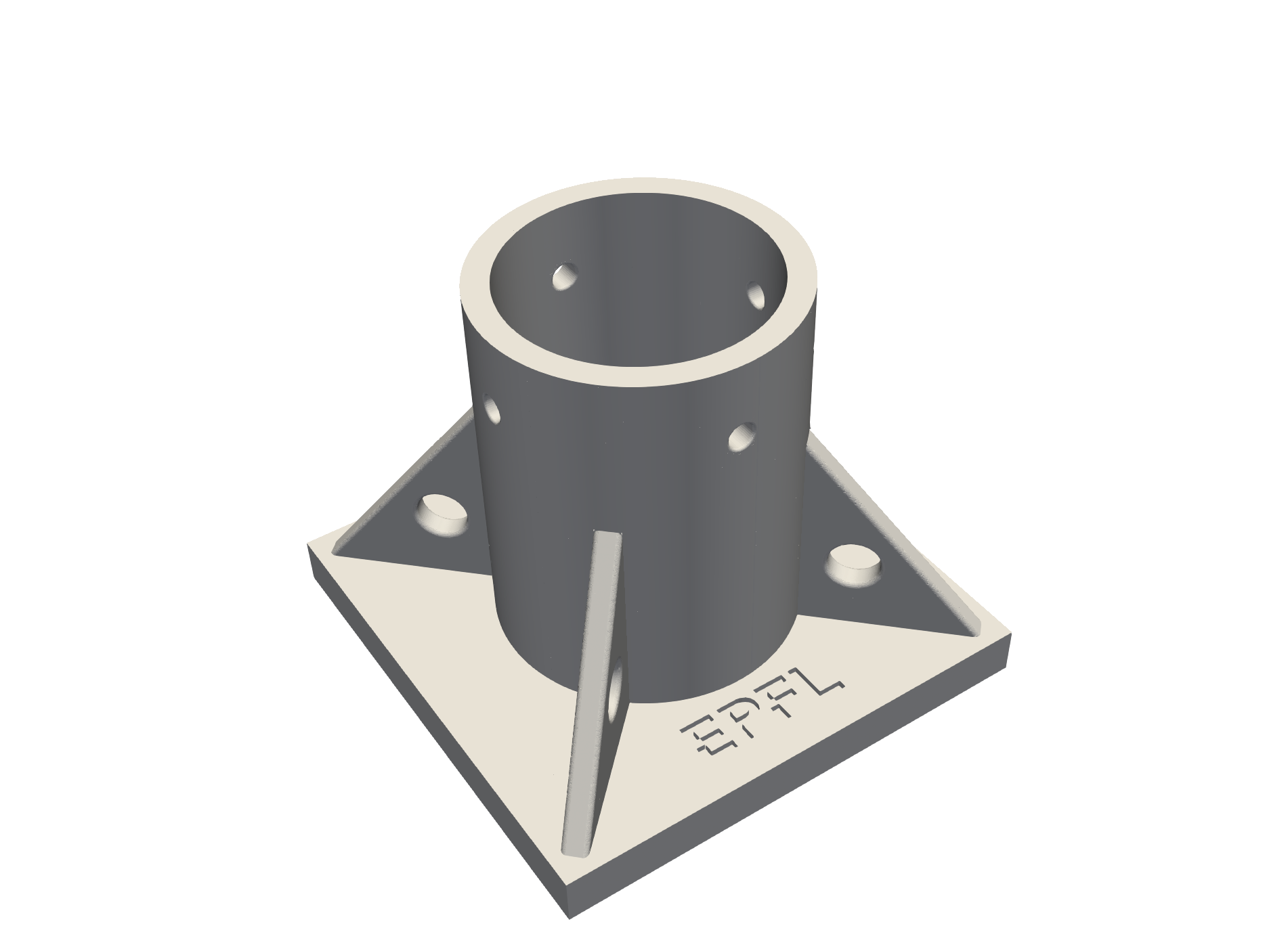}}; 
		\draw[->,
		>=stealth',
		auto,thick] (-3.5,-2.25) -- (-2.1,-2.25) ;
		\draw (-2.9,-2.6) node{{geometric} adaptivity} ;
		\draw[->,
		>=stealth',
		auto,thick] (3,1) -- (1.6,1) ;
		\draw (2.3,1.3) node{defeaturing} ;
		\draw (-6.5,-2.1) node{$\Omega^{(0)}_0$};
		\draw (1.2,-2.1) node{$\Omega^{(i)}_0$};
		\draw (6,-2.1) node{$\Omega$};
		\end{tikzpicture}
		\caption{Illustration of defeaturing and {geometric} adaptivity.} \label{figmf:defeatadaptive}
	\end{center}
\end{figure}

In this section, we aim at defining an adaptive analysis-aware defeaturing strategy in a geometry $\Omega$ containing $N_f\geq 1$ distinct complex features. More precisely, starting from a fully defeatured geometry $\Omega_0$, we want to precisely define a strategy that determines when and which geometrical features need to be reinserted in the geometrical model, among those that have been removed by defeaturing.
Note that the word \textit{defeaturing} may be misleading when thinking of an adaptive strategy: The geometry~$\Omega_0$ in which the problem is actually solved is (partially) defeatured, but the adaptive algorithm selects the features that need to be \textit{added} to the geometrical model, in order to solve the differential problem up to a given accuracy. The concept of geometric adaptivity is illustrated in Figure \ref{figmf:defeatadaptive}.

{In the sequel, we elaborate on each of the building blocks which compose one iteration of an iterative process:
\FloatBarrier
\begin{figure}[h!]
	\begin{center}
		\begin{tikzpicture}[
				start chain = going right,
				node distance=7mm,
				block/.style={shape=rectangle, draw,
						inner sep=1mm, align=center,
						minimum height=5mm, minimum width=15mm, on chain}]
			\node[block] (n1) {SOLVE};
			\node[block] (n2) {ESTIMATE};
			\node[block] (n3) {MARK};
			\node[block] (n4) {REFINE};

			\draw[->] (n1.east) --  + (0,0mm) -> (n2.west);
			\draw[->] (n2.east) --  + (0,0mm) -> (n3.west);
			\draw[->] (n3.east) --  + (0,0mm) -> (n4.west);
			\draw[<-] (n1.south) --  + (0,-5mm) -| (n4.south);
		\end{tikzpicture}
	\end{center}
\end{figure}
\vspace{-0.4cm}
\FloatBarrier}
To do so, let $i\in\mathbb N$ be the current iteration index of the adaptive geometric refinement strategy. For simplicity in this section, let us always write $u_\mathrm d$ the defeatured solution, even in the context of linear elasticity for which it should be $\boldsymbol u_\mathrm d$, or in the context of Stokes equations for which it should be $(\boldsymbol u_\mathrm d, p_\mathrm d)$. To begin the process, let $\Omega_0^{(0)}$ be the fully defeatured geometry defined as in~(\ref{eq:defomega0multi}). That is, $\Omega_0^{(0)}$ is the domain in which all features of $\Omega$ are removed: Their positive component is cut out, and their negative component is filled with material. Since some features will be reinserted during the adaptive process, we denote $\Omega_0^{(i)}$ the simplified geometry at the $i$-th iteration, and in general, we use the upper index $(i)$ to refer to objects at the same iteration. However, to alleviate the notation, we will drop the index $(i)$ when it is clear from the context. In particular, we will write $\Omega_0\equiv \Omega_0^{(i)}$.

\subsection{Solve and Estimate} \label{ssmf:solveestimate}
We first solve the defeatured problem~(\ref{eqid:weaksimplpb})\changes{/
(\ref{eq:weaksimplstokespb})} 
defined in the (partially) defeatured geometry $\Omega_0$. Then, we solve the local extension problem (\ref{eqid:weakfeaturepbmulti})\changes{/
(\ref{eq:weakfeaturestokespb})} for each feature having a non-empty positive component. We thus obtain the defeaturing solution $u_\mathrm d\equiv u_\mathrm d^{(i)}$ defined in (\ref{eqid:defudmulti}), as an approximation of the exact solution~$u$ of (\ref{eqid:weakoriginalpb})\changes{/
(\ref{eqid:weakoriginalstokespb})} at iteration~$i$.
Then, the defeaturing error is estimated by $\mathscr{E}(u_\mathrm d)$ defined in~(\ref{eq:multiestimatorpoisson})\changes{/
(\ref{eq:multiestimatorstokes}).}

\subsection{Mark} \label{ssmf:mark}
Recalling that $N_f \equiv N_f^{(i)}$ at the current iteration $i$, we select and mark some features $\left\{ F^{k_\mathrm m} \right\}_{k_\mathrm m\in I_\mathrm m}\subset \mathfrak F$ with $I_\mathrm m \subset \left\{1,\ldots,N_f^{(i)}\right\}$ to be added to the (partially) defeatured geometry $\Omega_0\equiv\Omega_0^{(i)}$. To do so, we employ in the following a maximum strategy, \changes{but any other convergent marking technique such as the D\"orfler strategy \cite{dorfler} could be used.}
That is, let us first recall definition~\eqref{eq:estforeachfeatpoisson}\changes{/\eqref{eq:estforeachfeatstokes}} of the single feature contributions $\mathscr E^k(u_\mathrm d)$ of the defeaturing error estimator $\mathscr E(u_\mathrm d)$, for $k=1,\ldots,N_f$. Then, after choosing a marking parameter $0<\theta\leq 1$, a feature $F^{k_\mathrm m}$ is marked, i.e., $k_\mathrm m \in I_\mathrm m$, if it verifies
\begin{align}
	\mathscr{E}^{k_\mathrm m}(u_\mathrm d) & \geq \theta \max_{k=1,\ldots,N_f} \left( \mathscr{E}^k(u_\mathrm d)\right). \label{eq:markmultifeat}
\end{align}
In other words, the set of marked features are the ones giving the most substantial contribution to the defeaturing error estimator. The smallest is $\theta$, the more features are selected, and viceversa.

\subsection{Refine} \label{ssmf:refine}
In this step, the defeatured geometry $\Omega_0^{(i)}$ is refined, meaning that the marked features $\left\{ F^{k} \right\}_{k\in I_\mathrm m}$ are inserted in the geometrical model. That is, the new partially defeatured geometrical model $\Omega_0^{(i+1)}$ at the next iteration is built as follows:
\begin{align}
	\Omega_0^{\left(i+\frac{1}{2}\right)} & = \Omega_0^{(i)} \setminus \overline{\bigcup_{k\in I_\mathrm m} F_\mathrm  n^k}, \label{eq:Omega0ip12abs}                                                      \\
	\Omega_0^{(i+1)}                      & = \mathrm{int}\left( \overline{\Omega_0^{\left(i+\frac{1}{2}\right)}} \cup \overline{\bigcup_{k\in I_\mathrm m} F_\mathrm p^k} \right). \label{eq:Omega0p1abs}
\end{align}
And thus in particular,
\begin{align*} 
	F_\mathrm n^{(i+1)}        & := F_\mathrm n^{(i)} \setminus \overline{\bigcup_{k\in I_\mathrm m} {F_\mathrm n^k}}, \quad F_\mathrm p^{(i+1)} := F_\mathrm p^{(i)} \setminus \overline{\bigcup_{k\in I_\mathrm m} {F_\mathrm p^k}}, \\
	\tilde F_\mathrm p^{(i+1)} & := \tilde F_\mathrm p^{(i)} \setminus \overline{\bigcup_{k\in I_\mathrm m} {\tilde F_\mathrm p^k}}, \quad \Omega_\star^{(i+1)} := \Omega \setminus \overline{F_\mathrm p^{(i+1)}},
\end{align*}
and as in definition (\ref{eq:defomega0multi}),
\begin{equation*} 
	\Omega_0^{(i+1)} = \text{int}\left( \overline{\Omega_\star^{(i+1)}} \cup \overline{F_\mathrm n^{(i+1)}} \right).
\end{equation*}

Once the mesh and the defeatured geometry have been refined, the modules SOLVE and ESTIMATE presented in Section \ref{ssmf:solveestimate} can be called again. To do so, we update $\Omega_0$ as $\Omega_0^{(i+1)}$, define $N_f^{(i+1)}:=N_f^{(i)} - \#I_\mathrm m$, update the set of features $\mathfrak F$ as $\mathfrak F \setminus \left\{ F^{k} \right\}_{k\in I_\mathrm m}$, and renumber the features from $1$ to $N_f^{(i+1)}$.
The adaptive loop is continued until a certain given tolerance on the error estimator $\mathscr{E}(u_\mathrm d)$ is reached, or until the set $\mathfrak F$ is empty, meaning that all the features have been added to the geometrical model.

\begin{remark}
	Note that a more precise geometric refinement strategy could be performed since $G_\mathrm p^k$ can be seen as a negative feature of $F_\mathrm p^k$ whose simplified domain is $\tilde F_\mathrm p^k$, for all $k=1,\ldots,N_f$. More precisely, one could consider separately the contributions to $\mathscr{E}^k(u_\mathrm d)$ given by
	\begin{itemize}
		\item $\gamma_\mathrm n^k$ and $\gamma_{0, \mathrm p}^k$, which indicate whether feature $F^k$ should be added to the defeatured geometrical model~$\Omega_0$;
		\item $\gamma_\setminussign^k$, which indicates whether the negative feature $G_\mathrm p^k$ of $F_\mathrm p^k$ should be removed from the simplified positive component $\tilde F_\mathrm p^k$ of $F^k$.
	\end{itemize}
	However, since this adds an extra complexity without introducing new conceptual ideas, this strategy is not further developed in the \review{remainder} of this article.
\end{remark}

\review{\begin{remark}For a given a design, engineers often consider more than one set of boundary conditions, or in other words, they analyze the same design under different loading scenarios. Then, to decide whether the design is valid, they consider an envelope of the results thanks to a max-like function. The proposed adaptive strategy could be adapted to this context by computing the defeaturing error estimator~(\ref{eq:multiestimatorpoisson})/(\ref{eq:multiestimatorstokes}) for each loading scenario. Then, one marks the union of the sets of features marked in each case, as the high-fidelity model used for the analysis should include the features that are important in all the considered loading scenarios.
\end{remark}}

\section{Numerical  considerations and experiments} \label{s:numexp}
In this section, we perform a few numerical experiments to illustrate the validity of the proposed \textit{a posteriori} estimators of the defeaturing error, introduced in Sections~\ref{s:poisson} \changes{and~\ref{s:linelaststokes}}. Thanks to these experiments, we also demonstrate that the adaptive procedure presented in Section~\ref{s:adaptivity} ensures the convergence of the defeaturing error in the energy norm.

\def\CC{{C\nolinebreak[4]\hspace{-.05em}\raisebox{.4ex}{\tiny\bf ++ }}}
For the numerical approximation of the differential problems treated in this section, we use isogeometric analysis (IGA) on very fine meshes, and multipatch and unfitted boundary techniques for the geometrical description of the features. More specifically, a code has been developed on top of GeoPDEs \cite{geopdes}, an open-source and free Octave/MATLAB package for the resolution of
PDEs using IGA. The local meshing process required for the integration of trimmed elements uses the in-house tools presented in \cite{antolinvreps,weimarrusigantolin}, that have been linked to GeoPDEs.
Finally, the $3$D numerical experiment of Section~\ref{sec:3Dtest} has been performed using an in-house C++ library for immersed problems that implements the folded decomposition technique presented in \cite{antolinweibuffa}.
The interested reader is referred to \cite{igabook} for a presentation of isogeometric analysis and advanced spline technologies. 

It is important to remark that, even if in this work we adopted spline based discretizations, the estimator and methodology presented in this work are discretization agnostic, being possible to use other techniques.

\subsection{Impact of some  feature properties on the defeaturing error} \label{ssmf:impactfeat}
In this section, we study the impact of some properties of the geometrical features on the defeaturing error and estimator. In particular, we study the influence of the size and shape of the features, of the distance between them, and of their number.

\subsubsection{Size of the features} \label{sssmf:sizeshape}
\begin{figure}
	\centering
	\begin{subfigure}[t]{0.48\textwidth}
		\begin{center}
			\begin{tikzpicture}[scale=3]
				\draw[thick] (0,0) -- (1,0) ;
				\draw[thick] (0,0) -- (0,1) ;
				\draw[thick] (1,0) -- (1,1) ;
				\draw[thick] (0,1) -- (1,1) ;
				\draw[c1,thick] (0.023, 0.023) circle (0.01);
				\draw (0.11,0.09) node{$F^1$};
				\draw[c1,thick] (0.88,0.88) circle (0.1);
				\draw (0.88,0.88) node{$F^2$};
				\draw (0.5,0.5) node{$\Omega$};
			\end{tikzpicture}
			\caption{Exact domain with two circular features.}\label{fig:circles}
		\end{center}
	\end{subfigure}
	~
	\begin{subfigure}[t]{0.48\textwidth}
		\begin{center}
			\begin{tikzpicture}[scale=3]
				\draw [thick](0,0) -- (1,0) ;
				\draw[thick] (0,0) -- (0,1) ;
				\draw[thick] (1,0) -- (1,1) ;
				\draw[thick] (0,1) -- (1,1) ;
				\draw[c1,thick] (0.013, 0.013) -- (0.033,0.013) -- (0.033,0.033) -- (0.013,0.033) -- cycle;
				\draw (0.11,0.09) node{$F^1$};
				\draw[c1,thick] (0.78,0.78) -- (0.98,0.78) -- (0.98,0.98) -- (0.78,0.98) -- cycle;
				\draw (0.88,0.88) node{$F^2$};
				\draw (0.5,0.5) node{$\Omega$};
			\end{tikzpicture}
			\caption{Exact domain with two square features.}\label{fig:squares}
		\end{center}
	\end{subfigure}
	\caption{Numerical test~\ref{sssmf:sizeshape} -- Exact geometries used for the comparison between features' sizes (not at scale).} \label{fig:circlessquares}
\end{figure}
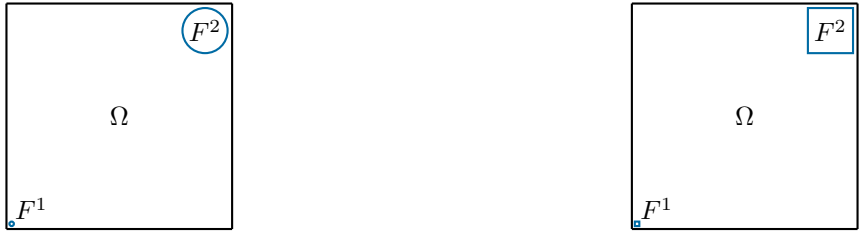

Let us consider a numerical experiment first studied in~\cite{paper1defeaturing}, in which a computational domain contains a very small but important feature, and a large feature whose presence or absence does not affect much the solution accuracy. More precisely, let $\Omega_0 := (0,1)^2$ be the fully defeatured geometry, and let $\Omega := \Omega_0 \setminus \big(F^1\cup F^2\big)$, where $F^1$ is a circular hole of radius $10^{-3}$ centered at $(1.1\cdot 10^{-3}, 1.1\cdot 10^{-3})$, and $F^2$ is a larger circular hole of radius $10^{-1}$ centered at $(8.9\cdot 10^{-1}, 8.9\cdot 10^{-1})$. The considered geometry is illustrated in Figure~\ref{fig:circles}. 
We consider Poisson's problem~(\ref{eqid:weakoriginalpb}) in the exact computational domain $\Omega$, and its defeatured version~(\ref{eqid:weaksimplpb}) defined in $\Omega_0$, with $f(x,y) := 128e^{-8(x+y)}$ in $\Omega_0$, $g_D(x,y) := e^{-8(x+y)}$ on $\Gamma_D := \big([0,1)\times \{0\} \big) \cup \big(\{0\}\times [0,1)\big),$
the bottom and left sides, $g(x,y) := -8e^{-8(x+y)}$ on $\partial \Omega_0 \setminus \overline{\Gamma_D}$, and finally $g\equiv 0$ on $\partial F^1\cup \partial F^2$. The exact and defeatured solutions $u$ and $u_0$ of these Poisson's problems have a high gradient in the region around the small feature $F^1$, while they are almost identically equal to zero in the region around the large feature $F^2$.

We then perform the same test, but with square holes instead of circular ones. That is, we consider the same defeatured geometry $\Omega_0=(0,1)^2$, and the same data to solve Poisson's equation~(\ref{eqid:weakoriginalpb}), but now $\Omega := \Omega_0 \setminus (\overline{F^1} \cup \overline{F^2})$, where $F^1$ and $F^2$ are squares centered at $(1.1\cdot 10^{-3}, 1.1\cdot 10^{-3})$ and $(8.9\cdot 10^{-1}, 8.9\cdot 10^{-1})$, respectively, and whose sides have length $2\cdot10^{-3}$ and $2\cdot10^{-1}$, respectively. The geometry is illustrated in Figure~\ref{fig:squares}, and as before, the solution has a high gradient close to the bottom left corner where $F^1$ is located, and is almost constantly equal to zero close to the top right corner where $F^2$ is located.

\begin{table}
	\centering
	\begin{tabular}{@{}cccccc@{}}
		\hline \rule{0pt}{12pt}
		Features                              & $\mathscr{E}^1(u_0)$ & $\mathscr{E}^2(u_0)$ & $\mathscr{E}(u_0)$  & $\vertiii{u-u_\mathrm d}_\Omega$ & $\eta_\mathrm{eff}\;$ \\[5pt]
		\hline \rule{0pt}{13pt}Circular holes & $5.03\cdot 10^{-2}$  & $7.86\cdot 10^{-6}$  & $5.03\cdot 10^{-2}$ & $1.45\cdot 10^{-2}$              & $3.47\;$              \\ 
		Square holes                          & $6.29\cdot 10^{-2}$  & $7.73\cdot 10^{-6}$  & $6.29\cdot 10^{-2}$ & $1.64\cdot 10^{-2}$              & $3.84\;$              \\ 
		\hline
	\end{tabular}
	\caption{Numerical test~\ref{sssmf:sizeshape} -- Results of the comparison between features' sizes.} \label{tbl:sizefeatshape}
\end{table}

The values of the defeaturing error estimator~(\ref{eq:multiestimatorpoisson}) and of the defeaturing error $\vertiii{u-u_\mathrm d}_\Omega$ are reported in Table~\ref{tbl:sizefeatshape}. In both geometries, independently of the shape of the features, $F^1$ is indeed more important than $F^2$ since the estimator for $F^1$ is four orders of magnitude larger than the estimator for $F^2$. This result was expected because of the solution's very high gradient close to $F^1$, and because of the homogeneous boundary conditions imposed on the feature's boundaries. In both cases, the proposed estimator well estimates the defeaturing error since the effectivity index $\eta_{\mathrm{eff}}:=\displaystyle\frac{\mathscr{E}(u_0)}{\vertiii{u-u_\mathrm d}_\Omega}$ is reasonably low, with values comparable to the single feature experiments performed in~\cite{paper1defeaturing}. These results perfectly agree with the theory developed in Section~\ref{s:poisson}, and illustrate how the proposed estimator does not only depend on geometrical considerations, but also on the analysis behind, i.e., on the considered differential problem that one wants to solve. Hence the name \textit{analysis-aware} defeaturing.

\subsubsection{Distance between features} \label{sec:featseparation}
The following numerical example is used to show that the separability Assumption~\ref{as:separatednew} is very weak, as one can consider features that are arbitrarily close to one another, as soon as the number of close features is bounded. Indeed, consider a geometry with either two square features, one positive and one negative, or one complex feature, as follows. 
Let $\Omega_0 := (0,1)^2$, let $\delta\in (-0.1,0.8)$, and let $\Omega_\delta := \mathrm{int}\left(\overline{\Omega_0}\cup \overline{F^1_\delta} \setminus \overline{F^2_\delta}\right)$ with
\begin{align*}
	F^1_\delta & := \left(0.4-\frac{\delta}{2}, 0.5-\frac{\delta}{2}\right) \times \left(1,1.1\right)\quad \text{and } \quad 
	F^2_\delta := \left(0.5+\frac{\delta}{2}, 0.6+\frac{\delta}{2}\right) \times \left(0.9,1\right), 
\end{align*}
as illustrated in Figure \ref{figmf:twoclosefeat}. That is,
\begin{itemize}
	\item if $\delta\leq0$, then ${F^1_\delta}\cup{F^2_\delta}$ needs to be considered as a single feature because of Assumption~\ref{as:separatednew}, where $F_\delta^1$ is the positive component of that feature, and $F_\delta^2$ is its negative component. In this case, we let $\gamma_{0,\delta}^1 := \gamma_{0,\mathrm p}$, $\gamma_{0,\delta}^2 := \gamma_{0,\mathrm n}$, $\gamma_{\delta}^1 := \gamma_{\mathrm p}$, and $\gamma_{\delta}^2 := \gamma_{\mathrm n}$;
	\item if $\delta>0$, then $F^1_\delta$ and $F^2_\delta$ are two distinct features satisfying Assumption~\ref{as:separatednew} and separated by a distance $\delta$, where $F^1_\delta$ is positive, and  $F^2_\delta$ is negative. In this case, we let $\gamma_{0,\delta}^1 := \gamma^1_{0,\mathrm p}=\gamma^1_{0}$, $\gamma_{0,\delta}^2 := \gamma^2_{0,\mathrm n} =\gamma^2_{0}$, $\gamma_{\delta}^1 := \gamma^1_{\mathrm p} = \gamma^1$, and $\gamma_{\delta}^2 := \gamma^2_{\mathrm n}=\gamma^2$.
\end{itemize}
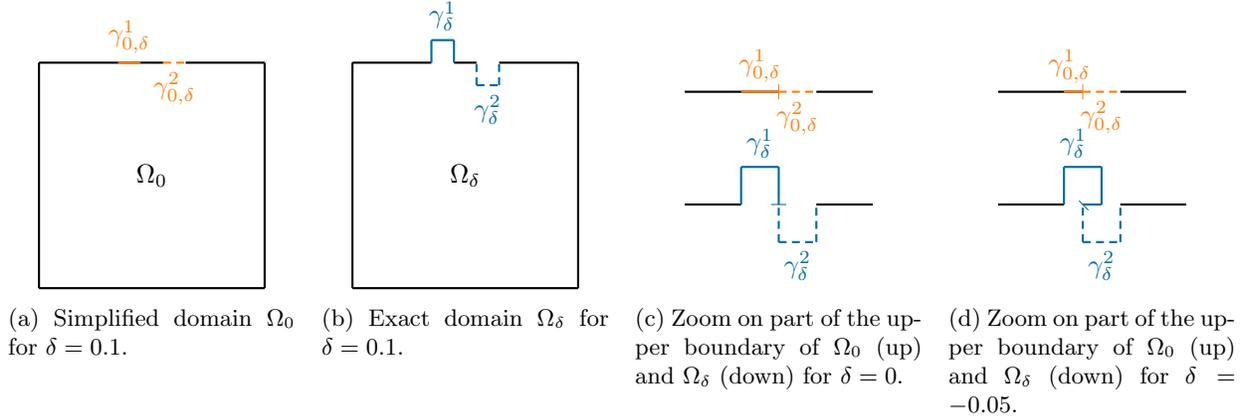
\begin{figure}
	\begin{center}
		\begin{subfigure}[t]{0.23\textwidth}
			\begin{center}
				\begin{tikzpicture}[scale=3]
					\draw[thick] (0,0) -- (1,0) ;
					\draw[thick] (0,0) -- (0,1) ;
					\draw[thick] (1,0) -- (1,1) ;
					\draw[thick] (0,1) -- (0.55,1) ;
					\draw[thick] (1,1) -- (0.65,1) ;
					\draw[c2,thick] (0.35,1) -- (0.45,1) ;
					\draw[c2,thick] (0.4,1) node[above]{$\gamma_{0,\delta}^1$} ;
					\draw[c2,thick,densely dashed] (0.55,1) -- (0.65,1) ;
					\draw[c2,thick] (0.6,1) node[below]{$\gamma_{0,\delta}^2$} ;
					\draw (0.5,0.5) node{$\Omega_0$} ;
				\end{tikzpicture}
				\caption{Simplified domain $\Omega_0$ for $\delta=0.1$.}
			\end{center}
		\end{subfigure}
		~
		\begin{subfigure}[t]{0.23\textwidth}
			\begin{center}
				\begin{tikzpicture}[scale=3]
					\draw[thick] (0,0) -- (1,0) ;
					\draw[thick] (0,0) -- (0,1) ;
					\draw[thick] (1,0) -- (1,1) ;
					\draw[thick] (0,1) -- (0.35,1) ;
					\draw[thick] (0.45,1) -- (0.55, 1) ;
					\draw[thick] (0.65,1) -- (1,1) ;
					\draw[c1,thick] (0.35,1.1)--(0.45,1.1);
					\draw[c1,thick] (0.35,1.1)--(0.35,1);
					\draw[c1,thick] (0.45,1.1)--(0.45,1);
					\draw[c1,thick,densely dashed] (0.55,0.9)--(0.65,0.9);
					\draw[c1,thick,densely dashed] (0.55,0.9)--(0.55,1);
					\draw[c1,thick,densely dashed] (0.65,0.9)--(0.65,1);
					\draw (0.5,0.5) node{$\Omega_\delta$} ;
					\draw[c1,thick] (0.4,1.1) node[above]{$\gamma_\delta^1$} ;
					\draw[c1,thick] (0.6,0.9) node[below]{$\gamma_\delta^2$} ;
				\end{tikzpicture}
				\caption{Exact domain $\Omega_\delta$ for $\delta=0.1$.}
			\end{center}
		\end{subfigure}
		~
		\begin{subfigure}[t]{0.23\textwidth}
			\begin{center}
				\begin{tikzpicture}[scale=5]
					\draw[thick] (0.25,1.3) -- (0.5,1.3);
					\draw[thick] (0.6,1.3) -- (0.75,1.3) ;
					\draw[c2,thick] (0.4,1.3) -- (0.5,1.3) ;
					\draw[c2,thick] (0.45,1.3) node[above]{$\gamma_{0,\delta}^1$} ;
					\draw[c2,thick,densely dashed] (0.5,1.3) -- (0.6,1.3) ;
					\draw[c2,thick] (0.55,1.3) node[below]{$\gamma_{0,\delta}^2$} ;
					\draw[thick] (0.25,1) -- (0.4,1) ;
					\draw[thick] (0.6,1) -- (0.75,1) ;
					\draw[c1,thick] (0.4,1.1)--(0.5,1.1);
					\draw[c1,thick] (0.4,1.1)--(0.4,1);
					\draw[c1,thick] (0.5,1.1)--(0.5,1);
					\draw[c1,thick,densely dashed] (0.5,0.9)--(0.6,0.9);
					\draw[c1,thick,densely dashed] (0.5,0.9)--(0.5,1);
					\draw[c1,thick,densely dashed] (0.6,0.9)--(0.6,1);
					\draw[c1] (0.45,1.1) node[above]{$\gamma_\delta^1$} ;
					\draw[c1] (0.55,0.9) node[below]{$\gamma_\delta^2$} ;
					\draw[c2] (0.5,1.28) -- (0.5,1.32) ;
					\draw[c1] (0.48,1) -- (0.52, 1) ;
				\end{tikzpicture}
				\caption{Zoom on part of the upper boundary of $\Omega_0$ (up) and $\Omega_\delta$ (down) for $\delta=0$.}
			\end{center}
		\end{subfigure}
		~
		\begin{subfigure}[t]{0.23\textwidth}
			\begin{center}
				\begin{tikzpicture}[scale=5]
					\draw[thick] (0.25,1.3) -- (0.475,1.3);
					\draw[thick] (0.575,1.3) -- (0.75,1.3) ;
					\draw[c2,thick] (0.425,1.3) -- (0.475,1.3) ;
					\draw[c2,thick,densely dashed] (0.475,1.3) -- (0.575,1.3) ;
					\draw[c2] (0.44,1.3) node[above]{$\gamma_{0,\delta}^1$} ;
					\draw[c2] (0.525,1.3) node[below]{$\gamma_{0,\delta}^2$} ;
					\draw[c2] (0.475,1.28) -- (0.475,1.32);
					\draw[thick] (0.25,1) -- (0.425,1) ;
					\draw[c1,thick] (0.475,1) -- (0.525, 1) ;
					\draw[thick] (0.575,1) -- (0.75,1) ;
					\draw[c1] (0.465,1.01) -- (0.49,0.985);
					\draw[c1,thick] (0.425,1.1)--(0.525,1.1);
					\draw[c1,thick] (0.425,1.1)--(0.425,1);
					\draw[c1,thick] (0.525,1.1)--(0.525,1);
					\draw[c1,thick,densely dashed] (0.475,0.9)--(0.575,0.9);
					\draw[c1,thick,densely dashed] (0.475,0.9)--(0.475,1);
					\draw[c1,thick,densely dashed] (0.575,0.9)--(0.575,1);
					\draw[c1] (0.45,1.1) node[above]{$\gamma_\delta^1$} ;
					\draw[c1] (0.525,0.9) node[below]{$\gamma_\delta^2$} ;
				\end{tikzpicture}
				\caption{Zoom on part of the upper boundary of $\Omega_0$ (up) and $\Omega_\delta$ (down) for $\delta=-0.05$.}
			\end{center}
		\end{subfigure}
		\caption{Numerical test~\ref{sec:featseparation} -- Simplified domain $\Omega_0$ and exact domains $\Omega_\delta$ for different values of $\delta$.} \label{figmf:twoclosefeat}
	\end{center}
\end{figure}
\begin{table}[t]
	\centering
	\begin{tabular}{c c c c}
		\hline \rule{0pt}{12pt}
		$\delta$                                        & $\mathscr{E}(u_\mathrm d)$ & $\vertiii{u-u_\mathrm d}_{\Omega_\delta}$ & $\eta_{\mathrm{eff}}$ \\ [0.5ex]
		\hline \rule{0pt}{12pt}$\,\,\,2.0\cdot 10^{-1}$ & $1.58$                     & $1.49$                                    & $1.73$                \\
		$\,\,\,2.0\cdot 10^{-4}$                        & $2.84$                     & $1.68$                                    & $1.69$                \\
		$0.0 \cdot 10^0$                                & $2.84$                     & $1.68$                                    & $1.69$                \\
		$-1.0\cdot 10^{-3}$                             & $27.0$                     & $15.1$                                    & $1.78$                \\
		$-9.9\cdot 10^{-2}$                             & $24.5$                     & $14.3$                                    & $1.71$                \\
		\hline
	\end{tabular}
	\caption{Numerical test~\ref{sec:featseparation} -- Values of the defeaturing error and estimator for different values of $\delta$. The cases in which $\delta>0$ correspond to separate features, while the cases $\delta<0$ correspond to features with overlapping boundaries.} \label{table:twodifffeat}
\end{table}
Let us consider Poisson problem~(\ref{eqid:weakoriginalpb}) with $f\equiv 0$ in $\Omega$, $g_D(x,y):= 40\cos(\pi x)+10\cos(5\pi x)$ on
$\Gamma_D := (0,1)\times \{0\},$ 
and $g\equiv 0$ on $\Gamma_N:=\partial \Omega_\delta\setminus \overline{\Gamma_D}.$
We solve the defeatured Poisson problem (\ref{eqid:weaksimplpb}) with the same data, and we take $g_0\equiv 0$ on
\begin{align*}
	\gamma_{0,\delta}^1                    & = \left(0.4-\displaystyle\frac{\delta}{2}, 0.5-\displaystyle\frac{|\delta|}{2}\right) \times \{1\} \quad 
	\text{ and } \quad \gamma_{0,\delta}^2 = \left(0.5+\displaystyle\frac{\delta}{2}, 0.6+\displaystyle\frac{\delta}{2}\right) \times \{1\}.  
\end{align*}
Finally, we solve the Dirichlet extension problem~(\ref{eqid:weakfeaturepbmulti}) in $\tilde F^1_\delta = F^1_\delta$.
We choose different values of $\delta$ in order to consider different cases:
\begin{itemize}
	\item with $\delta=2\cdot 10^{-1}$, the distance between the features and the distance between $\gamma_{0,\delta}^1$ and $\gamma_{0,\delta}^2$ are of the same order of magnitude as the measures of $\gamma_{0,\delta}^1$ and $\gamma_{0,\delta}^2$;
	\item with $\delta=2\cdot 10^{-4}$, the distance between $\gamma_{0,\delta}^1$ and $\gamma_{0,\delta}^2$ is several orders of magnitude smaller than the measures of $\gamma_{0,\delta}^1$ and $\gamma_{0,\delta}^2$;
	\item with $\delta=0$, the boundaries of the feature components intersect in one single point;
	\item with $\delta=-1\cdot 10^{-3}$, the measure of the intersection between the boundaries of the feature components is several orders of magnitude smaller than the measures of the boundaries of the features;
	\item with $\delta=-9.9\cdot 10^{-2}$, the measure of the intersection between the boundaries of the feature components is of the same order of magnitude as the measures of the boundaries of the features.
\end{itemize}
The results are presented in Table~\ref{table:twodifffeat}, and we indeed see that the defeaturing estimator approximates well the defeaturing error in all the different presented cases. In particular, we observe that the effectivity index $\eta_{\mathrm{eff}}:=\displaystyle\frac{\mathscr{E}(u_\mathrm d)}{\vertiii{u-u_\mathrm d}_\Omega}$ is not influenced by the distance separating the positive and negative components of the feature(s). This confirms the fact that Assumption~\ref{as:separatednew} in not very restrictive in practice.

\subsubsection{Number of features} \label{sssmf:nbfeat}	
Under Assumption~\ref{as:separatednew}, the effectivity index of the defeaturing error estimator should not depend on the number of features that are present in the original geometry $\Omega$. To verify this, let $\Omega_0:= (0,1)^2$ be the fully defeatured domain, and let $\Omega := \Omega_0 \setminus \bigcup_{k=1}^{N_f} \overline{F^k}$, where $N_f=27$, and the features $F^k$ are circular holes of random radii in the interval $(0,0.01)$ which are randomly distributed in $\Omega_0$, under the condition that Assumption~\ref{as:separatednew} is satisfied. For the sake of reproducibility, the values of the radii and centers of the features are reported in Table~\ref{tbl:radiusescenters}. The exact domain $\Omega$ with all the $27$ features is represented in Figure~\ref{fig:geom81holes}.

We want to find a good approximation of the solution of Poisson's problem~(\ref{eqid:weakoriginalpb}) in $\Omega$, whose exact solution is shown in Figure~\ref{fig:sol81holes}, being $f(x,y) := -18e^{-3(x+y)}$ in $\Omega_0$, $g_D(x,y) := e^{-3(x+y)}$ on the bottom and left boundaries, i.e., on
$\Gamma_D := \big( [0,1)\times\{0\} \big) \cup \big( \{0\}\times [0,1) \big),$ 
$g(x,y):=-3e^{-3(x+y)}$ on $\partial \Omega_0 \setminus \overline{\Gamma_D}$, and $g \equiv 0$ on $\partial F^k$ for $k=1,\ldots,N_f$.
Thus we perform the adaptive algorithm introduced in Section~\ref{s:adaptivity} starting from the fully defeatured domain $\Omega_0^{(0)}:=\Omega_0$, with marking parameter $\theta = 0.95$.
I.e., at every iteration we include the $5\%$ of the features whose error contributes the most.
We recursively solve the partially defeatured problem~(\ref{eqid:weaksimplpb}) in $\Omega_0^{(i)}$ at each iteration $i\geq 0$, and we call $u_0^{(i)}$ its solution.

{\begin{table}
	\centering
	{\small\def\arraystretch{1.2}
		\begin{tabular}{@{}cccccccccc@{}}
			\hline
			Feature index $k$        & $1$    & $2$    & $3$    & $4$    & $5$    & $6$    & $7$    & $8$    & $9$    \\
			\hline
			Radius $[\cdot 10^{-2}]$ & $8.13$ & $6.64$ & $3.89$ & $7.40$ & $8.18$ & $6.00$ & $0.85$ & $9.22$ & $0.54$ \\
			Center $[\cdot 10^{-2}]$ & $0.98$ & $2.84$ & $5.46$ & $7.16$ & $8.99$ & $0.67$ & $3.12$ & $4.95$ & $7.06$ \\
			                         & $0.93$ & $1.24$ & $0.57$ & $0.93$ & $1.04$ & $3.40$ & $3.03$ & $3.08$ & $2.48$ \\
			\hline\hline
			Feature index $k$        & $10$   & $11$   & $12$   & $13$   & $14$   & $15$   & $16$   & $17$   & $18$   \\
			\hline
			Radius $[\cdot 10^{-2}]$ & $5.27$ & $1.19$ & $3.80$ & $8.13$ & $2.44$ & $8.84$ & $7.13$ & $3.78$ & $2.49$ \\
			Center $[\cdot 10^{-2}]$ & $8.86$ & $0.67$ & $3.28$ & $5.01$ & $7.44$ & $8.93$ & $1.10$ & $2.44$ & $5.45$ \\
			                         & $2.90$ & $5.35$ & $4.46$ & $5.09$ & $4.88$ & $5.07$ & $6.93$ & $6.78$ & $7.73$ \\
			\hline\hline
			Feature index $k$        & $19$   & $20$   & $21$   & $22$   & $23$   & $24$   & $25$   & $26$   & $27$   \\
			\hline
			Radius $[\cdot 10^{-2}]$ & $2.53$ & $6.67$ & $0.50$ & $6.85$ & $6.20$ & $7.47$ & $8.77$ & $2.00$ & $1.00$ \\
			Center $[\cdot 10^{-2}]$ & $7.27$ & $9.21$ & $0.22$ & $3.26$ & $5.01$ & $7.06$ & $8.99$ & $4.00$ & $1.00$ \\
			                         & $7.33$ & $6.96$ & $8.24$ & $9.15$ & $9.10$ & $8.78$ & $8.98$ & $7.00$ & $9.00$ \\
			\hline
		\end{tabular}}
	\caption{Numerical test~\ref{sssmf:nbfeat} -- Data of the $27$ circular features.} \label{tbl:radiusescenters}
\end{table}
	\begin{figure}
			\begin{center}
				\begin{subfigure}[t!]{0.3\textwidth}
					\begin{center}
						\begin{tikzpicture}[scale=4.8]
							\draw[thick] (0,0) -- (0,1) -- (1,1) -- (1,0) -- cycle;
							\draw[c1,thick] (0.097984328699318,0.092663690863534) circle (0.081318234590090);
							\draw[c1,thick] (0.283884293965162,0.123563491365580) circle (0.066396825383881);
							\draw[c1,thick] (0.546112220339815,0.056574896706227) circle (0.038919282299261);
							\draw[c1,thick] (0.715916527014347,0.093265457598229) circle (0.074000758141615);
							\draw[c1,thick] (0.898581968051930,0.104463986172862) circle (0.081763485970873);
							\draw[c1,thick] (0.067305141757723,0.339769841694328) circle (0.060034479317649);
							\draw[c1,thick] (0.311761232757188,0.303173986566693) circle (0.008499706450270);
							\draw[c1,thick] (0.495145556524009,0.307616854246403) circle (0.092235800109005);
							\draw[c1,thick] (0.705998218664498,0.248238323585443) circle (0.005359781531656);

							\draw[c1,thick] (0.886286875701112,0.290351799644088) circle (0.052702495267737);
							\draw[c1,thick] (0.067359418419866,0.534641171336466) circle (0.011885327661975);
							\draw[c1,thick] (0.328109457773039,0.446027904596783) circle (0.038014299024100);
							\draw[c1,thick] (0.500590433180793,0.509270838208646) circle (0.081283253301511);
							\draw[c1,thick] (0.743939933090058,0.487966072170030) circle (0.024409591557169);
							\draw[c1,thick] (0.893169207667069,0.507196247921792) circle (0.088442265909357);
							\draw[c1,thick] (0.110235870451937,0.693080838594426) circle (0.071264679475915);
							\draw[c1,thick] (0.244342487141727,0.677497408280121) circle (0.037814843808214);
							\draw[c1,thick] (0.545240931016100,0.773012684854573) circle (0.024891961099593);

							\draw[c1,thick] (0.726683377942314,0.732602609237550) circle (0.025285374523927);
							\draw[c1,thick] (0.920762287773320,0.695958574121097) circle (0.066724359452694);
							\draw[c1,thick] (0.022384764760657,0.823728658042989) circle (0.004986188504724);
							\draw[c1,thick] (0.325708044929628,0.914763687652304) circle (0.068528856760891);
							\draw[c1,thick] (0.500755877184029,0.910427640482701) circle (0.062027812215100);
							\draw[c1,thick] (0.705821367824273,0.878409519462713) circle (0.074668462996873);
							\draw[c1,thick] (0.899163331742718,0.898273740857340) circle (0.087725564834567);
							\draw[c1,thick] (0.4,0.7) circle (0.020000000000000);
							\draw[c1,thick] (0.1,0.9) circle (0.010000000000000);

							\draw[c1,thick] (0.097984328699318,0.092663690863534) node{\tiny $1$};
							\draw[c1,thick] (0.283884293965162,0.123563491365580) node{\tiny $2$};
							\draw[c1,thick] (0.546112220339815,0.056574896706227) node{\tiny $3$};
							\draw[c1,thick] (0.715916527014347,0.093265457598229) node{\tiny $4$};
							\draw[c1,thick] (0.898581968051930,0.104463986172862) node{\tiny $5$};
							\draw[c1,thick] (0.067305141757723,0.339769841694328) node{\tiny $6$};
							\draw[c1,thick] (0.3,0.303173986566693) node[right]{\tiny $7$}; 
							\draw[c1,thick] (0.495145556524009,0.307616854246403) node{\tiny $8$};
							\draw[c1,thick] (0.715,0.248238323585443) node[left]{\tiny $9$}; 

							\draw[c1,thick] (0.886286875701112,0.290351799644088) node{\tiny $10$};
							\draw[c1,thick] (0.067359418419866,0.534641171336466) node[right]{\tiny $11$};
							\draw[c1,thick] (0.328109457773039,0.446027904596783) node{\tiny $12$};
							\draw[c1,thick] (0.500590433180793,0.509270838208646) node{\tiny $13$};
							\draw[c1,thick] (0.743939933090058,0.487966072170030) node[left]{\tiny $14$};
							\draw[c1,thick] (0.893169207667069,0.507196247921792) node{\tiny $15$};
							\draw[c1,thick] (0.110235870451937,0.693080838594426) node{\tiny $16$};
							\draw[c1,thick] (0.244342487141727,0.677497408280121) node{\tiny $17$};
							\draw[c1,thick] (0.545240931016100,0.773012684854573) node[left]{\tiny $18$};

							\draw[c1,thick] (0.726683377942314,0.732602609237550) node[left]{\tiny $19$};
							\draw[c1,thick] (0.920762287773320,0.695958574121097) node{\tiny $20$};
							\draw[c1,thick] (0.01,0.823728658042989) node[right]{\tiny $21$}; 
							\draw[c1,thick] (0.325708044929628,0.914763687652304) node{\tiny $22$};
							\draw[c1,thick] (0.500755877184029,0.910427640482701) node{\tiny $23$};
							\draw[c1,thick] (0.705821367824273,0.878409519462713) node{\tiny $24$};
							\draw[c1,thick] (0.899163331742718,0.898273740857340) node{\tiny $25$};
							\draw[c1,thick] (0.4,0.7) node[right]{\tiny $26$};
							\draw[c1,thick] (0.1,0.9) node[right]{\tiny $27$};
						\end{tikzpicture}
						\caption{Exact domain $\Omega$ with $27$ features.} \label{fig:geom81holes}
					\end{center}
				\end{subfigure}
				~
				\begin{subfigure}[t!]{0.315\textwidth}
					\begin{center}
						\begin{tikzpicture}[scale=4.8]
							\draw[thick] (0,0) -- (0,1) -- (1,1) -- (1,0) -- cycle;
							\draw[c1,thick] (0.097984328699318,0.092663690863534) circle (0.081318234590090);
							\draw[c1,thick] (0.283884293965162,0.123563491365580) circle (0.066396825383881);
							\draw[c1,thick] (0.546112220339815,0.056574896706227) circle (0.038919282299261);
							\draw[c1,thick] (0.715916527014347,0.093265457598229) circle (0.074000758141615);
							\draw[c1,thick] (0.898581968051930,0.104463986172862) circle (0.081763485970873);
							\draw[c1,thick] (0.067305141757723,0.339769841694328) circle (0.060034479317649);
							\draw[c1,thick] (0.495145556524009,0.307616854246403) circle (0.092235800109005);
							\draw[c1,thick] (0.500590433180793,0.509270838208646) circle (0.081283253301511);
							\draw[c1,thick] (0.110235870451937,0.693080838594426) circle (0.071264679475915);

							\draw[c1,thick] (0.097984328699318,0.092663690863534) node{\tiny $1$};
							\draw[c1,thick] (0.283884293965162,0.123563491365580) node{\tiny $2$};
							\draw[c1,thick] (0.546112220339815,0.056574896706227) node{\tiny $3$};
							\draw[c1,thick] (0.715916527014347,0.093265457598229) node{\tiny $4$};
							\draw[c1,thick] (0.898581968051930,0.104463986172862) node{\tiny $5$};
							\draw[c1,thick] (0.067305141757723,0.339769841694328) node{\tiny $6$};
							\draw[c1,thick] (0.495145556524009,0.307616854246403) node{\tiny $8$};
							\draw[c1,thick] (0.500590433180793,0.509270838208646) node{\tiny $13$};
							\draw[c1,thick] (0.110235870451937,0.693080838594426) node{\tiny $16$};
						\end{tikzpicture}
						\caption{Partially defeatured domain $\Omega^{(i)}_0$ at iteration $i=7$.} \label{fig:Omega08}
					\end{center}
				\end{subfigure}
				~
				\begin{subfigure}[t!]{0.345\textwidth}
					\begin{center}
						\includegraphics[scale=0.135]{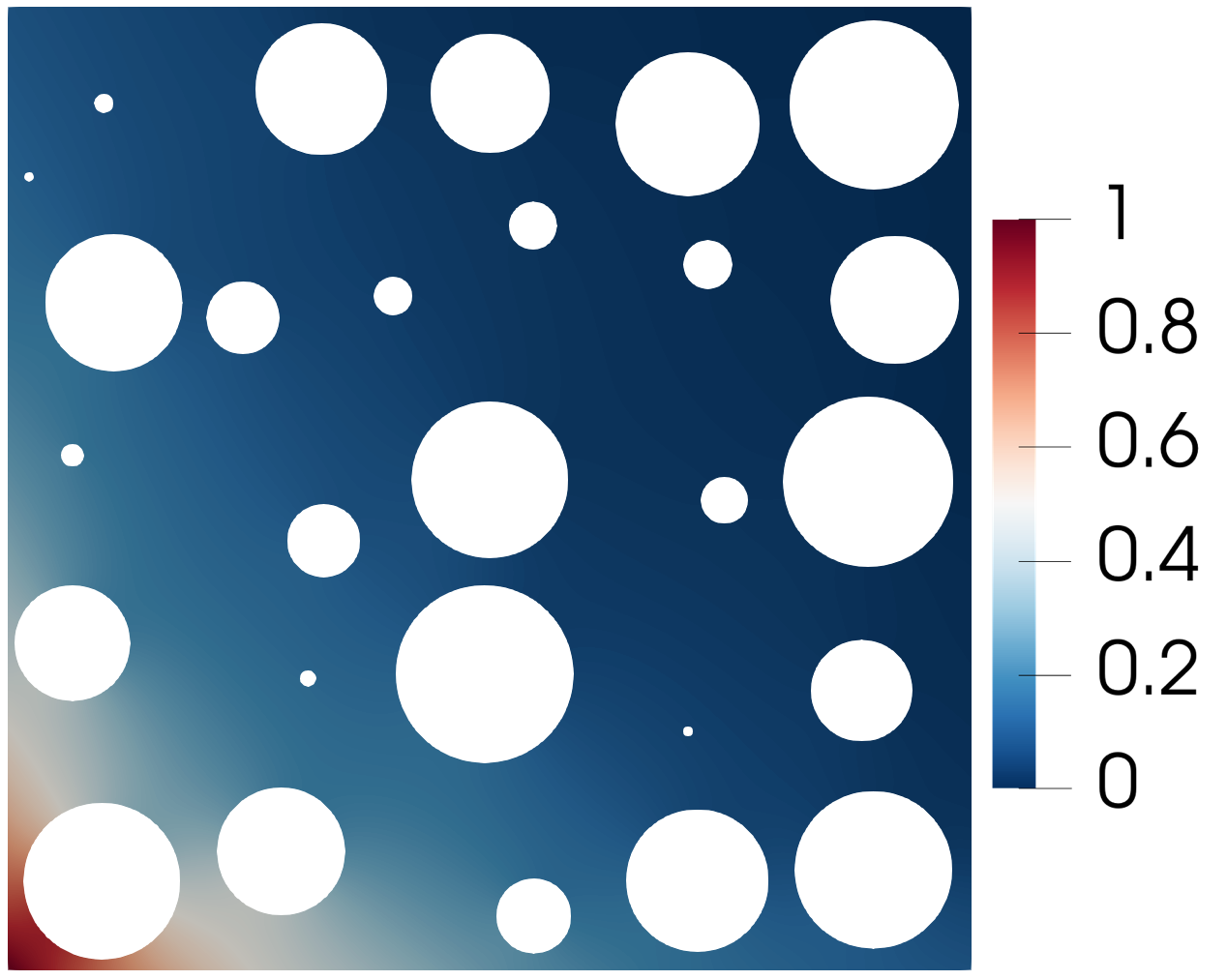}
						\caption{Exact solution in $\Omega$.\\{\color{white}.}} \label{fig:sol81holes}
					\end{center}
				\end{subfigure}
				\caption{Numerical test~\ref{sssmf:nbfeat} -- Geometry with $27$ features, considered exact solution, and corresponding partially defeatured geometrical model at iteration $i=7$.}
			\end{center}
			\end{figure}
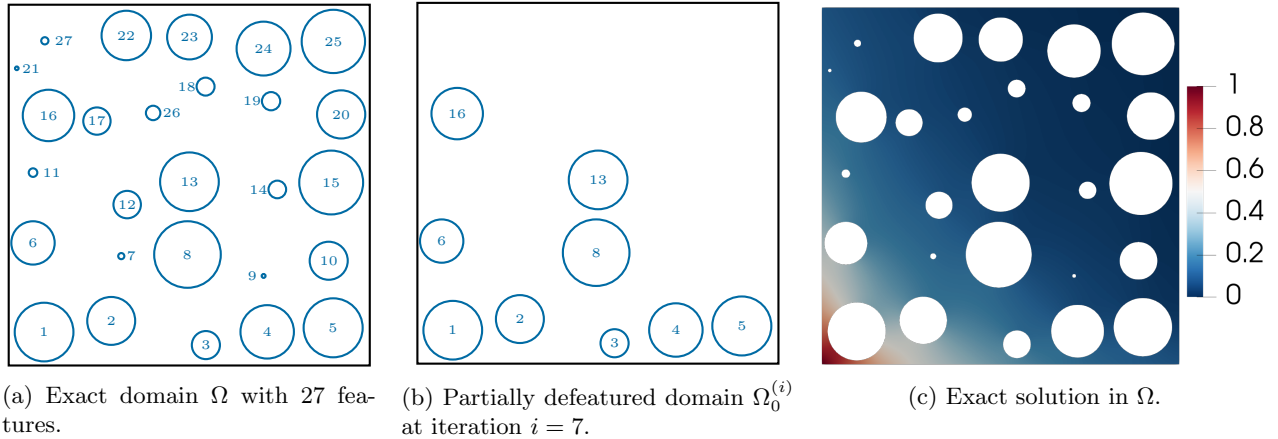
			\begin{figure}
			\centering
			\begin{tikzpicture}[scale=0.8]
				\begin{axis}[ymode=log, legend style={draw=none, row sep=5pt}, xlabel=Number of features added to $\Omega_0^{(i)}$,width=12cm,height=8cm, grid style=dotted,grid]
					\addplot[mark=o, c2, thick] table [x=km1, y=err, col sep=comma] {\dataPath/cvnbfeat.csv};
					\addplot[mark=x, c1, thick] table [x=km1, y=est, col sep=comma] {\dataPath/cvnbfeat.csv};
					\legend{$\vertiii{u-u^{(i)}_0}_\Omega$,$\mathscr E\Big(u_0^{(i)}\Big)$}; 
				\end{axis}
			\end{tikzpicture}
			\caption{Numerical test~\ref{sssmf:nbfeat} -- Behavior of the defeaturing error and estimator with respect to the number of features in the defeatured geometrical models $\Omega_0^{(i)}$. Each marker corresponds to the value at one iteration.} \label{fig:cvnbfeat}
		\end{figure}
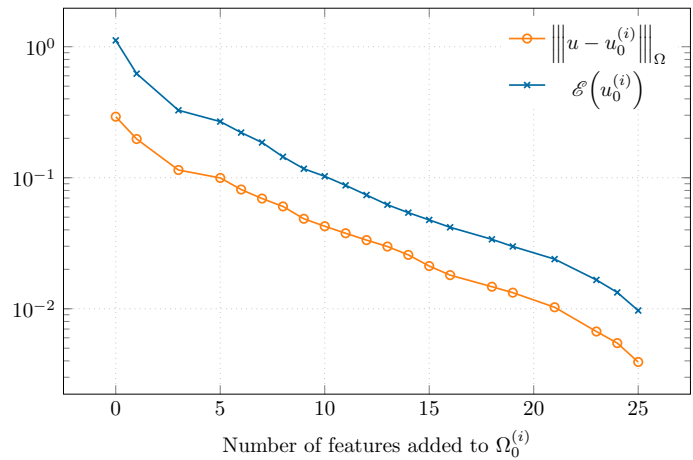}

The results are presented in Figure~\ref{fig:cvnbfeat}, and the sets of added features at each iteration are the following:
$\{1\}$, $\{2,6\}$, $\{4,16\}$, $\{8\}$, $\{3\}$, $\{5\}$, $\{13\}$, $\{12\}$, $\{17\}$, $\{22\}$, $\{11\}$, $\{10\}$, $\{15\}$, $\{23\}$, $\{24,16\}$, $\{7\}$, $\{20,27\}$, $\{18,25\}$, $\{21\}$, $\{14\}$, $\{19,9\}$. For instance, the error is divided by $10$ when $9$ out of the $27$ features are inserted in the partially defeatured geometrical model, i.e., a third of total number of features; this happens at iteration $i=7$, and $\Omega_0^{(7)}$ is represented in Figure~\ref{fig:Omega08}.
On the other hand, features $19$, $20$, $24$, and $25$, that are placed near the right top corner, where the solution is gradient is very small, are only activated in the last iterations.
We remark that the iteration index is directly linked to the number $N_f^{(i)}$ at each iteration $i$, that is, to the number of features that are still missing in the simplified geometrical model $\Omega_0^{(i)}$ with respect to the $27$ features in $\Omega$.
Moreover, we can see that the effectivity index is independent of the number of features that are not in the simplified geometrical model in which the problem is solved. Indeed, $\eta_{\mathrm{eff}}$ remains almost constant at each iteration, between $2.1$ and $3.8$. This result perfectly agrees with the theory developed in this paper, in particular the reliability and efficiency results of Theorems~\ref{thm:upperbound} and~\ref{thm:lowerbound}.
\subsubsection{Shape of the features} \label{ss:shapesstokes}
The presence of a feature in the computational domain may greatly but locally perturb the solution, for instance because of a sharp or even re-entrant corner. When such feature is removed, the defeatured solution becomes smoother. In this section, the presented numerical illustration demonstrates that the proposed estimator is able to capture the correct behavior of the error independently of the shape of the feature, and with a low effectivity index.

\begin{figure}[h!]
	\begin{center}
		\begin{subfigure}{\textwidth}
			\centering
			\includegraphics[scale=0.2,trim=575 1590 575 0, clip]{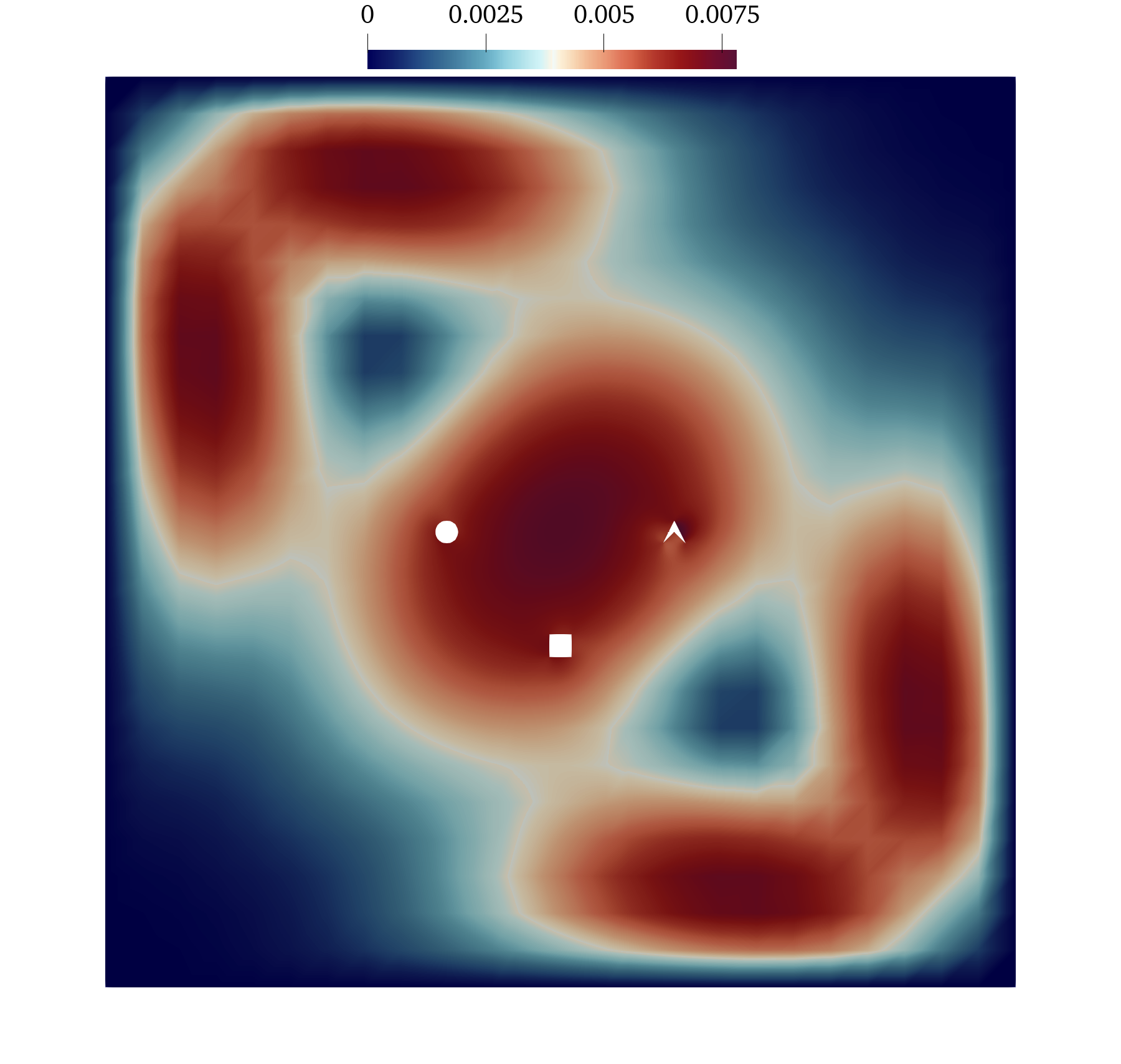}

			\includegraphics[scale=0.12,trim=125 125 125 125, clip]{images/shapes_vex_withsolution.png}
			\includegraphics[scale=0.12,trim=125 125 125 125, clip]{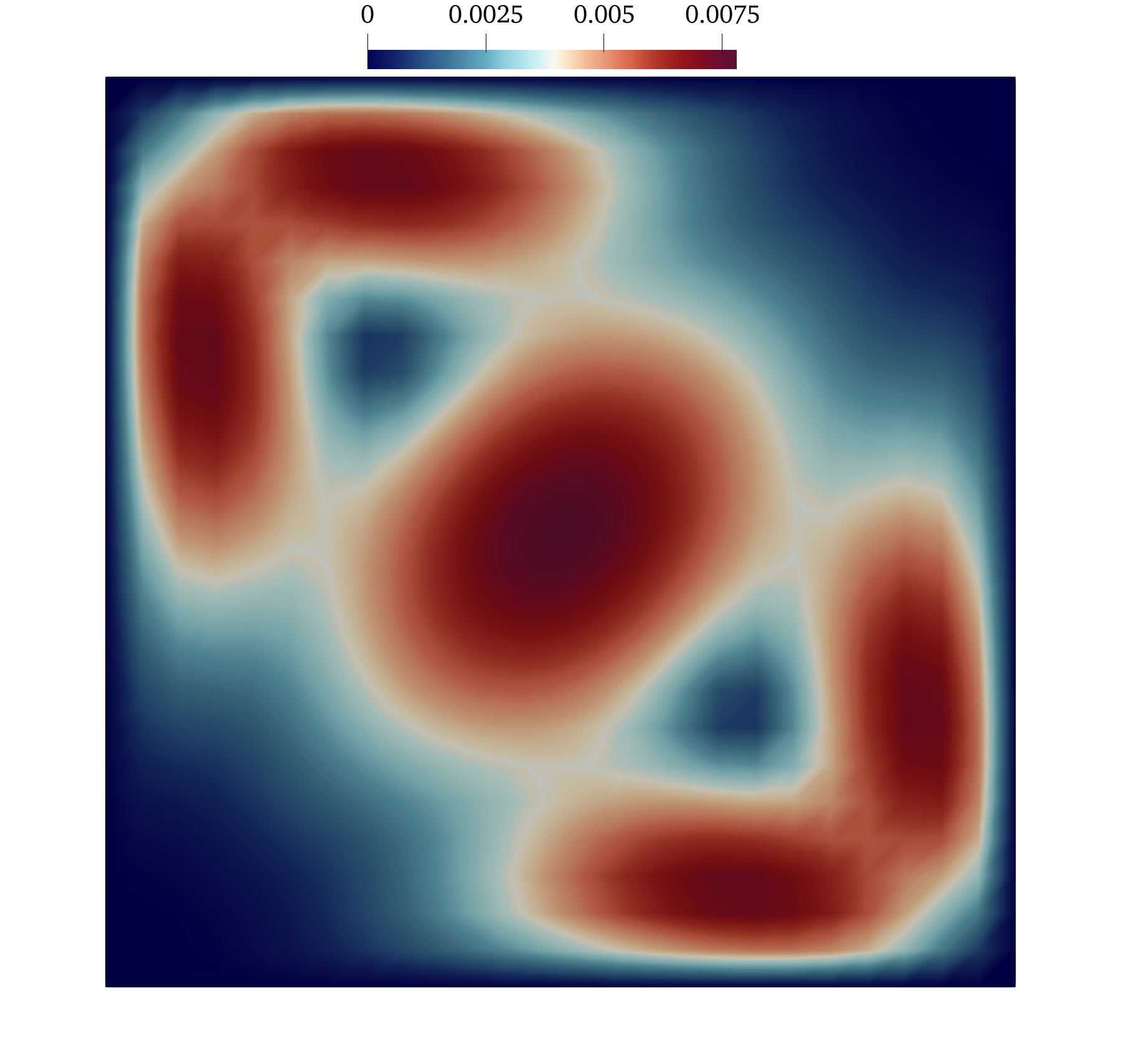}
			\caption{Magnitude of $\boldsymbol u$ (left) and  $\boldsymbol u_\mathrm d$ (right) in the domains $\Omega$ and $\Omega_0$, respectively.}
			\label{f:geomshapesstokes}
		\end{subfigure}
		~
		\begin{subfigure}{\textwidth}
			\centering
			\includegraphics[scale=0.2,trim=575 1590 575 0, clip]{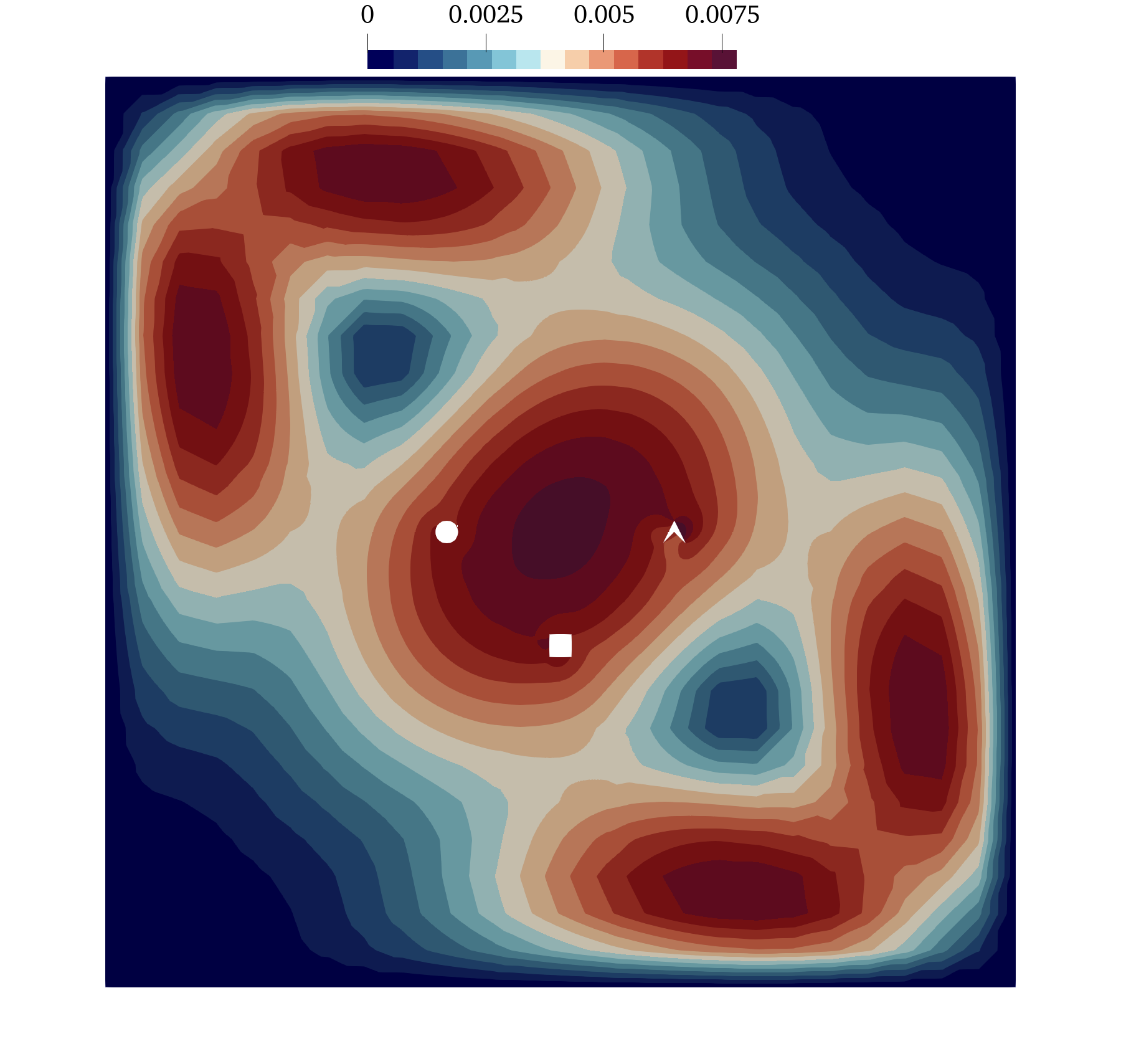}

			\includegraphics[scale=0.3,trim=600 600 600 600, clip]{images/shapes_vex_withsolution_coarsecolors.png}
			\includegraphics[scale=0.3,trim=600 600 600 600, clip]{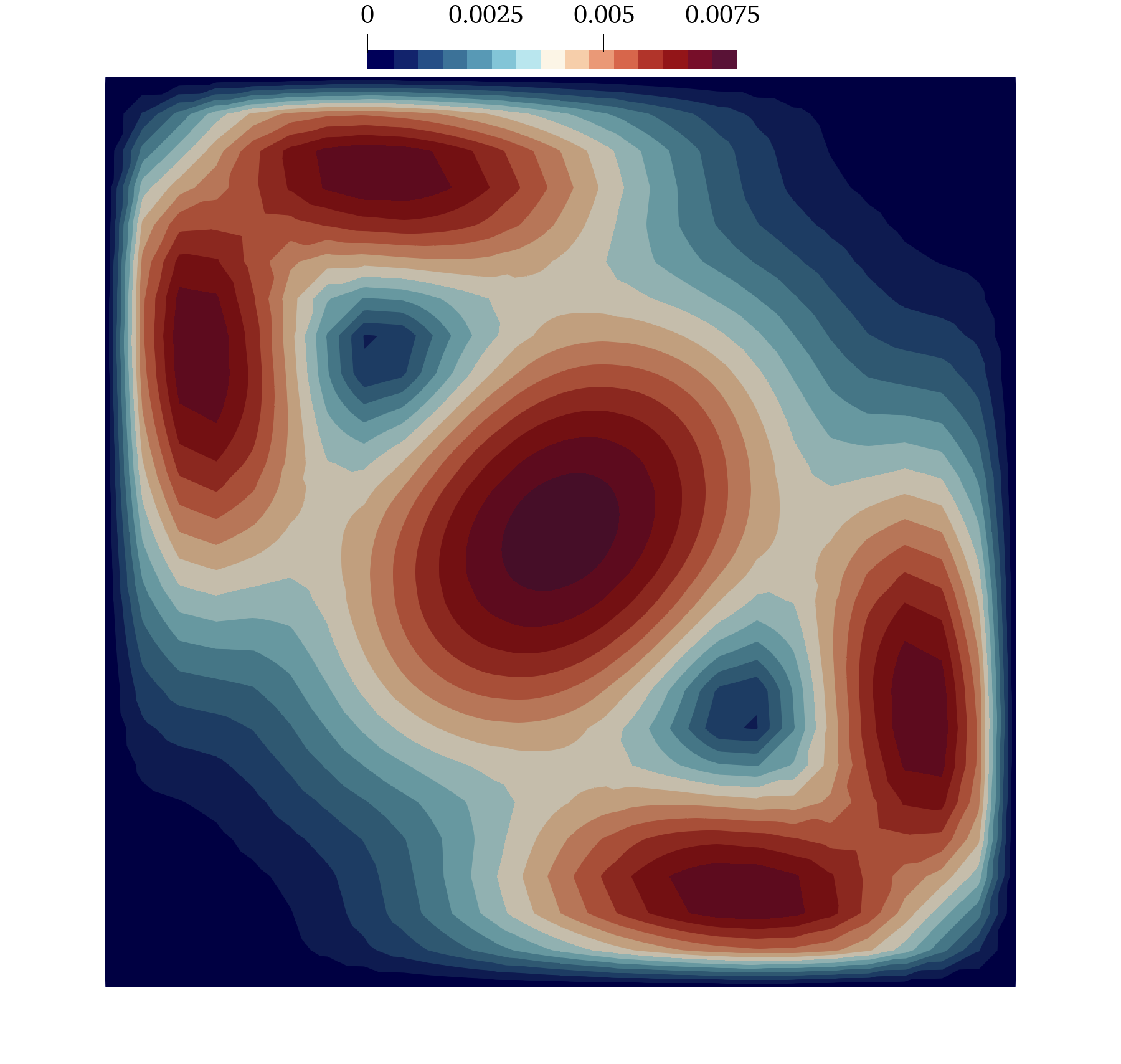}
			\caption{Zoom of the central region: Magnitude of $\boldsymbol u$ (left) and  $\boldsymbol u_\mathrm d$ (right) in the domains $\Omega$ and $\Omega_0$, respectively.}
			\label{f:geomshapesstokes_zoom}
		\end{subfigure}
		\caption{Numerical test~\ref{ss:shapesstokes} -- Magnitude of the velocity in the exact $\Omega$ and fully defeatured $\Omega_0$ domains.} \label{fig:displacementshapesstokes}
	\end{center}
\end{figure}

To do so, let $\Omega_0:=(0,1)^2$ be the fully defeatured domain, and let $\Omega:= \Omega_0\setminus\bigcup_{k=1}^3 \overline{F^k}$ be the exact computational domain containing three holes: $F^1$ is a circle of radius $0.0125$ centered at $(0.375,0.5)$, $F^2$ is a square of side $0.0250$ centered at $(0.5,0.375)$, and $F^3$ is a non-convex quadrilateral creating three re-entrant corners in $\Omega$, whose vertices are placed in $(0.625,0.5125)$, $(0.6125, 0.4875)$, $(0.625, 0.5)$ and $(0.6375,0.4875)$. See Figure~\ref{f:geomshapesstokes} (left) for an illustration of the domain $\Omega$.

In this experiment, we aim at finding the exact solution $(\boldsymbol u, p)$ of the Stokes' problem~(\ref{eqid:weakoriginalstokespb}) defined in $\Omega$, where $\boldsymbol f(x,y):=\exp\big[4\big( (x-0.5)^2+(y-0.5)^2 \big)\big]$ and $f_c\equiv 0$ in $\Omega$, $\boldsymbol g_D \equiv 0$ in $\partial \Omega_0$, and $\boldsymbol g\equiv 0$ in $\partial F^k$ for all $k=1,2,3$. Instead of solving problem~(\ref{eqid:weakoriginalstokespb}), we tackle the approximate one~(\ref{eqid:simplstokespb}) in $\Omega_0$ with the same data; in particular, $\boldsymbol f$ is naturally extended in the features, and we obtain the defeatured solution $(\boldsymbol{u}_0,p_0)$.

The magnitude of the velocity fields $\boldsymbol{u}$ and $\boldsymbol{u}_\mathrm d$ is shown in Figure~\ref{fig:displacementshapesstokes}.
As it can be seen in the zoom-in (Figure~\ref{f:geomshapesstokes_zoom}), in this case, the sharp non-convex feature $F^3$ introduces a localized perturbation that is similar to ones produced by the sharp convex $F^1$ and smooth $F^2$ features, as all of them are located at solution regions with similar gradient magnitudes.
Indeed, their contributions $\mathscr{E}^k(\boldsymbol u_\mathrm d, p_\mathrm d)$, for $k=1,2,3$, reported in Table~\ref{tbl:featurecontribshapes}, are of the same order.
From those individual contributions, through~\eqref{eq:estforeachfeatstokes}, the total defeaturing error estimate is computed to be $\mathscr{E}(\boldsymbol u_\mathrm d, p_\mathrm d)= 1.0895\cdot 10^{-1}$, and the corresponding energy error is $\vertiii{(\boldsymbol u - \boldsymbol u_\mathrm d, p-p_\mathrm d)}_\Omega = 7.5381\cdot 10^{-2}$. Therefore, the effectivity index of the proposed estimator is
$$\eta_{\mathrm{eff}} := \frac{\mathscr E(\boldsymbol{u}_\mathrm d, p_\mathrm d)}{\vertiii{(\boldsymbol u - \boldsymbol u_\mathrm d, p-p_\mathrm d)}_\Omega} = 1.4454,$$
which is notably very low (at the same level as the ones observed in \cite{paper1defeaturing} for geometries with a single feature): The estimator is able to estimate the effect of different features, independently of their shape.

\begin{table}
	\setlength{\tabcolsep}{4pt}
	\centering
	{\def\arraystretch{1.2}
		\begin{tabular}{@{}cccc@{}}
			\hline
			Feature index $k$                                    & $1$                   & $2$                   & $3$                   \\\hline
			$\mathscr E^k(\boldsymbol u_\mathrm d, p_\mathrm d)$ & $2.1713\cdot 10^{-2}$ & $2.4578\cdot 10^{-2}$ & $1.5601\cdot 10^{-2}$ \\
			\hline
		\end{tabular}}
	\caption{Numerical test \ref{ss:shapesstokes} -- Feature contributions $\mathscr E^k(\boldsymbol u_\mathrm d, p_\mathrm d)$ to the multi-feature estimator $\mathscr E(\boldsymbol u_\mathrm d, p_\mathrm d)$.} \label{tbl:featurecontribshapes}
\end{table}

\subsection{Lid-driven cavity} \label{ss:lidcavity}
For this next numerical experiment, let us consider Stokes' problem in a lid-driven cavity~\cite{kuhlmann2019lid} in which three holes are located. More precisely, we consider an exact domain $\Omega = \Omega_0\setminus (\overline{F^1}\cup \overline{F^2}\cup \overline{F^3})$ for which $\Omega_0 := (0,1)^2$ is the fully defeatured domain, and $F^1$, $F^2$, and $F^3$ are three circular holes of radius $0.01$ centered, respectively, at $(0.011,0.989)$, $(0.5,0.75)$, and $(0.011,0.011)$. For this test, we let $\Gamma_D:= \partial \Omega_0$, $\Gamma_N:= \partial F^1 \cup \partial F^2\cup \partial F^3$, $\boldsymbol{f}\equiv\boldsymbol{0}$ in $\Omega_0$, $\boldsymbol g \equiv \boldsymbol{0}$ on $\Gamma_N$, and
$$\boldsymbol{g}_D \equiv \begin{cases}
		(1,0) & \text{ on } (0,1)\times \{1\} \text{ the top boundary,} \\
		(0,0) & \text{ everywhere else on }\Gamma_D.
	\end{cases}$$
Then, let $(\boldsymbol{u}, p)$ be the solution of the exact Stokes' problem~(\ref{eqid:weakoriginalstokespb}), and let $(\boldsymbol{u}_\mathrm d, p_\mathrm d)\equiv(\boldsymbol{u}_0, p_0)$ be the solution of the corresponding Stokes' problem~(\ref{eqid:simplstokespb}) in the fully defeatured geometry $\Omega_0$. We compute the estimator $\mathscr E(\boldsymbol{u_\mathrm d}, p_\mathrm d)$ defined in~\eqref{eq:estforeachfeatstokes} by computing each feature contribution $\mathscr E^k(\boldsymbol{u}_\mathrm d, p_\mathrm d)$ for $k=1,2,3$.

\begin{figure}
	\begin{center}
		\begin{subfigure}{\textwidth}
			\centering
			\includegraphics[scale=0.2,trim=575 1590 575 0, clip]{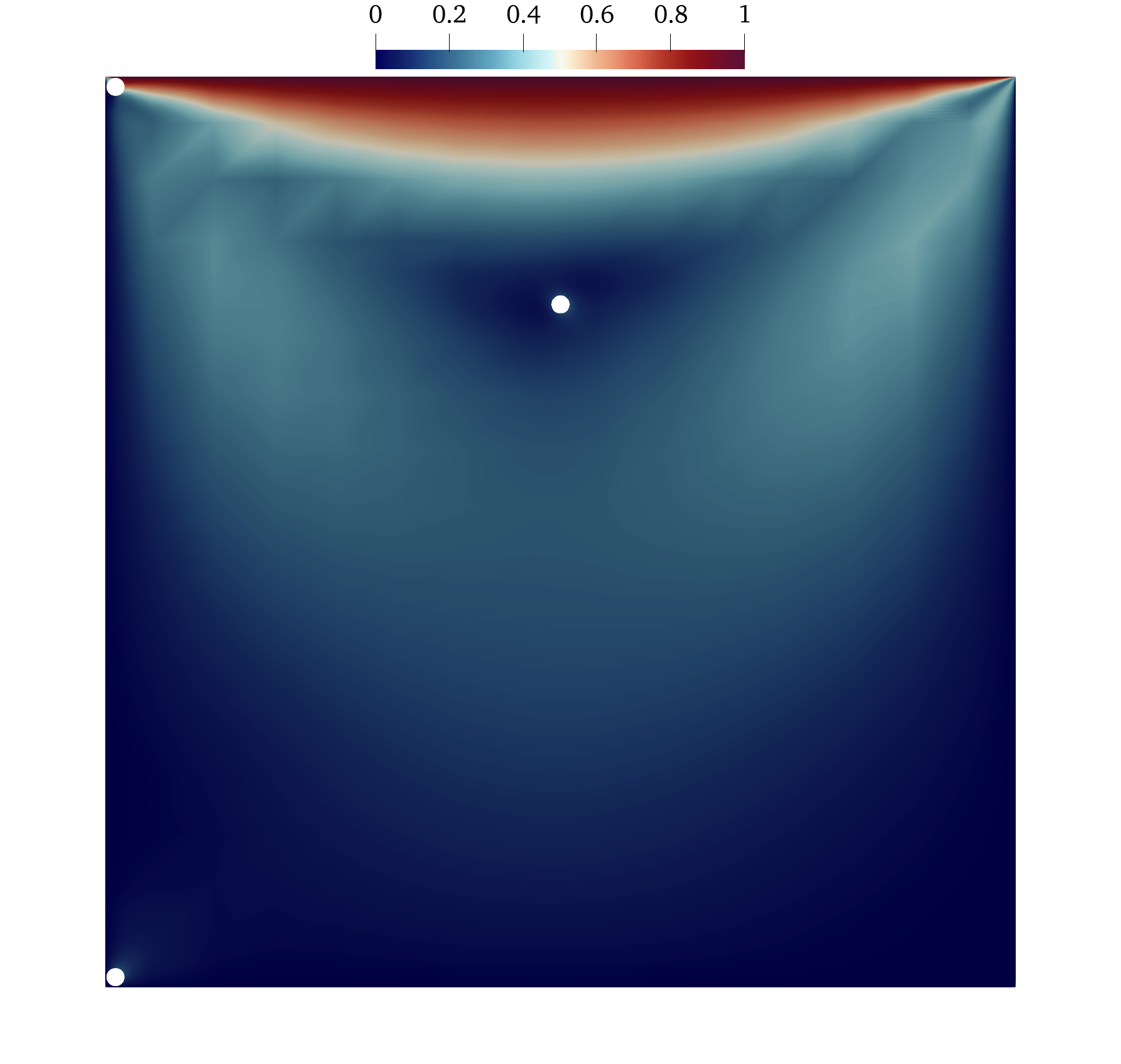}

			\includegraphics[scale=0.15,trim=125 125 125 125, clip]{images/lidcavity_vex_withsolution.png}
			\includegraphics[scale=0.15,trim=125 125 125 125, clip]{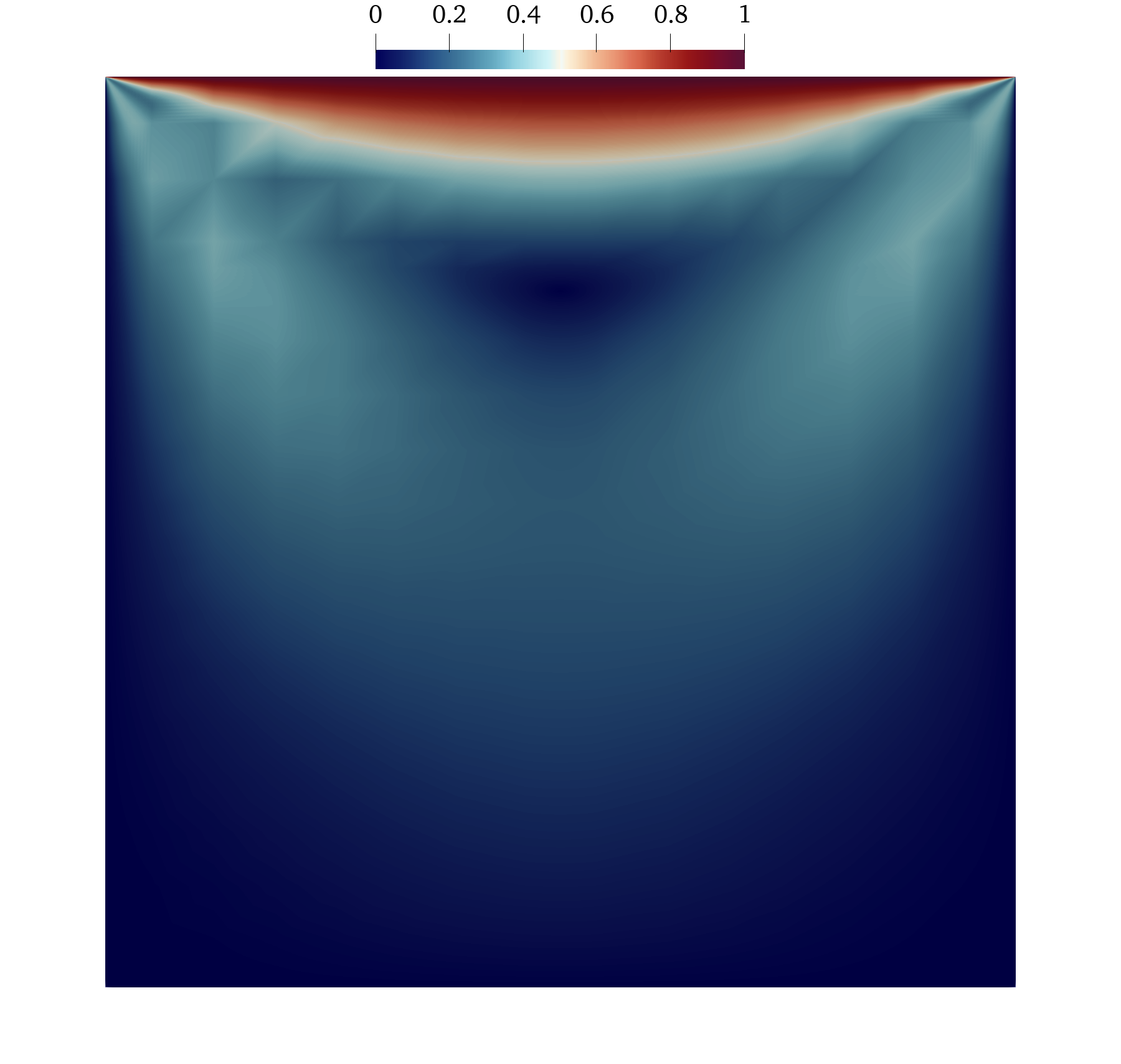}
			\caption{Magnitude of $\boldsymbol u$ (left) and  $\boldsymbol u_0$ (right) in the domains $\Omega$ and $\Omega_0$, respectively.}
			\label{f:lidcavity}
		\end{subfigure}
		~
		\begin{subfigure}{0.48\textwidth}
			\centering
			\includegraphics[scale=0.2,trim=575 1590 575 0, clip]{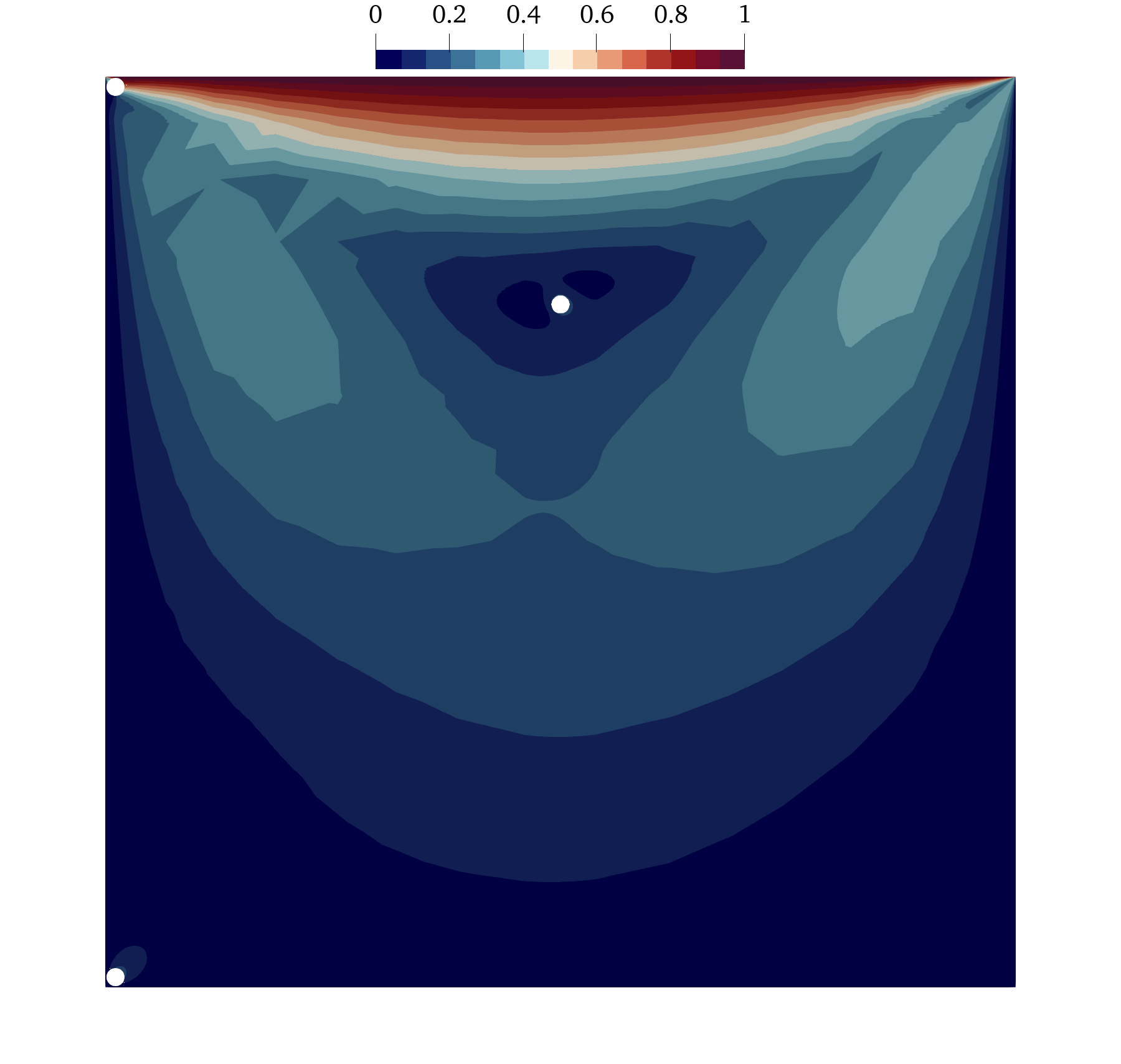}
			\includegraphics[scale=0.06,trim=0 0 0 0, clip]{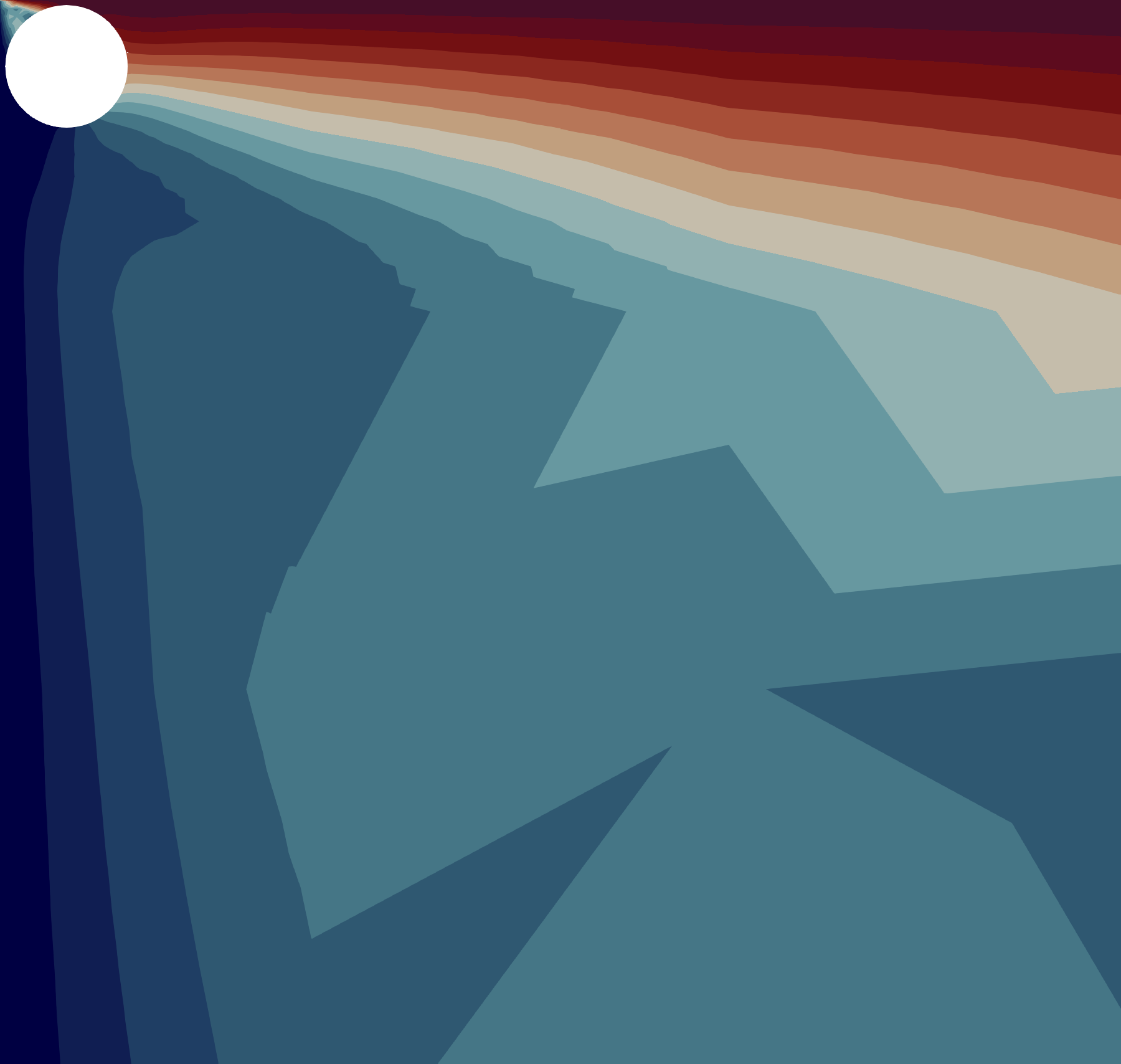}
			\includegraphics[scale=0.06,trim=0 0 0 0, clip]{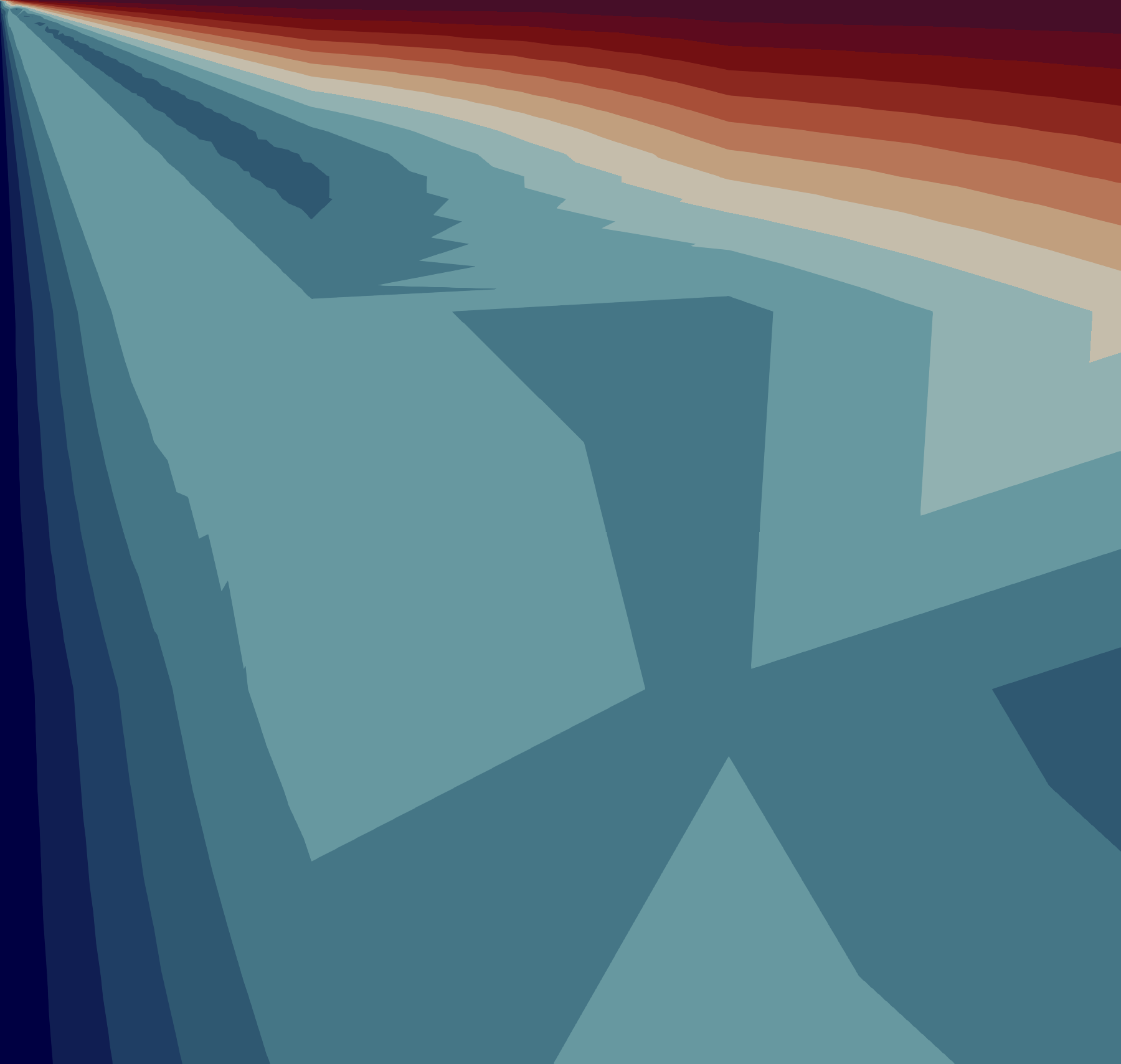}
			\caption{Magnitude of $\boldsymbol u$ (left) and $\boldsymbol u_0$ (right) around $F^1$.}
			\label{f:lidcavity_f1}
		\end{subfigure}
		~
		\begin{subfigure}{0.48\textwidth}
			\centering
			\includegraphics[scale=0.2,trim=575 1590 575 0, clip]{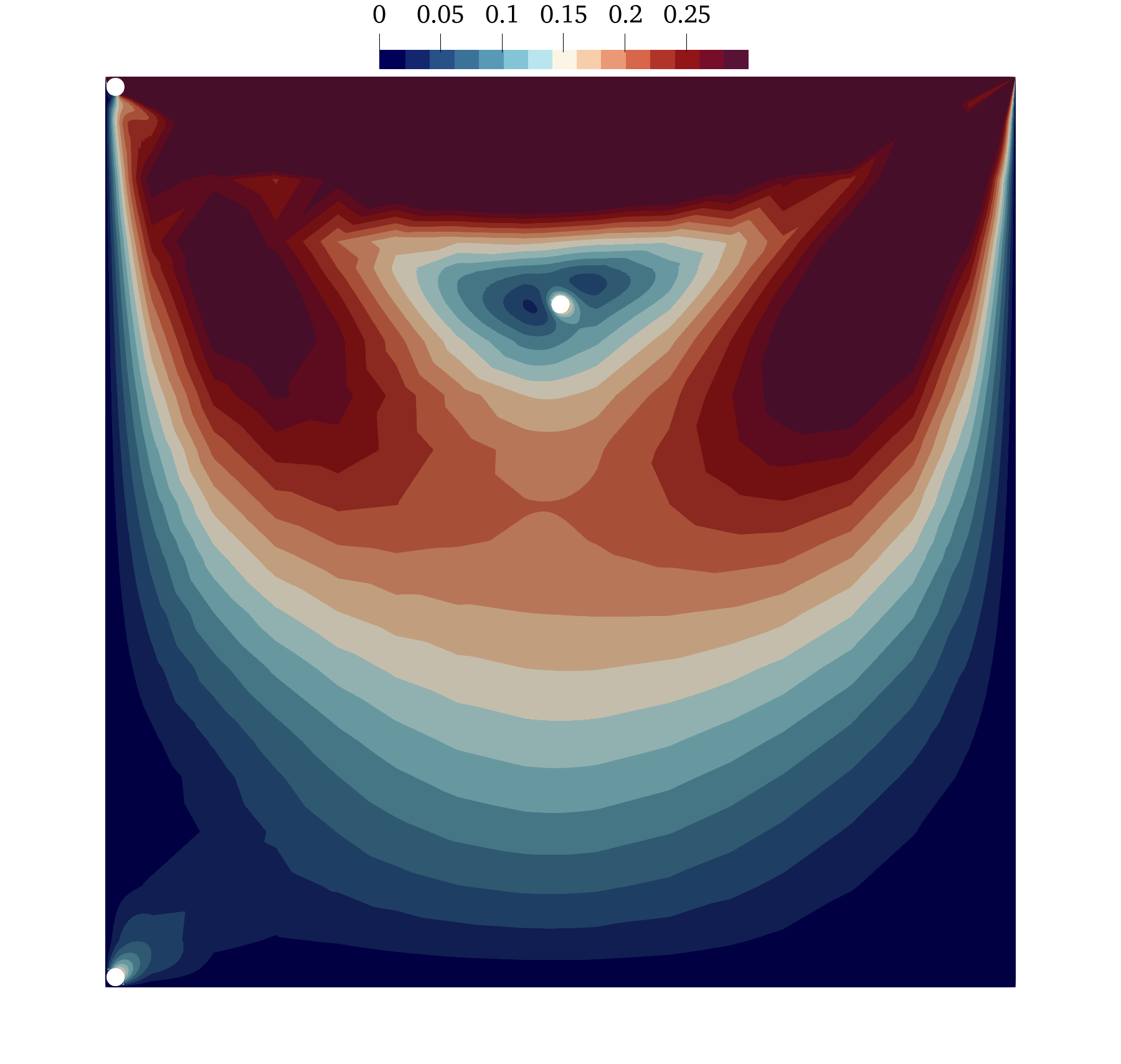}
			\includegraphics[scale=0.06,trim=0 0 0 0, clip]{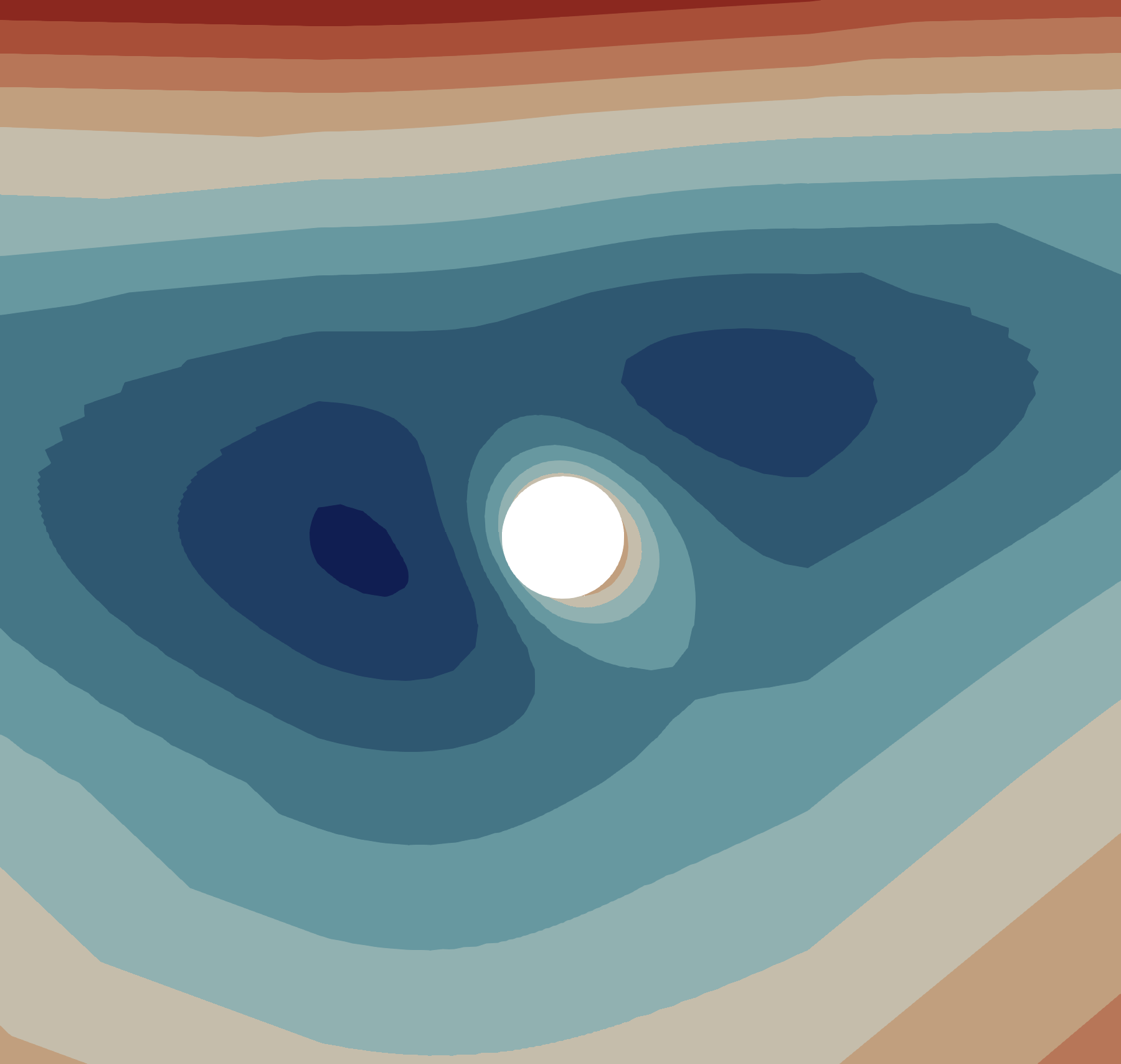}
			\includegraphics[scale=0.06,trim=0 0 0 0, clip]{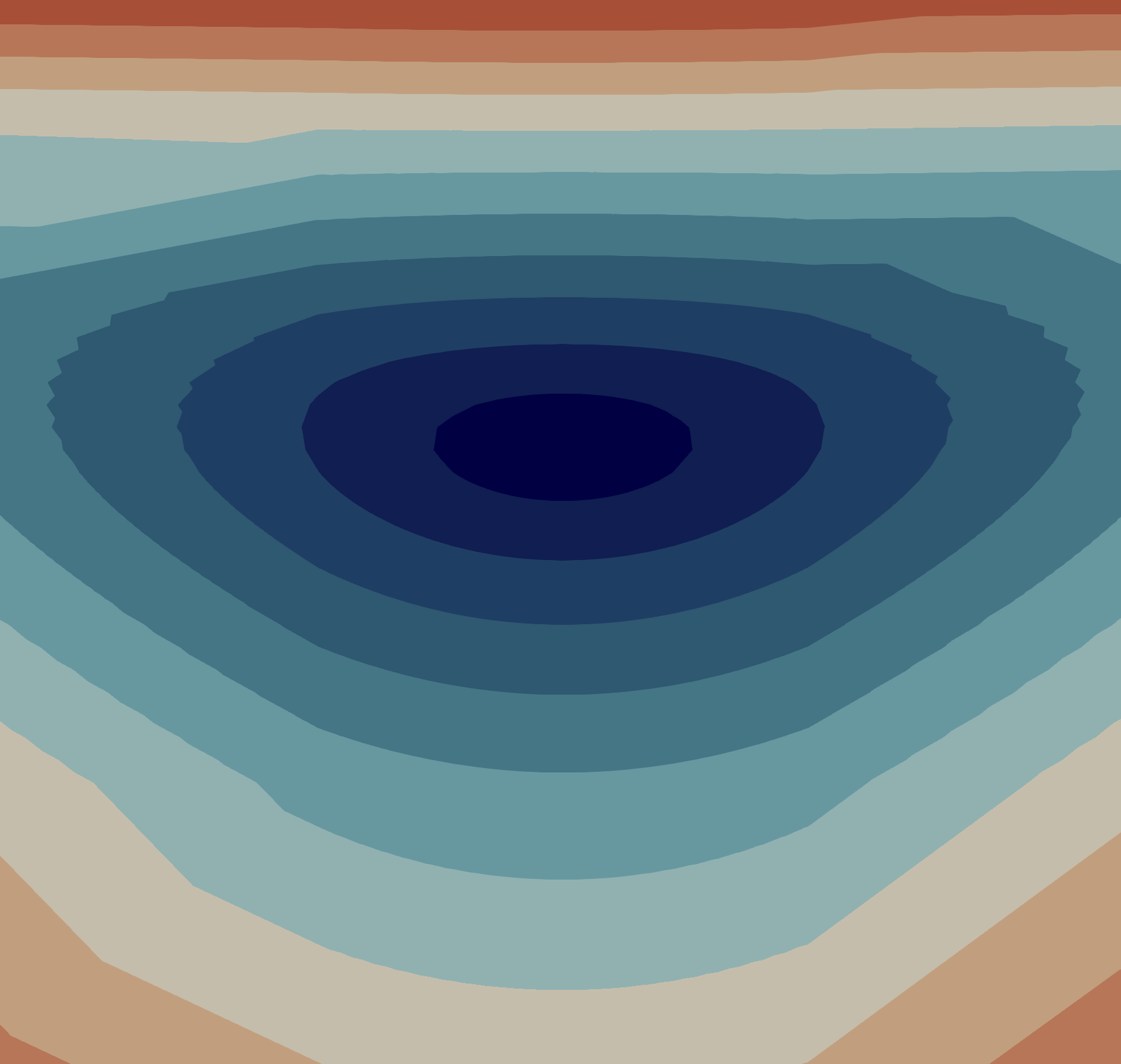}
			\caption{Magnitude of $\boldsymbol u$ (left) and $\boldsymbol u_0$ (right) around $F^2$.}
			\label{f:lidcavity_f2}
		\end{subfigure}
		~
		\begin{subfigure}{0.48\textwidth}
			\centering
			\includegraphics[scale=0.2,trim=575 1590 575 0, clip]{images/lidcavity_vex_withsolution_coarsecolors1.png}
			\includegraphics[scale=0.06,trim=0 0 0 0, clip]{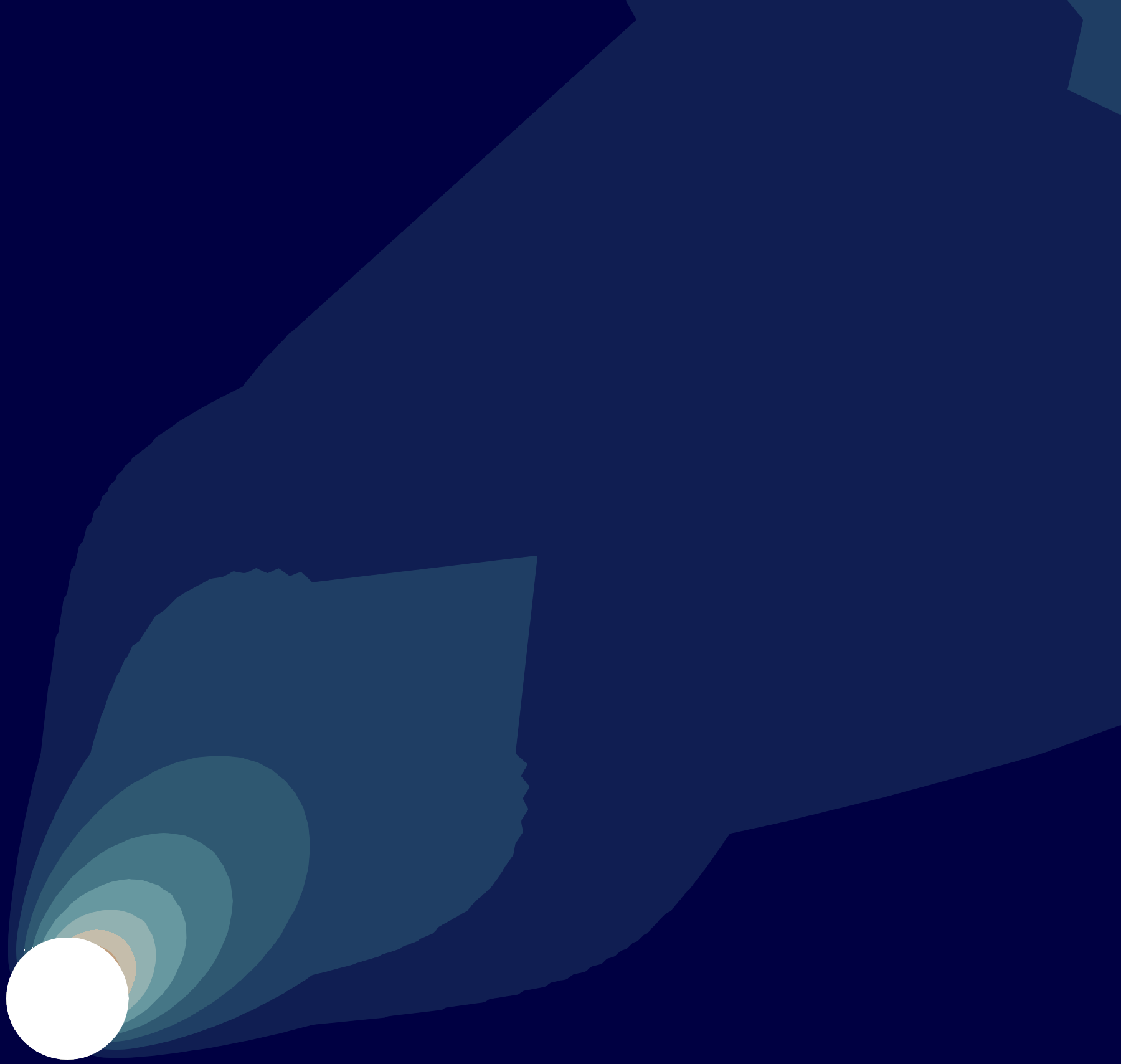}
			\includegraphics[scale=0.06,trim=0 0 0 0, clip]{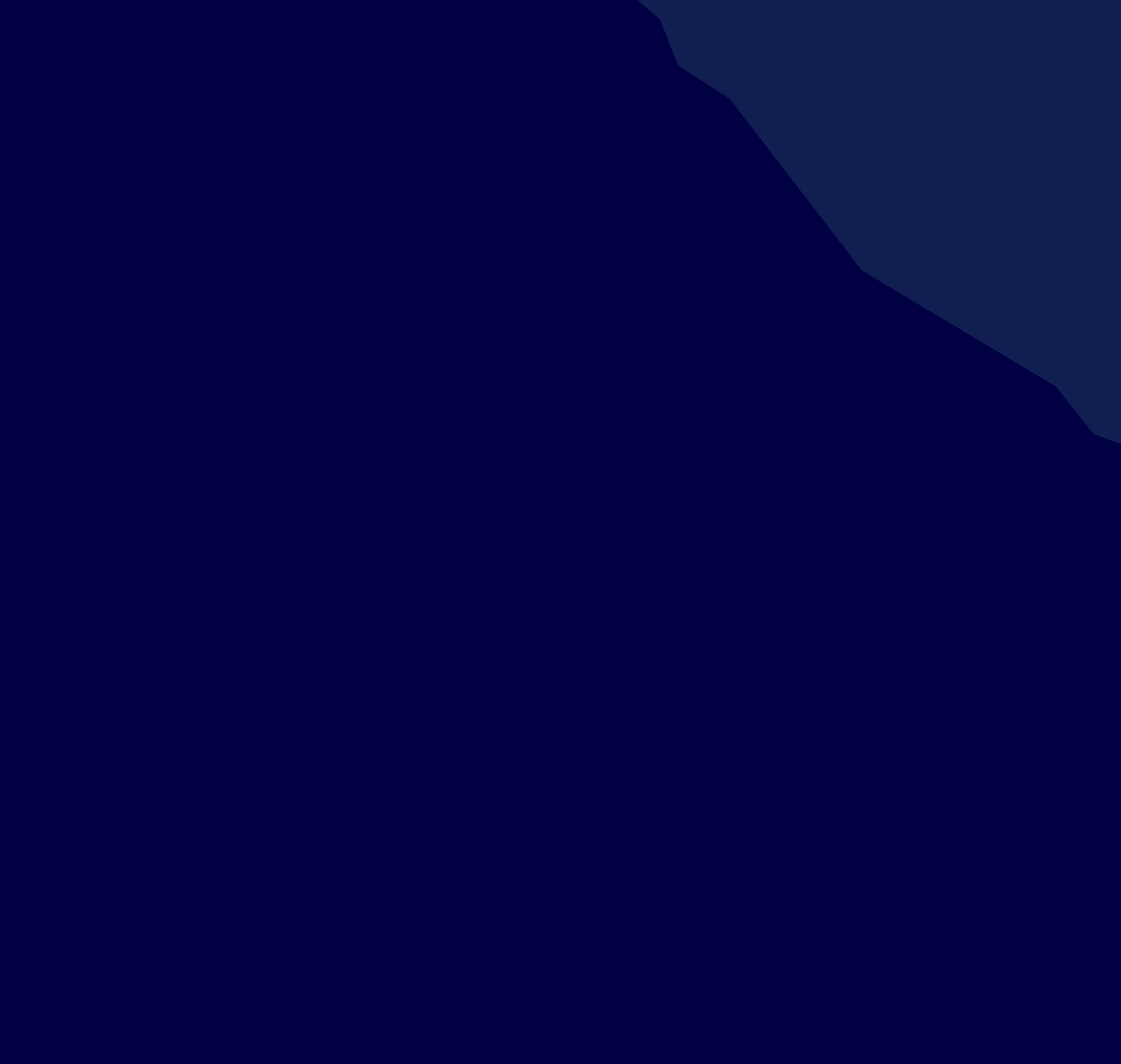}
			\caption{Magnitude of $\boldsymbol u$ (left) and $\boldsymbol u_0$ (right) around $F^3$.}
			\label{f:lidcavity_f3}
		\end{subfigure}
		\caption{Numerical test~\ref{ss:lidcavity} -- Magnitude of the velocity in the exact domain $\Omega$ and in the corresponding fully defeatured domain $\Omega_0$.} \label{fig:displacementlidcavity}
	\end{center}
\end{figure}

Results are presented in Figure~\ref{fig:displacementlidcavity}, in which the magnitude of the velocity fields $\boldsymbol{u}$ and $\boldsymbol{u}_0$ is shown,
and Table~\ref{tbl:featurecontribldc} reports each feature contribution's $\mathscr{E}^k(\boldsymbol u_\mathrm d, p_\mathrm d)$ for $k=1,2,3$.
We observe that the presence of the features changes the fluid velocity inside the cavity, especially for features $F^2$ and $F^3$. This is expected because of the homogeneous Neumann boundary conditions on the features.
However, features only change locally the velocity field of the fluid around them. Thus, feature $F^1$ brings the largest contribution to the estimator while feature $F^3$ brings the smallest one. This is coherent with the theory since the velocity has a large gradient close to the moving boundary at the top, while it is almost constantly equal to zero around features $F^2$ and $F^3$ (note the different scales in Figures~\ref{f:lidcavity_f1}-\ref{f:lidcavity_f3}). We therefore expect $F^1$ to contribute the most to the overall defeaturing error, even if at a first glance it is the hole whose absence changes the least the velocity of the fluid (recall Figure~\ref{f:lidcavity_f1}).

\begin{figure}
	\centering
	\includegraphics[scale=0.3,trim=0 80 0 50, clip]{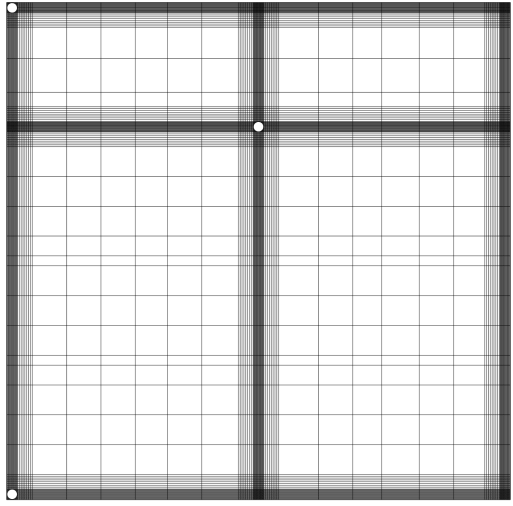}
	\caption{Numerical test \ref{ss:lidcavity} -- Mesh used to compute an overkilled solution of the lid driven cavity Stokes' problem in the exact domain $\Omega$.} \label{fig:meshldc}
\end{figure}

The total defeaturing error estimate is equal to $\mathscr{E}(\boldsymbol u_\mathrm d, p_\mathrm d)= 35.622$, and $\vertiii{(\boldsymbol u - \boldsymbol u_\mathrm d, p-p_\mathrm d)}_\Omega = 31.308$. Therefore, the effectivity index of the proposed estimator is equal to
$$\eta_{\mathrm{eff}} := \frac{\mathscr E(\boldsymbol{u}_\mathrm d, p_\mathrm d)}{\vertiii{(\boldsymbol u - \boldsymbol u_\mathrm d, p-p_\mathrm d)}_\Omega} = 1.1378,$$
which is notably very low, as in the previous experiment.
Note that to compute the defeaturing error, a very fine mesh around the holes had to be taken in $\Omega$ in order to be able to neglect the component of the error coming from the numerical approximation of the problem, as represented in Figure~\ref{fig:meshldc}.
This is not required when the holes are filled as in the defeatured geometry $\Omega_0$: This shows the potential of defeaturing, in terms of memory and computational time savings.
Far from the features and high solution gradients, a coarser grid is considered.

\begin{table}
	\setlength{\tabcolsep}{4pt}
	\centering
	{\def\arraystretch{1.2}
		\begin{tabular}{@{}cccc@{}}
			\hline
			Feature index $k$                                    & $1$                & $2$                   & $3$                   \\\hline
			$\mathscr E^k(\boldsymbol u_\mathrm d, p_\mathrm d)$ & $1.1874\cdot 10^1$ & $6.7065\cdot 10^{-2}$ & $2.1130\cdot 10^{-2}$ \\
			\hline
		\end{tabular}}
	\caption{Numerical test \ref{ss:lidcavity} -- Feature contributions $\mathscr E^k(\boldsymbol u_\mathrm d, p_\mathrm d)$ to the multi-feature estimator $\mathscr E(\boldsymbol u_\mathrm d, p_\mathrm d)$.} \label{tbl:featurecontribldc}
\end{table}

\subsection{Three-dimensional elastic structure} \label{sec:3Dtest}
\begin{figure}
	\begin{subfigure}[t]{0.9\textwidth}
		\centering
		\includegraphics[scale=0.465]{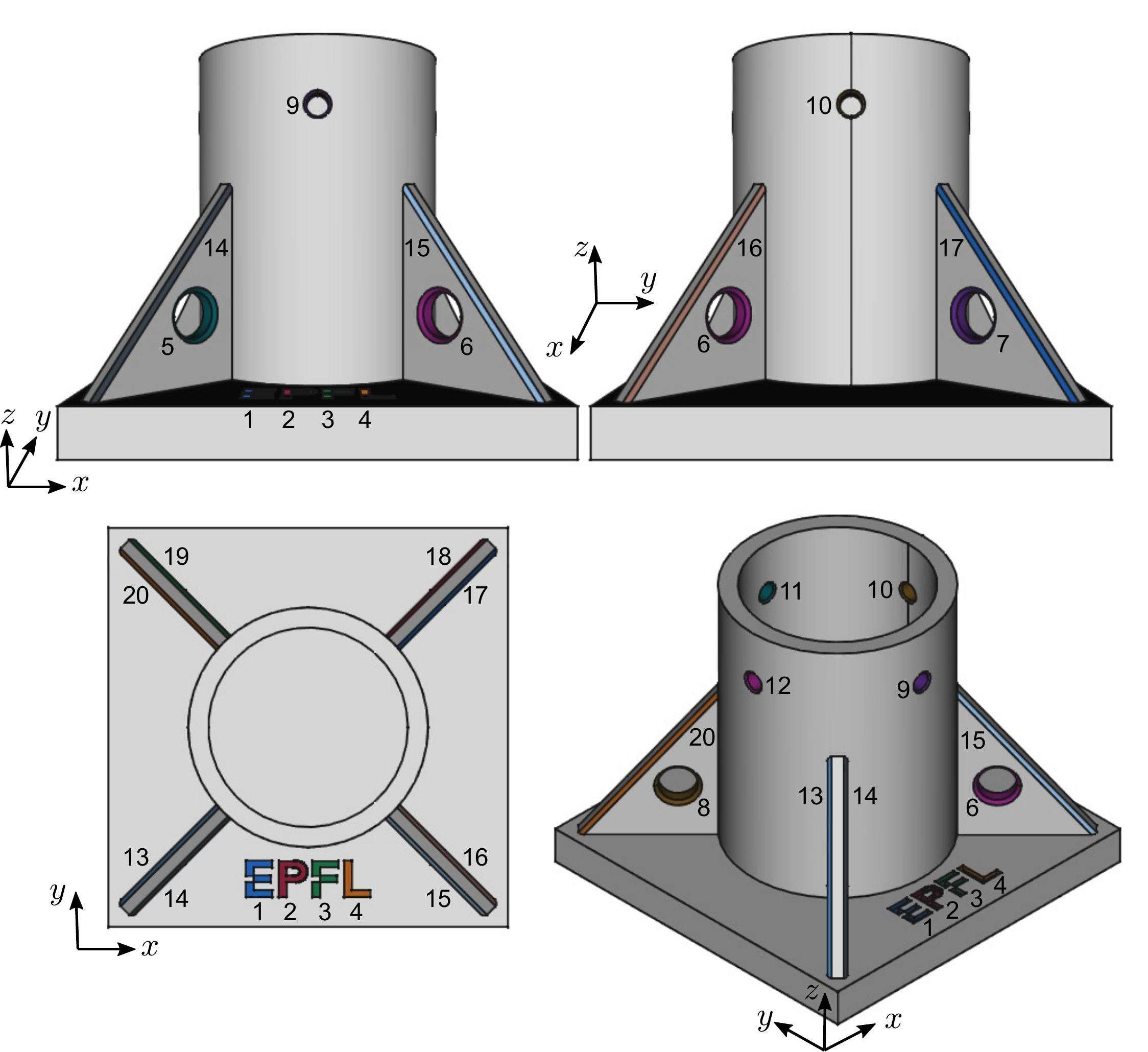}
		\caption{Exact geometry $\Omega$ and numbering of the $20$ features (in color).} \label{fig:exact3D}
	\end{subfigure}
	~
	\begin{subfigure}[t]{0.9\textwidth}
		\centering
		\includegraphics[scale=0.465]{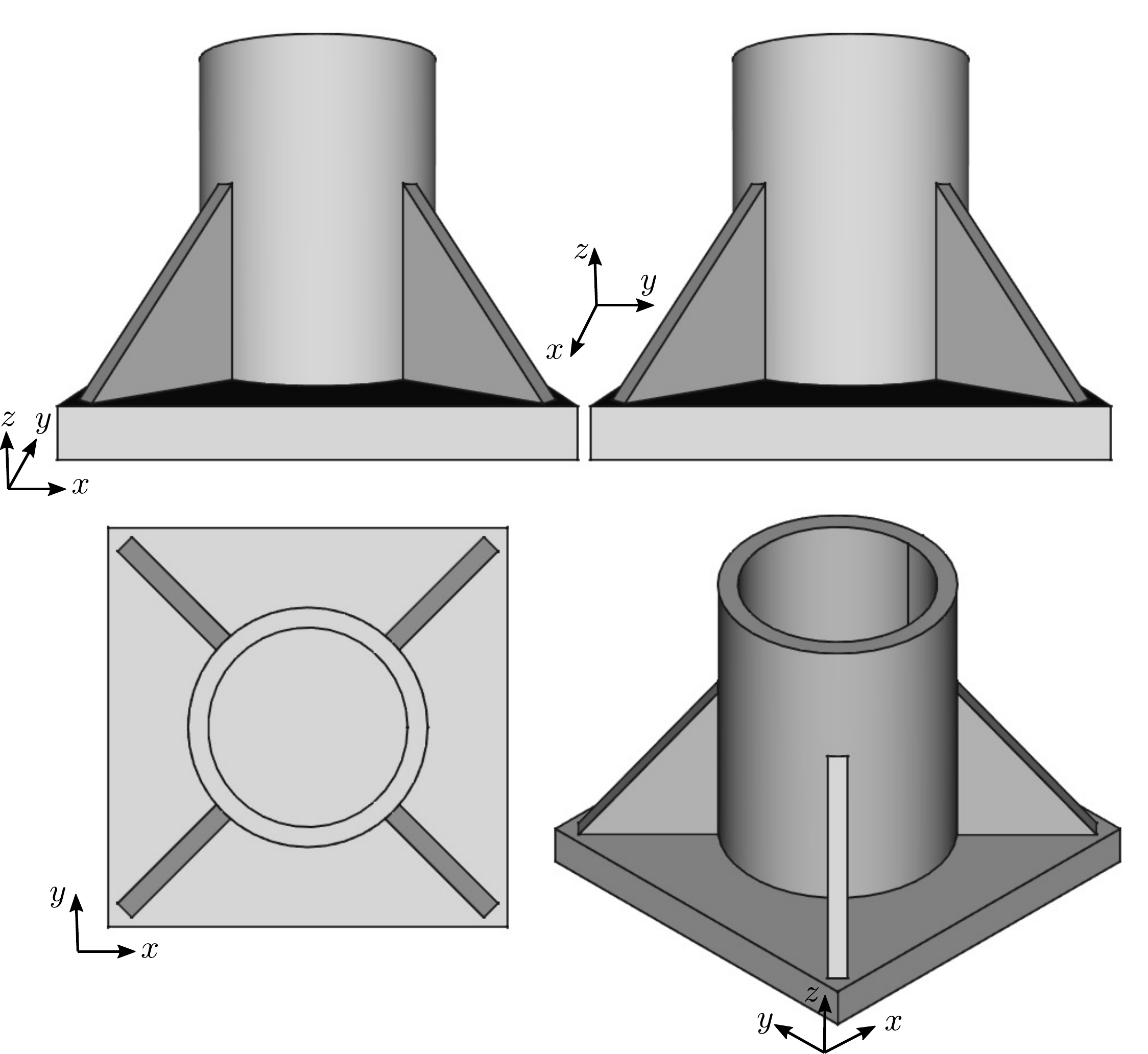}
		\caption{Defeatured geometry $\Omega_0$}\label{fig:defeat3D}
	\end{subfigure}
	\caption{Numerical test \ref{sec:3Dtest} -- Exact and defeatured $3$D domains; the colored boundaries correspond to $\gamma^k$, for each feature $k=1,\ldots,20$ as numbered in (\ref{sub@fig:exact3D}).} \label{fig:structures}
\end{figure}

For this last numerical experiment, let us consider the exact domain $\Omega$ and the corresponding defeatured domain $\Omega_0$ represented in Figure~\ref{fig:structures}. More precisely, the base has dimensions $200\times200\times20$ [mm], and the cylinder has a height of $150$ [mm]. Moreover, and in particular, the exact domain contains $20$ features numbered as illustrated in Figure~\ref{fig:exact3D}:
\begin{itemize}
	\item $F^1$ to $F^4$ are the four letters of the carved ``EPFL'' logo, in order (see also Figure~\ref{figmf:defeatadaptive} in which these features are more clearly visible).
	\item $F^5$ to $F^8$ are the four holes in the stiffeners, counted counter-clockwise beginning from the one on the left of the ``EPFL'' logo.
	\item $F^9$ to $F^{12}$ are the four holes in the vertical part of the structure, counted counter-clockwise beginning from the one above the ``EPFL'' logo.
	\item $F^{13}$ to $F^{20}$ are the eight rounds present on the left and right diagonal angles of the stiffeners, counted counter-clockwise beginning from the left round of the stiffener on the left of the ``EPFL'' logo.
\end{itemize}
Rounds, holes, and carved logos are three of the most typical features that finite element analysis practitioners encounter in CAD designs. These features are interesting to analyze, since they are usual candidates to be removed before creating a finite element mesh.

Taking the origin at the bottom lower left corner of the structure, let $\Gamma_D$ be the bottom of the structure and $\Gamma_N := \partial \Omega \setminus \overline{\Gamma_D}$, let $\boldsymbol{f} = \boldsymbol{0}$ [N$\cdot$ mm$^{-3}$], $\boldsymbol{g}_D = \boldsymbol{0}$ [mm], and
$$\boldsymbol{g} \,\mathrm{[MPa]} = \begin{cases} \boldsymbol{0}                    & \mathrm{on } \,\Gamma_N\setminus \overline{\Gamma_{\mathrm{top}}} \\
              \boldsymbol{e}_x = (1, 0, 0)^\top & \mathrm{on } \,\Gamma_{\mathrm{top}},\end{cases}$$
where $\Gamma_{\mathrm{top}}$ is the top face of the cylinder. Then,  let $\boldsymbol u\in \boldsymbol H_{\boldsymbol 0,\Gamma_D}^1(\Omega)$ be the solution of the linear elasticity problem given \changes{by~(\ref{eqid:weakoriginalstokespb}) in which the pressure terms and the divergence condition are removed, and} where the material properties correspond to steel. That is, the Lam\'e parameters $\lambda$ and $\mu$ are expressed in terms of the Young modulus $E=210$ [GPa] and Poisson's ration $\nu = 0.3$ [--] as
\begin{equation*}
	\lambda = \frac{E\nu}{(1+\nu)(1-2\nu)} \quad \text{ and } \quad \mu = \frac{E}{2(1+\nu)}.
\end{equation*}
Now, let us extend $\boldsymbol f$ by $\boldsymbol 0$ in all features so that $\boldsymbol f=\boldsymbol{0}$ [N$\cdot$ m$^{-3}$] in $\Omega_0$, and let $\boldsymbol{g}_0 = \boldsymbol{0}$ [Pa] on $\gamma_0 := \partial \Omega_0\setminus \overline{\partial \Omega}$. Then we compute the defeatured solution $\boldsymbol u_\mathrm d\equiv\boldsymbol u_0\in \boldsymbol H^1_{\boldsymbol 0, \Gamma_D}(\Omega_0)$ given by \changes{problem~(\ref{eq:weaksimplstokespb}) in which again the pressure terms and the divergence condition are removed}. Finally, we compute the estimator $\mathscr{E}(\boldsymbol u_\mathrm d)$ defined \changes{in~(\ref{eq:estforeachfeatstokes}) with $p_\mathrm d\equiv 0$} by computing each feature contribution $\mathscr{E}^k(\boldsymbol u_\mathrm d)$ for $k=1,\ldots,20$.

\begin{figure}
	\centering
	\def\svgwidth{0.8\linewidth}
	\input{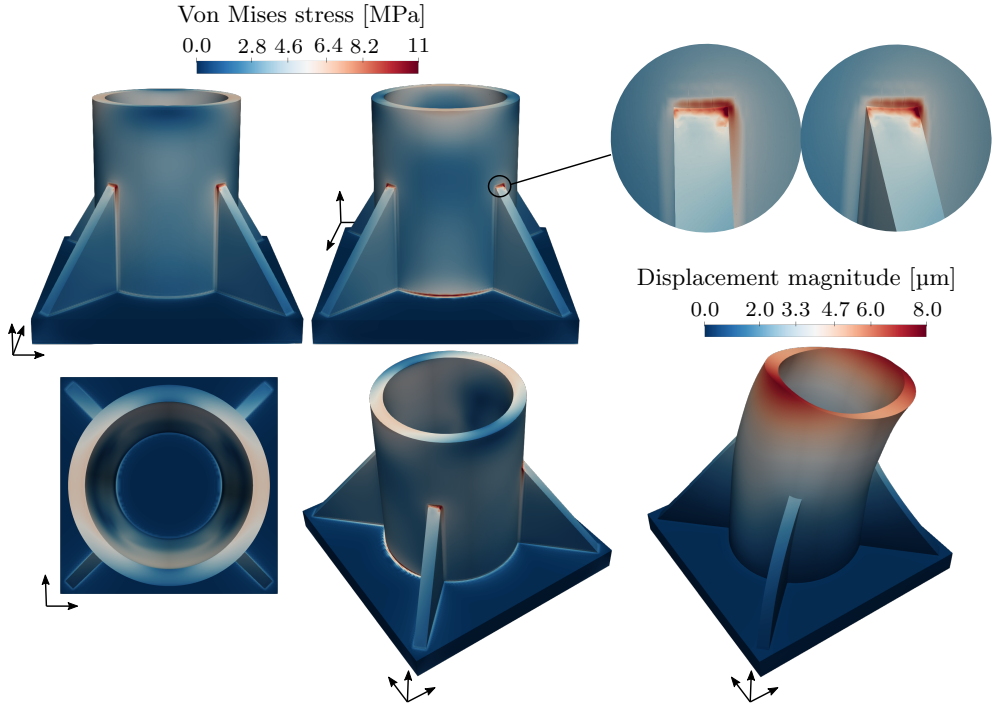}
	\caption{Numerical test \ref{sec:3Dtest} -- Defeatured solution in the defeatured geometry $\Omega_0$. The views correspond to those of Figure \ref{fig:structures}, and the deformed configuration is magnified $[\times 5\cdot10^3]$ for visualization purposes.} \label{fig:defeatsolution3D}
\end{figure}

A rather fine mesh is used in order to reduce the error derived from numerical approximation. More precisely, the bounding box of $\Omega_0$ is meshed with $n_{\mathrm{el}}=128$ elements per direction, and B-splines of degree $2$ and regularity $1$ are used. Results are presented in Table~\ref{tbl:featurecontrib}, where we report each feature's contribution $\mathscr{E}^k(\boldsymbol u_\mathrm d)$ for $k=1,\ldots,20$. The obtained total error estimator is equal to $\mathscr{E}(\boldsymbol u_\mathrm d)= 1.716\cdot 10^{-6}$[J]. Moreover, the magnitude of the solution displacements $\boldsymbol{u}$ and $\boldsymbol{u}_\mathrm d$, and the corresponding von Mises stress distributions are shown in Figures~\ref{fig:defeatsolution3D} and~\ref{fig:exactsolution3D}, respectively. 

\begin{table}
	\setlength{\tabcolsep}{4pt}
	\centering
	{\def\arraystretch{1.2}
		\begin{tabular}{@{}ccccccccccc@{}}
			\hline
			Feature index $k$                                          & $1$     & $2$     & $3$     & $4$     & $5$     & $6$     & $7$     & $8$     & $9$     & $10$    \\
			$\mathscr E^k(\boldsymbol u_\mathrm d)$ [$\cdot 10^{-8}$J] & $1.949$ & $2.904$ & $3.032$ & $1.278$ & $69.28$ & $69.28$ & $69.28$ & $69.24$ & $61.96$ & $26.06$ \\
			\hline
			Feature index $k$                                          & $11$    & $12$    & $13$    & $14$    & $15$    & $16$    & $17$    & $18$    & $19$    & $20$    \\
			$\mathscr E^k(\boldsymbol u_\mathrm d)$ [$\cdot 10^{-8}$J] & $61.96$ & $26.06$ & $8.797$ & $14.84$ & $14.69$ & $8.871$ & $8.980$ & $14.74$ & $14.75$ & $9.001$ \\
			\hline
		\end{tabular}}
	\caption{Numerical test \ref{sec:3Dtest} -- Feature contributions $\mathscr E^k(\boldsymbol u_\mathrm d)$ to the multi-feature estimator $\mathscr E(\boldsymbol u_\mathrm d)$.} \label{tbl:featurecontrib}
\end{table}

We can first see that the absence of features $F^5$ to $F^{8}$ in the defeatured geometry significantly affects the solution in the stiffeners. This is indeed reflected in the estimator: The estimator contributions of those four features is very large, corresponding to around half of the total error estimator. On the other hand, the solution is basically constant around the ``EPFL'' logo, no deformation is observed around it. We can therefore expect that the absence of features $F^1$ to $F^4$ in the defeatured geometry is not affecting much the accuracy of the solution. This is indeed observed in the estimator contributions of those features, as $\mathscr{E}^1(\boldsymbol u_\mathrm d)$ to $\mathscr{E}^4(\boldsymbol u_\mathrm d)$ are the lowest contributions of the estimator, corresponding to around $1\%-2\%$ of $\mathscr{E}(\boldsymbol u_\mathrm d)$. This is a typical situation that simulation practitioners encounter daily: Carved logos and trademarks are usually defeatured before creating a finite element mesh, since they complicate the meshing process and increase the number of elements (see, e.g., Figure~\ref{fig:complexmesh_org}), but they contribute little to the accuracy of the problem's solution. The proposed estimator identifies them straightaway.

{\begin{figure}
	\centering
	\def\svgwidth{0.8\linewidth}
	\input{\dataPath/iteration10_4.pdf_tex}
	\caption{Numerical test \ref{sec:3Dtest} -- Exact solution in the exact geometry $\Omega$. The views correspond to those of Figure \ref{fig:structures}, and the deformed configuration is magnified $[\times 5\cdot10^3]$ for visualization purposes.} \label{fig:exactsolution3D}
	\centering
	\def\svgwidth{0.8\linewidth}
	\input{\dataPath/iteration5_4.pdf_tex}
	\caption{Numerical test \ref{sec:3Dtest} -- Partially defeatured solution in the partially defeatured geometry obtained at iteration $4$. The letters of the ``EPFL'' logo and the four rounds $F^{13}$, $F^{16}$, $F^{17}$, and $F^{20}$ are the missing features in this geometry. The views correspond to those of Figure \ref{fig:structures}, and the deformed configuration is magnified $[\times 5\cdot10^3]$ for visualization purposes.} \label{fig:intermediatesolution3D}
\end{figure}}

Let us now run the adaptive algorithm introduced in Section~\ref{s:adaptivity} with $\theta = 0.99$ as marking parameter, until all features are added to the geometrical model. We call $\boldsymbol u_\mathrm d^{(i)}$ the solution of the defeatured problem at iteration $i$.
In Table \ref{tbl:adaptivedefeatresults}, we report the indices of the features that are added to the defeatured geometrical model at each iteration, together with the value of the estimator $\mathscr{E}\big(\boldsymbol{u}_\mathrm d^{(i)}\big)$. The magnitude of the solution displacement and the corresponding von Mises stress distribution at iteration $4$ are represented in Figure \ref{fig:intermediatesolution3D}. Comparing the values of the estimator at each iteration and the von Mises stress distributions around each feature, we can see that the features that are added to the geometrical model at each iteration seem to be the ones that are affecting the most the solution accuracy, as one would expect.
In Table \ref{tbl:adaptivedefeatresults} we can also see that to reduce the error estimator by $90\%$, it is enough to consider $12$ out of the total $20$ features of $\Omega$ (see iteration $4$, whose solution is represented in Figure \ref{fig:intermediatesolution3D}).

For instance, the holes $F^9$ and $F^{11}$ are added before the holes $F^{10}$ and $F^{12}$ during the adaptive process, because of the direction in which the structure is bending due to the applied traction along the $x$-direction; this is reflected by the variation of the von Mises stresses that are larger in $F^9$ and $F^{11}$ than in $F^{10}$ and $F^{12}$.
We can also see that larger stresses are present nearer the rounds $F^{14}$, $F^{15}$, $F^{18}$, and $F^{19}$ than around the other four rounds. This is again coming from the direction of the bending. And very interestingly, the estimator is able to capture this effect, as rounds $F^{14}$, $F^{15}$, $F^{18}$, and $F^{19}$ are introduced in the defeatured geometry after iteration $3$, while the other rounds are introduced later, after iterations $4$ and $5$.
See, for instance, the stress distribution in the connection between features $F^{17}$ and $F^{18}$ and the main cylinder (zoom-in regions in Figures~\ref{fig:defeatsolution3D}, \ref{fig:exactsolution3D}, and \ref{fig:intermediatesolution3D}).
As it can be appreciated, the stress concentration is higher around feature $F^{18}$, fact that also reveals the estimator value in Table~\ref{tbl:featurecontrib}, and the fact that feature $F^{18}$ is activated before than $F^{17}$ (see Table~\ref{tbl:addedfeatsdiffnel} and Figure~\ref{fig:intermediatesolution3D}).

Rounds are other typical examples of features that are candidates to be removed. However, in this case, the situation is usually less clear. Indeed, \review{on the one} hand, rounds complicate the meshing process and increase the number of elements in the model. But on the other hand, depending on the boundary conditions, removing rounds may lead to the creation of singularities in the solution. The proposed estimator is able to determine the impact of removing those rounds.

\begin{table}
	\centering
	{\def\arraystretch{1.2}
		\begin{tabular}{@{}ccccccccccc@{}}
			\hline
			Iteration $i$                                                            & $0$       & $1$     & $2$     & $3$           & $4$     \\
			\hline
			Marked features                                                          & $5,6,7,8$ & $9,11$  & $10,12$ & $14,15,18,19$ & $17,20$ \\
			$\mathscr{E}\big(\boldsymbol{u}_\mathrm d^{(i)}\big)$ [$\cdot 10^{-8}$J] & $171.6$   & $99.75$ & $49.19$ & $32.08$       & $16.64$ \\[3pt]
			\hline\hline
			Iteration $i$                                                            & $5$       & $6$     & $7$     & $8$           & $9$     \\
			\hline
			Marked features                                                          & $13,16$   & $3$     & $2$     & $1$           & $4$     \\
			$\mathscr{E}\big(\boldsymbol{u}_\mathrm d^{(i)}\big)$ [$\cdot 10^{-8}$J] & $12.21$   & $5.056$ & $3.958$ & $2.516$       & $1.345$ \\[3pt]
			\hline
	\end{tabular}}
	\caption{Numerical test \ref{sec:3Dtest} -- Results of the adaptive defeaturing strategy for $n_\mathrm{el}=128$ elements.} \label{tbl:adaptivedefeatresults}
\end{table}

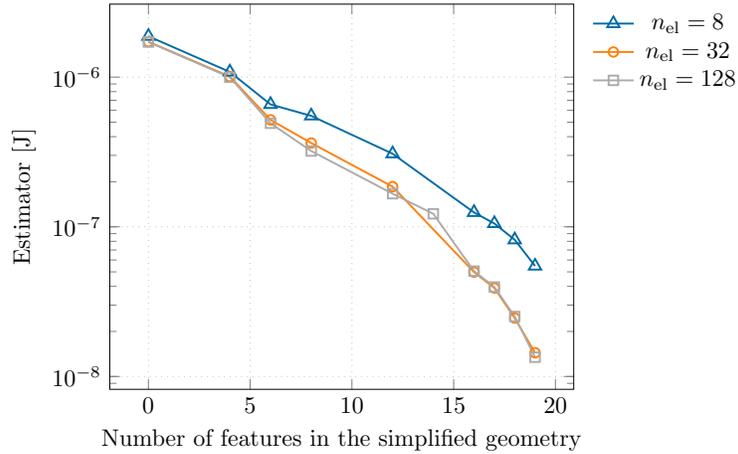
\begin{figure}
	\begin{center}
		\begin{tikzpicture}[scale=0.9]
			\begin{axis}[ymode=log, legend style={at={(1.2,1)}, legend columns=1, anchor=north, draw=none}, xlabel=Number of features in the simplified geometry, ylabel={Estimator $[\mathrm{J}]$}, grid style=dotted,grid]
				\addplot[mark=triangle, c1, thick, mark options={solid,scale=1.5}] table [x=feat8, y=est8scaled, col sep=comma] {\dataPath/solution_8_32.csv};
				\addplot[mark=o, c2, thick] table [x=feat32, y=est32scaled, col sep=comma] {\dataPath/solution_8_32.csv};
				\addplot[mark=square, c3, thick] table [x=feat128, y=est128scaled, col sep=comma] {\dataPath/solution_128.csv};
				\legend{$n_\mathrm{el} = 8$, $n_\mathrm{el} = 32$, $n_\mathrm{el} = 128$};
			\end{axis}
		\end{tikzpicture}
		\caption{Numerical test \ref{sec:3Dtest} -- Results of the adaptive defeaturing strategy for different discretization parameters.} \label{fig:adaptivedefeatresults}
	\end{center}
\end{figure}

\begin{table}
	\setlength{\tabcolsep}{4.5pt}
	\centering
	{\def\arraystretch{1.2}
		\begin{tabular}{@{}ccccccccccc@{}}
			\hline
			Iteration $i$       & $0$                & $1$       & $2$        & $3$                    & $4$                    & $5$        & $6$ & $7$ & $8$ & $9$      \\
			\hline
			$n_\mathrm{el}=8$   & $5$, $6$, $7$, $8$ & $9$, $11$ & $10$, $12$ & $14$, $15$, $18$, $19$ & $13$, $16$, $17$, $20$ & $1$        & $2$ & $3$ & $4$ & $\times$ \\
			$n_\mathrm{el}=32$  & $5$, $6$, $7$, $8$ & $9$, $11$ & $10$, $12$ & $14$, $15$, $18$, $19$ & $13$, $16$, $17$, $20$ & $3$        & $2$ & $1$ & $4$ & $\times$ \\
			$n_\mathrm{el}=128$ & $5$, $6$, $7$, $8$ & $9$, $11$ & $10$, $12$ & $14$, $15$, $18$, $19$ & $17$, $20$             & $13$, $16$ & $3$ & $2$ & $1$ & $4$      \\
			\hline
		\end{tabular}}
	\caption{Numerical test \ref{sec:3Dtest} -- Marked features at each iteration on different mesh refinements.} \label{tbl:addedfeatsdiffnel}
\end{table}

Finally, the numerical error is not considered in this article. However, it is interesting to note that the estimator is still able to drive the proposed adaptive strategy on a coarser mesh. Indeed, this algorithm has been performed on multiple meshes containing a different number $n_{\mathrm{el}}$ of elements in each space direction of the bounding box of $\Omega_0$. More precisely, we have considered $n_{\mathrm{el}}=8$, $32$, and $128$. In all three cases, the convergence of the estimator $\mathscr{E}\big(\boldsymbol{u}_\mathrm d^{(i)}\big)$ is reported in Figure~\ref{fig:adaptivedefeatresults}, and the features chosen at each iteration are reported in Table~\ref{tbl:addedfeatsdiffnel}. We can observe that except for features whose error contributions are very close to one another, the adaptive algorithm is able to correctly choose the important features, even on the coarsest mesh.

\section{Conclusions} \label{s:ccl}
In the context of the Poisson's, linear elasticity, and Stokes equations, we have studied the accuracy impact of removing \changes{distinct Neumann} features from geometries in which the solution of a PDE is sought. In particular, we have generalized the \textit{a posteriori} estimator of the energy norm of the defeaturing error from \cite{paper1defeaturing} to two- and three-dimensional geometries containing an arbitrary number of negative, positive, or generally complex features.
The proposed estimator has the following properties:
\begin{itemize}
  \item it is not only driven by geometrical considerations, but also by the PDE at hand;
  \item it is able to weight the impact of defeaturing in the energy \review{norm,} and its effectivity index is independent of the size of the geometrical features and of their number;
  \item it is able to determine whether the defeaturing error comes from the choice of defeaturing data (right hand side and Neumann boundary conditions), or if it comes from the importance of the presence of the feature itself;
  \item it is rigorously proven to be reliable and efficient up to oscillations;
  \item it is naturally decomposed into single feature contributions;
  \item it is simple, computationally cheap, and embarrassingly parallel, as it only requires the evaluation of fluxes through boundary pieces and the resolution of small problems at feature level.
  \item it has been tested on an extensive set of numerical experiments: In all of them, the estimator acts as an excellent approximation of the defeaturing error.
\end{itemize}
Note however that our framework does not include the case of a geometry whose boundary is complex everywhere as considered for instance in \cite{hiptmairshapeapprox,HEYDAROV2022102157,shapederpaper}, because of Assumption~\ref{as:separatednew}.

Then, with the help of the proposed error estimator, we have been able to design an adaptive geometric refinement strategy taking into account the defeaturing errors. More precisely, starting from a fully defeatured geometry, features are iteratively added to the geometrical model when their absence is responsible for most of the solution accuracy loss. That is, the strategy is able to build a (partially) defeatured geometric model containing few features, for which the defeaturing error is below a prescribed tolerance. Presented numerical experiments have demonstrated the convergence of the defeaturing error during the adaptive loop. In a subsequent work \cite{adaptivedefeaturing}, the proposed adaptive strategy will be combined with a mesh refinement strategy in the case in which a finite element method (and in particular isogeometric analysis \cite{igabasis,igabook}) is used to approximately solve the PDE at hand.


\section*{Acknowledgments}
The authors acknowledge the support of the European Research Council, via the ERC AdG project CHANGE n.694515.
Pablo Antol\'in also acknowledges the support of the Swiss National Science Foundation through the project “Design-through-Analysis (of PDEs): the litmus test” n.40B2-0 187094 (BRIDGE
Discovery 2019). \changes{Ondine Chanon also acknowledges the support of the Swiss National Science Foundation through the project n.P500PT\_210974.} 
Prof.\ Annalisa Buffa and Dr.\ Rafael V\'azquez are also gratefully acknowledged for the fruitful discussions on the subject.

\section*{Declarations}

\paragraph{Conflict of interest}
The authors declare that they have no conflict of interest.

\addcontentsline{toc}{chapter}{Bibliography}
\bibliography{bib2}
\bibliographystyle{ieeetr}

\appendix
\changesbis{
	\section{Proofs of reliability and efficiency}\label{app:proofs}
	The proofs of Theorems~\ref{thm:upperbound} and~\ref{thm:lowerbound} are given in this appendix in the framework of Stokes' equations. Note that the proofs in the context of linear elasticity and Poisson's problems are very similar. From the proofs for Stokes' equations, one only needs to remove the pressure terms and the divergence condition everywhere to obtain the proofs for the linear elasticity problem. One moreover needs to consider the scalar problems equivalent to the vectorial problems to obtain the proofs for Poisson's equation. \\

\review{To ease the notation in the following analysis, we respectively denote the Neumann boundaries of $\Omega_0$ and of $\tilde F_\mathrm p^k$ by
$\Gamma_N^0 := \left(\Gamma_N\setminus\gamma\right)\cup\gamma_0$, 
$\tilde \Gamma_N^k := \gamma_\intersign^k \cup \tilde \gamma^k$, for all 
$k=1,\ldots,N_f$,
and we let
$\tilde \Gamma_N := \displaystyle\bigcup_{k=1}^{N_f} \tilde \Gamma_N^k = \gamma_\intersign \cup \tilde \gamma$.} 

\subsection{Reliability}
In this section, we prove Theorem~\ref{thm:upperbound} in the context of Stokes' equations. That is, \review{under Assumptions~\ref{as:separatednew} and~\ref{assu:notrimmingdirichlet}}, we prove that the error indicator defined in~(\ref{eq:multiestimatorstokes}) is reliable, i.e., it is an upper bound for the defeaturing error.

\begin{proof}
	Let $\boldsymbol e_{\boldsymbol u} :=\boldsymbol u-\boldsymbol u_\mathrm d \in \boldsymbol H^1_{\boldsymbol 0,\Gamma_D}(\Omega)$ and $e_{p} :=p-p_\mathrm d \in L^2(\Omega)$. Let us consider the exact problem~(\ref{eqid:weakoriginalstokespb}) restricted to $\Omega_\star=\Omega\setminus \overline{F_\mathrm p}$ with the natural Neumann boundary condition $\boldsymbol \sigma(\boldsymbol u)\mathbf{n}_0- p\mathbf{n}_0$ on $\gamma_{0,\mathrm p}$, and let us consider the simplified problem~(\ref{eqid:simplstokespb}) also restricted to $\Omega_\star$, with the natural Neumann boundary condition $\boldsymbol{\sigma}(\boldsymbol u_\mathrm d)\mathbf{n}-p_\mathrm d\mathbf{n}$ on $\gamma_\mathrm n$. Then, combining both differential problems, for all $(\boldsymbol v_0,q_0)\in \boldsymbol H^1_{\boldsymbol 0,\Gamma_D}(\Omega_\star)\times L^2(\Omega_\star)$,
	\begin{align} 
	\int_{\Omega_\star}\boldsymbol{\sigma}(\boldsymbol e_{\boldsymbol u}) : \boldsymbol{\varepsilon}(\boldsymbol v_0) \,\mathrm dx - \int_{\Omega_\star} e_p\boldsymbol \nabla \cdot \boldsymbol v_0 \,\mathrm dx = &\int_{\gamma_\mathrm n} \big(\boldsymbol g-\boldsymbol{\sigma}(\boldsymbol u_\mathrm d)\mathbf n + p_\mathrm d\mathbf n\big) \cdot \boldsymbol v_0\,\mathrm ds + \int_{\gamma_{0,\mathrm p}} \big( \boldsymbol{\sigma}(\boldsymbol u)\mathbf n_0 - p\mathbf n_0 - \boldsymbol g_0\big) \cdot \boldsymbol v_0\,\mathrm ds, \nonumber \\
	-\int_{\Omega_\star} q_0\boldsymbol \nabla \cdot \boldsymbol e_{\boldsymbol u} \,\mathrm dx = \,&0. \label{eq:multieinterstokes}
	\end{align}
	In a very similar fashion, we can deduce that for all $k=1,\ldots,N_f$ and all $(\boldsymbol v^k, q^k)\in \boldsymbol H^1\big(F_\mathrm p^k\big)\times L^2\big(F_\mathrm p^k\big)$,
	\begin{align} 
	\int_{F_\mathrm p^k} \boldsymbol{\sigma}(\boldsymbol e_{\boldsymbol u}) : \boldsymbol{\varepsilon}(\boldsymbol v^k) \,\mathrm dx - \int_{F_\mathrm p^k} e_p\boldsymbol \nabla \cdot \boldsymbol v^k \,\mathrm dx = &\int_{\gamma_{0,\mathrm p}^k} \left[\big(\boldsymbol{\sigma}(\boldsymbol u)-\boldsymbol{\sigma}(\boldsymbol u_\mathrm d)\big)\mathbf n^k - (p-p_\mathrm d)\mathbf n^k\right] \cdot \boldsymbol v^k\,\mathrm ds \nonumber \\
	&+ \int_{\gamma_\setminussign^k}\big(\boldsymbol g-\boldsymbol{\sigma}(\boldsymbol u_\mathrm d)\mathbf n^k + p_\mathrm d\mathbf n^k \big) \cdot \boldsymbol v^k\,\mathrm ds, \nonumber \\
	-\int_{F_\mathrm p^k} q^k\boldsymbol \nabla \cdot \boldsymbol e_{\boldsymbol u} \,\mathrm dx = \,&0. \label{eq:multieFistokes}
	\end{align}
	Therefore, let $(\boldsymbol v, q)\in \boldsymbol H^1_{\boldsymbol 0,\Gamma_D}(\Omega)\times L^2(\Omega)$, 
	since $\mathbf n = \mathbf n^{k_\boundarypiece}$ on all $\boundarypiece \in \Sigma_\setminussign$, then 
	\begin{align}
	\mathfrak a(\boldsymbol e_{\boldsymbol u}, \boldsymbol v) + \mathfrak b(\boldsymbol v, e_p) &= \int_{\Omega} \boldsymbol{\sigma}(\boldsymbol e_{\boldsymbol u}) : \boldsymbol{\varepsilon}(\boldsymbol v) \,\mathrm{d}x - \int_\Omega e_p\boldsymbol \nabla \cdot \boldsymbol v \,\mathrm dx = \sum_{\boundarypiece \in \Sigma} \int_\boundarypiece \boldsymbol d_\boundarypiece \cdot \boldsymbol v\,\mathrm ds, \label{eq:upperdecompstokes} \\ 
	\mathfrak b(\boldsymbol e_{\boldsymbol u}, q) &= -\int_\Omega q \boldsymbol \nabla \cdot \boldsymbol e_{\boldsymbol u}\,\mathrm dx = 0. \label{eq:upperdecompstokes2}
	\end{align}
	The right hand side of (\ref{eq:upperdecompstokes}) can be rewritten 
	\begin{align}
	\sum_{\boundarypiece \in \Sigma} \int_\boundarypiece \boldsymbol d_\boundarypiece \cdot \boldsymbol v\,\mathrm ds
	&= \sum_{\boundarypiece\in\Sigma} \left[ \int_\boundarypiece \left( \boldsymbol d_\boundarypiece -\overline{\boldsymbol d_\boundarypiece}^\boundarypiece\right) \cdot \left(\boldsymbol v-\overline{\boldsymbol v}^\boundarypiece\right) \,\mathrm ds + \overline{\boldsymbol d_\boundarypiece}^\boundarypiece \cdot \int_\boundarypiece \boldsymbol v \,\mathrm ds \right]. \label{eq:upperdecomplinelast}
	\end{align}
	For each $\boundarypiece\in\Sigma$, the first terms of~(\ref{eq:upperdecompstokes}) can be estimated thanks to the Poincar\'e inequality and trace inequalities, using the domains $\Omega^k$ defined in \review{condition~\ref{it:conditionbseparated} of} Assumption~\ref{as:separatednew} for $k=1,\ldots,N_f$. That is,
	\begin{align}
	&\sum_{\boundarypiece\in\Sigma} \int_\boundarypiece \left( \boldsymbol d_\boundarypiece -\overline{\boldsymbol d_\boundarypiece}^\boundarypiece\right) \cdot \left(\boldsymbol v-\overline{\boldsymbol v}^\boundarypiece\right) \,\mathrm ds \nonumber \\
	\lesssim &\sum_{\boundarypiece\in\Sigma} \left\|\boldsymbol d_\boundarypiece-\overline{\boldsymbol d_\boundarypiece}^\boundarypiece\right\|_{0,\boundarypiece} \left\|\boldsymbol v-\overline{\boldsymbol v}^\boundarypiece\right\|_{0,\boundarypiece} \lesssim \sum_{\boundarypiece\in\Sigma} \left|\boundarypiece\right|^{\frac{1}{2(n-1)}} \left\|\boldsymbol d_\boundarypiece-\overline{\boldsymbol d_\boundarypiece}^\boundarypiece\right\|_{0,\boundarypiece} |\boldsymbol v|_{\frac{1}{2},\boundarypiece} \nonumber \\
	\lesssim &\sum_{\boundarypiece\in\Sigma_\mathrm n\cup\Sigma_\setminussign} \left|\boundarypiece\right|^{\frac{1}{2(n-1)}} \left\|\boldsymbol d_\boundarypiece-\overline{\boldsymbol d_\boundarypiece}^\boundarypiece\right\|_{0,\boundarypiece} \|\boldsymbol v\|_{1,\Omega^{k_\boundarypiece}} + \sum_{\boundarypiece\in\Sigma_{0,\mathrm p}} \left|\boundarypiece\right|^{\frac{1}{2(n-1)}} \left\|\boldsymbol d_\boundarypiece-\overline{\boldsymbol d_\boundarypiece}^\boundarypiece\right\|_{0,\boundarypiece} \|\boldsymbol v\|_{1,\Omega^{k_\boundarypiece}\cap\Omega_\star}. \label{eq:firsttermlowerlinelast}
	\end{align}
	Then, the last terms of (\ref{eq:upperdecomplinelast}) can be estimated thanks to \cite[Appendix~A.2]{paper1defeaturing} and trace inequalities, that is, 
	\begin{align}
	\overline{\boldsymbol d_\boundarypiece}^\boundarypiece \cdot \int_\boundarypiece \boldsymbol v \,\mathrm ds \lesssim \,&\,|\boundarypiece|^\frac{1}{2}\left\| \overline{\boldsymbol d_\boundarypiece}^\boundarypiece \right\|_{\ell^2} \left\|\boldsymbol v\right\|_{0,\boundarypiece} \nonumber \\
	\lesssim &\sum_{\boundarypiece\in\Sigma_\mathrm n\cup \Sigma_\setminussign} c_\boundarypiece |\boundarypiece|^\frac{n}{2(n-1)} \left\| \overline{\boldsymbol d_\boundarypiece}^\boundarypiece \right\|_{\ell^2} \left\|\boldsymbol v\right\|_{\frac{1}{2},\partial \Omega^{k_\boundarypiece}} + \sum_{\boundarypiece\in\Sigma_{0,\mathrm p}} c_\boundarypiece |\boundarypiece|^\frac{n}{2(n-1)} \left\| \overline{\boldsymbol d_\boundarypiece}^\boundarypiece \right\|_{\ell^2} \left\|\boldsymbol v\right\|_{\frac{1}{2},\partial \left(\Omega^{k_\boundarypiece}\cap\Omega_\star\right)} \nonumber \\
	\lesssim &\sum_{\boundarypiece\in\Sigma_\mathrm n\cup \Sigma_\setminussign} c_\boundarypiece |\boundarypiece|^\frac{n}{2(n-1)} \left\| \overline{\boldsymbol d_\boundarypiece}^\boundarypiece \right\|_{\ell^2} \left\|\boldsymbol v\right\|_{1, \Omega^{k_\boundarypiece}} + \sum_{\boundarypiece\in\Sigma_{0,\mathrm p}} c_\boundarypiece |\boundarypiece|^\frac{n}{2(n-1)} \left\| \overline{\boldsymbol d_\boundarypiece}^\boundarypiece \right\|_{\ell^2} \left\|\boldsymbol v\right\|_{1,\Omega^{k_\boundarypiece}\cap\Omega_\star}. \label{eq:secondtermlowerlinelast}
	\end{align}
	Thus combining~(\ref{eq:upperdecomplinelast}), (\ref{eq:firsttermlowerlinelast}) and~(\ref{eq:secondtermlowerlinelast}), defining $\hat{\mathscr{E}}(\boldsymbol u_\mathrm d, p_\mathrm d) := \mu^\frac{1}{2}\mathscr{E}(\boldsymbol u_\mathrm d, p_\mathrm d)$, for all $\boldsymbol v \in \boldsymbol H^1_{\boldsymbol 0,\Gamma_D}(\Omega)$, 
	\begin{equation}
	\sum_{\boundarypiece \in \Sigma} \int_\boundarypiece \boldsymbol d_\boundarypiece \cdot \boldsymbol v\,\mathrm ds \lesssim \hat{\mathscr{E}}(\boldsymbol u_\mathrm d, p_\mathrm d) \|\boldsymbol \nabla \boldsymbol v\|_{0,\Omega}. \label{eq:estdsigma}
	\end{equation}
	Now, remark that if we take $\boldsymbol v = \boldsymbol e_{\boldsymbol u} \in \boldsymbol H^1_{\boldsymbol 0,\Gamma_D}(\Omega)$ and $q=e_p\in L^2(\Omega)$, then equation~(\ref{eq:upperdecompstokes2}) reads $\mathfrak b(\boldsymbol e_{\boldsymbol u},e_p) = 0$, and thus using~(\ref{eq:estdsigma}), equation~(\ref{eq:upperdecompstokes}) rewrites
	\begin{align*}
	\mathfrak a(\boldsymbol e_{\boldsymbol u}, \boldsymbol e_{\boldsymbol u}) &= \int_{\Omega} \boldsymbol{\sigma}(\boldsymbol e_{\boldsymbol u}) : \boldsymbol{\varepsilon}(\boldsymbol e_{\boldsymbol u}) \,\mathrm{d}x - \int_\Omega e_p\boldsymbol \nabla \cdot \boldsymbol e_{\boldsymbol u} \,\mathrm dx = \sum_{\boundarypiece \in \Sigma} \int_\boundarypiece \boldsymbol d_\boundarypiece \cdot \boldsymbol e_{\boldsymbol u}\,\mathrm ds \lesssim \hat{\mathscr{E}}(\boldsymbol u_\mathrm d, p_\mathrm d) \|\boldsymbol \nabla \boldsymbol e_{\boldsymbol u}\|_{0,\Omega}. 
	\end{align*}
	Using the coercivity of $\mathfrak a(\boldsymbol \cdot,\boldsymbol \cdot)$ in $\boldsymbol H^1_{\boldsymbol 0,\Gamma_D}(\Omega)$ equipped with the norm $\|\boldsymbol \nabla \cdot \|_{0,\Omega}$, then 
	\begin{align} 
	\mathfrak a(\boldsymbol e_{\boldsymbol u}, \boldsymbol e_{\boldsymbol u}) &\lesssim \hat{\mathscr{E}}(\boldsymbol u_\mathrm d, p_\mathrm d) \|\boldsymbol{\nabla} \boldsymbol e_{\boldsymbol u}\|_{0,\Omega}\lesssim \mu^{-\frac{1}{2}}\hat{\mathscr{E}}(\boldsymbol u_\mathrm d, p_\mathrm d) \big( \mathfrak a(\boldsymbol e_{\boldsymbol u}, \boldsymbol e_{\boldsymbol u}) \big)^\frac{1}{2} 
	\label{eq:usingcoercivitya}
	\end{align}
	so that if we simplify on both sides, 
	\begin{equation}
	\big( \mathfrak a(\boldsymbol e_{\boldsymbol u}, \boldsymbol e_{\boldsymbol u}) \big)^\frac{1}{2} \lesssim \mathscr{E}(\boldsymbol u_\mathrm d, p_\mathrm d). \label{eq:esteu}
	\end{equation}
	
	Finally, since $\mathfrak b(\boldsymbol\cdot,\cdot)$ satisfies the inf-sup condition, using (\ref{eq:upperdecompstokes}) and (\ref{eq:estdsigma}), using the continuity of $\mathfrak a(\boldsymbol \cdot, \boldsymbol\cdot)$ in $\boldsymbol H^1_{\boldsymbol 0,\Gamma_D}(\Omega)$, and then its coercivity as in (\ref{eq:usingcoercivitya}), then
	\begin{align}
	\|e_p\|_{0,\Omega} \lesssim \sup_{\substack{\boldsymbol v\in \boldsymbol H^1_{\boldsymbol 0,\Gamma_D}(\Omega)\\\boldsymbol v\neq \boldsymbol 0}} \frac{\mathfrak b(\boldsymbol v, e_p)}{\|\boldsymbol \nabla \boldsymbol v\|_{0,\Omega}} &= \sup_{\substack{\boldsymbol v\in \boldsymbol H^1_{\boldsymbol 0,\Gamma_D}(\Omega)\\\boldsymbol v\neq \boldsymbol 0}} \frac{\displaystyle\sum_{\boundarypiece \in \Sigma} \displaystyle\int_\boundarypiece \boldsymbol d_\boundarypiece \cdot \boldsymbol v\,\mathrm ds - \mathfrak a(\boldsymbol e_{\boldsymbol u}, \boldsymbol v)}{\|\boldsymbol \nabla \boldsymbol v\|_{0,\Omega}} \nonumber \\
	&\leq \sup_{\substack{\boldsymbol v\in \boldsymbol H^1_{\boldsymbol 0,\Gamma_D}(\Omega)\\\boldsymbol v\neq \boldsymbol 0}} \frac{\displaystyle\sum_{\boundarypiece \in \Sigma} \displaystyle\int_\boundarypiece \boldsymbol d_\boundarypiece \cdot \boldsymbol v\,\mathrm ds }{\|\boldsymbol \nabla \boldsymbol v\|_{0,\Omega}} - \inf_{\substack{\boldsymbol v\in \boldsymbol H^1_{\boldsymbol 0,\Gamma_D}(\Omega)\\\boldsymbol v\neq \boldsymbol 0}} \frac{\mathfrak a(\boldsymbol e_{\boldsymbol u}, \boldsymbol v)}{\|\boldsymbol \nabla \boldsymbol v\|_{0,\Omega}}\nonumber \\
	&\lesssim \hat{\mathscr E}(\boldsymbol u_\mathrm d, p_\mathrm d) + \|\boldsymbol \nabla \boldsymbol e_{\boldsymbol u}\|_{0,\Omega} \lesssim \hat{\mathscr E}(\boldsymbol u_\mathrm d, p_\mathrm d) + \mu^{\frac{1}{2}}\big( \mathfrak a(\boldsymbol e_{\boldsymbol u}, \boldsymbol e_{\boldsymbol u}) \big)^\frac{1}{2} \lesssim \hat{\mathscr E}(\boldsymbol u_\mathrm d, p_\mathrm d). \label{eq:estep}
	\end{align}
	We can therefore conclude by combining (\ref{eq:esteu}) and (\ref{eq:estep}). 
\end{proof}

\subsection{Efficiency}
In this section, we prove Theorem~\ref{thm:lowerbound} in the context of Stokes' equations. That is, we prove that the error indicator defined in~(\ref{eq:multiestimatorstokes}) is efficient \review{under Assumptions~\ref{as:separatednew} and~\ref{assu:notrimmingdirichlet}}, i.e., it is a lower bound for the defeaturing error, up to oscillations. \review{The proof is given in the general case for $n=3$, but under the following additional assumption if $n=2$:
\begin{assumption}\label{as:compatcond2d}
	If $n=2$, assume that the compatibility conditions~(\ref{eq:compatcondmultilinelast}) and~(\ref{eq:compatcondmultilinelast2}) are satisfied. 
\end{assumption}
Note however that the numerical experiments show that the result also holds without compatibility conditions if $n=2$. The same restrictions were considered in the efficiency proofs of~\cite{paper1defeaturing} in the two-dimensional case.}\\

In the following, for all $D\subset \mathbb R^n$ and $\Lambda\subset \partial D$, we denote by $H_{00}^\frac{1}{2}(\Lambda)$ the space of functions $w\in L^2(\Lambda)$ such that the zero-extension $w^\star$ of $w$ to $\partial D$ belongs to $H^\frac{1}{2}(\partial D)$, and we denote by $H_{00}^{-\frac{1}{2}}(\Lambda)$ its dual space.
For the sake of simplicity, let us strengthen Definition~\ref{defn:pwsmoothshaperegweak}, even though the proofs could be easily generalized to the weaker setting.
\begin{definition}\label{defn:pwsmoothshaperegstrong}
	An $(n-1)$-dimensional \review{subset} $\Lambda$ of $\mathbb R^n$ is \emph{regular} if $\Lambda$ is piecewise shape regular and composed of flat elements, that is, if there is $N_\Lambda\in\mathbb{N}$ such that for all $\ell_1,\ell_2=1,\ldots,N_\Lambda$ with $\ell_1\neq \ell_2$, $\Lambda = \mathrm{int}\left(\displaystyle\bigcup_{\ell=1}^{N_\Lambda} \overline{\Lambda^\ell}\right)$, $\Lambda^{\ell_1} \cap \Lambda^{\ell_2} = \emptyset$, $|\Lambda|\lesssim \left|\Lambda^{\ell_1}\right|$ and $\Lambda^{\ell_1}$ is flat, i.e. it is a straight line if $n=2$ or a flat square or triangle if $n=3$.
\end{definition}

\review{That is, we suppose that the boundaries $\boundarypiece\in\Sigma$ of the features are shape regular as in Definition~\ref{defn:pwsmoothshaperegstrong} instead of Definition~\ref{defn:pwsmoothshaperegweak}, for simplicity. This allows us to easily define the following Cl\'ement operator. Whenever some boundary $\Lambda$ is regular, then for all $m\in\mathbb N$, let $\mathbb Q_{m,0}^\mathrm{pw}(\Lambda)$ be the space of continuous piecewise polynomials of degree at most $m$ on each variable, that vanish at the boundary $\partial \Lambda$. Then we define
\begin{equation} \label{eq:clement}
\Pi_{m,\Lambda} : L^2(\Lambda) \to \mathbb Q_{m,0}^\mathrm{pw}(\Lambda)
\end{equation}
as the extension of the Cl\'ement operator \cite{clement} developed in \cite{bernardigirault} on $\Lambda$.}\\

In this context, we define the oscillations appearing in the upper bound of Theorem~\ref{thm:lowerbound} as follows:
\begin{definition}\label{def:oscmultistokes}
	For any $m\in \mathbb{N}$, let $\boldsymbol\Pi_{m}$ be such that $\boldsymbol\Pi_{m}\vert_{\boundarypiece} \equiv \boldsymbol\Pi_{m,\boundarypiece}$ for all $\boundarypiece\in\Sigma$, where $\boldsymbol\Pi_{m,\boundarypiece}$ is the component-wise extensions of the Cl\'ement operator defined in (\ref{eq:clement}), and let $\boldsymbol d^k$ be such that $\boldsymbol d^k\vert_\boundarypiece \equiv \boldsymbol d_\boundarypiece$ on all $\boundarypiece\in\Sigma^k$, for all $k=1,\ldots,N_f$. Then we define
	\begin{align*}
	&\mathrm{osc}(\boldsymbol u_\mathrm d, p_\mathrm d)^2 := \sum_{k=1}^{N_f} \left( \mathrm{osc}^k(\boldsymbol u_\mathrm d, p_\mathrm d)\right)^2,
	\text{ with } \, \mathrm{osc}^k(\boldsymbol u_\mathrm d, p_\mathrm d) := \left|\Gamma^k\right|^\frac{1}{2(n-1)} \left\| \boldsymbol d^k - \boldsymbol\Pi_{m}\left(\boldsymbol d^k\right)\right\|_{0,\Gamma^k} \, \text{ for } k=1,\ldots,N_f. 
	\end{align*}
\end{definition}


\review{Let us finally give the proof of Theorem~\ref{thm:lowerbound} under Assumption~\ref{as:compatcond2d}.}
\begin{proof}
	Let $\boldsymbol e_{\boldsymbol u} :=\boldsymbol u-\boldsymbol u_\mathrm d \in \boldsymbol H^1_{\boldsymbol 0,\Gamma_D}(\Omega)$, let $e_{p} :=p-p_\mathrm d \in L^2(\Omega)$, let $k\in\{1,\ldots,N_f\}$, and let $\Omega^k_\star:= \Omega_\star \cap \Omega^k$, where $\Omega^k$ is the domain associated to feature $F^k$ defined in \review{condition~\ref{it:conditionbseparated} of} Assumption~\ref{as:repeatseparatednew}. 
	Then, let us consider the exact problem (\ref{eqid:weakoriginalstokespb}) restricted to $\Omega^k_\star$ with the natural Neumann boundary condition $\boldsymbol{\sigma}(\boldsymbol u)\mathbf{n}_0 -p\mathbf{n}_0$ on $\gamma_{0,\mathrm p}^k$ and the natural Dirichlet boundary condition $\mathrm{tr}_{\partial \Omega_\star^k\setminus \partial\Omega_\star}(\boldsymbol u)$ on $\partial \Omega^k_\star \setminus \partial \Omega_\star$, and let us consider the simplified problem~(\ref{eqid:simplstokespb}) also restricted to $\Omega_\star^k$, with the natural Neumann boundary condition $\boldsymbol{\sigma}(u_\mathrm d)\mathbf{n}-p_\mathrm d\mathbf n$ on $\gamma_\mathrm n^k$, and the natural Dirichlet boundary condition $\mathrm{tr}_{\partial \Omega_\star^k\setminus \partial\Omega_\star}(\boldsymbol u_\mathrm d)$ on $\partial \Omega^k_\star \setminus \partial \Omega_\star$. Then, combining both differential problems, for all $(\boldsymbol v_0,q_0)\in \boldsymbol H^1_{\boldsymbol 0,\partial \Omega^k_\star \setminus \left[ \Gamma_N\cup\gamma_{0,\mathrm p}^k\right]}(\Omega^k_\star)\times L^2(\Omega^k_\star)$, \vspace{-0.2cm}
	\begin{align}
	\int_{\Omega^k_\star} \boldsymbol{\sigma}(\boldsymbol e_{\boldsymbol u}) : \boldsymbol{\varepsilon}(\boldsymbol v_0) \,\mathrm dx - \int_{\Omega^k_\star} e_p \boldsymbol \nabla \cdot \boldsymbol v_0\,\mathrm dx = &\int_{\gamma_\mathrm n^k} \hspace{-0.05cm}\big(\boldsymbol g-\boldsymbol{\sigma}(\boldsymbol u_\mathrm d)\mathbf n + p_\mathrm d\mathbf n\big) \cdot \boldsymbol v_0\,\mathrm ds + \int_{\gamma_{0,\mathrm p}^k} \hspace{-0.05cm}\big( \boldsymbol{\sigma}(\boldsymbol u)\mathbf n_0 - p\mathbf n_0 - \boldsymbol g_0 \big)\cdot \boldsymbol v_0 \,\mathrm ds, \nonumber\\
	-\int_{\Omega^k_\star} q_0\boldsymbol \nabla \cdot \boldsymbol e_{\boldsymbol u} \,\mathrm dx = \,&0. \label{eq:omegakinteromega0stokes}
	\end{align}
	Let $(\boldsymbol v,q)\in \boldsymbol H^1_{\boldsymbol 0,\partial \Omega^k\setminus \Gamma_N}(\Omega^k)\times L^2(\Omega^k)$, and recall that $\Omega_\star = \Omega\setminus \overline{F_\mathrm p}$, so that $\Omega^k = \mathrm{int}\left(F_\mathrm p^k \cup \Omega_\star^k\right)$. Consequently, reusing equation~(\ref{eq:multieFistokes}),
	\begin{align}
	&\int_{\Omega^k} \boldsymbol{\sigma}(\boldsymbol e_{\boldsymbol u}) : \boldsymbol{\varepsilon}(\boldsymbol v) \,\mathrm dx - \int_{\Omega^k} e_p \boldsymbol \nabla \cdot \boldsymbol v\,\mathrm dx \nonumber \\
	= &\int_{\gamma_\mathrm n^k} \big( \boldsymbol g-\boldsymbol{\sigma}(\boldsymbol u_\mathrm d)\mathbf n + p_\mathrm d\mathbf n \big)\cdot \boldsymbol v\,\mathrm ds + \int_{\gamma_{0,\mathrm p}^k} \big(\hspace{-1mm}-\boldsymbol g_0-\boldsymbol{\sigma}(\boldsymbol u_\mathrm d)\mathbf n^k + p_\mathrm d\mathbf n^k \big) \cdot \boldsymbol v \,\mathrm ds + \int_{\gamma_\setminussign^k} \big( \boldsymbol g-\boldsymbol{\sigma}(\boldsymbol u_\mathrm d)\mathbf n^k + p_\mathrm d\mathbf n^k \big)\cdot\boldsymbol v\,\mathrm ds \nonumber\\
	= &\sum_{\boundarypiece \in \Sigma^k} \int_\boundarypiece \boldsymbol d_\boundarypiece \cdot \boldsymbol v \,\mathrm ds, \qquad \qquad \text{ and } \qquad -\int_{\Omega^k} q\boldsymbol \nabla \cdot \boldsymbol e_{\boldsymbol u} \,\mathrm dx = 0. \label{eq:beforecontinuitystokes}
	\end{align}
	Now, let 
	$\mathfrak a^k(\boldsymbol \cdot, \boldsymbol \cdot):\boldsymbol H^1_{\boldsymbol 0,\partial \Omega^k\setminus\Gamma_N}(\Omega^k)\times\boldsymbol H^1_{\boldsymbol 0,\partial \Omega^k\setminus\Gamma_N}(\Omega^k)\to\mathbb R$ and $\mathfrak b^k(\boldsymbol \cdot, \cdot):\boldsymbol H^1_{\boldsymbol 0,\partial \Omega^k\setminus\Gamma_N}(\Omega^k)\times L^2(\Omega^k)\to\mathbb R$
	be defined by \vspace{-0.2cm}
	\begin{align*}
	\mathfrak a^k(\boldsymbol w,\boldsymbol v) &= \int_{\Omega^k} \boldsymbol{\sigma}(\boldsymbol w) : \boldsymbol{\varepsilon}(\boldsymbol v) \,\mathrm dx, && \forall \boldsymbol{w},\boldsymbol{v}\in \boldsymbol H^1_{\boldsymbol 0,\partial \Omega^k\setminus\Gamma_N}(\Omega^k), \\
	\mathfrak b^k(\boldsymbol v, q) &= -\int_{\Omega^k} q\boldsymbol{\nabla} \cdot \boldsymbol v \,\mathrm dx, && \forall \boldsymbol{v}\in \boldsymbol H^1_{\boldsymbol 0,\partial \Omega^k\setminus\Gamma_N}(\Omega^k), \forall q\in L^2(\Omega^k).
	\end{align*}
	Note that $\mathfrak a^k(\boldsymbol \cdot, \boldsymbol \cdot)$ and $\mathfrak b^k(\boldsymbol \cdot,\cdot)$ are continuous with respect to the norms $\|\boldsymbol \nabla \boldsymbol \cdot \|_{0,\Omega^k}$ for $\boldsymbol H^1_{\boldsymbol 0,\partial \Omega^k\setminus\Gamma_N}(\Omega^k)$, and $\|\cdot\|_{0,\Omega^k}$ for $L^2(\Omega^k)$. 
	Thus using~(\ref{eq:beforecontinuitystokes}), for all $\boldsymbol v\in \boldsymbol H^1_{\boldsymbol 0,\partial \Omega^k\setminus\Gamma_N}(\Omega^k)$, 
	\begin{equation} \label{eq:errOmegakstokes}
	\sum_{\boundarypiece \in \Sigma^k} \int_\boundarypiece \boldsymbol d_\boundarypiece \cdot \boldsymbol v \,\mathrm ds = \mathfrak a^k(\boldsymbol e_{\boldsymbol u}, \boldsymbol v) + \mathfrak b^k(\boldsymbol v, e_p)  \lesssim \mu \left(\|\boldsymbol \nabla \boldsymbol e_{\boldsymbol u}\|_{0,\Omega^k} + \|e_p\|_{0,\Omega^k}\right) \|\boldsymbol \nabla \boldsymbol v\|_{0,\Omega^k}.
	\end{equation}
	
	Furthermore, let $\boldsymbol H^{(k)} := \left\{ \boldsymbol v\in \boldsymbol H_{00}^\frac{1}{2}\left(\Gamma^k\right) : \boldsymbol v|_\boundarypiece \in \boldsymbol H_{00}^\frac{1}{2}(\boundarypiece), \text{ for all } \boundarypiece\in\Sigma^k \right\}$ equipped with the norm 
	\begin{equation*} 
	\|\cdot\|_{\boldsymbol H^{(k)}} := \left( \displaystyle\sum_{\boundarypiece\in\Sigma^k} \|\cdot\|_{\boldsymbol H_{00}^{1/2}(\boundarypiece)}^2 \right)^\frac{1}{2},
	\end{equation*}
	and let $\left(\boldsymbol H^{(k)}\right)^*$ be its dual space equipped with the dual norm $\|\cdot\|_{\left(\boldsymbol H^{(k)}\right)^*}$. For all $\boldsymbol w\in \boldsymbol H^{(k)}$, let us define piecewise $\boldsymbol u_{\boldsymbol w}\in \boldsymbol H^1_{\boldsymbol 0,\partial \Omega^k\setminus\left(\gamma^k_\mathrm n \cup \gamma^k_\setminussign\right)}\left(\Omega^k\right) \subset \boldsymbol H^1_{\boldsymbol 0,\partial \Omega^k\setminus \Gamma_N}\left(\Omega^k\right)$ as the unique solution of
	\begin{align*}
	\begin{cases}
	-\boldsymbol\nabla \cdot \tensorsigma\left(\boldsymbol u_{\boldsymbol w}|_{F_\mathrm p^k} \right) = \boldsymbol 0 & \text{in } F_\mathrm p^k \\
	\boldsymbol u_{\boldsymbol w}|_{F_\mathrm p^k} = \left(\boldsymbol w|_{\gamma_{0,\mathrm p}^k\cup\gamma_\setminussign^k}\right)^\star & \text{on } \partial F_\mathrm p^k, 
	\end{cases}
	\qquad
	\begin{cases}
	-\boldsymbol\nabla \cdot \tensorsigma\left(\boldsymbol u_{\boldsymbol w}|_{\Omega^k_\star} \right) = \boldsymbol 0 & \text{in } \Omega_\star^k \\
	\boldsymbol u_{\boldsymbol w}|_{\Omega_\star^k} = \left(\boldsymbol w|_{\gamma_{0,\mathrm p}^k\cup\gamma_\mathrm n^k}\right)^\star & \text{on } \partial \Omega_\star^k, 
	\end{cases}
	\end{align*} 
	where $\left(\boldsymbol w|_{\gamma_{0,\mathrm p}^k\cup\gamma_\setminussign^k}\right)^\star$ and $\left(\boldsymbol w|_{\gamma_{0,\mathrm p}^k\cup\gamma_\mathrm n^k}\right)^\star$ are the extensions by $\boldsymbol 0$ of $\boldsymbol w|_{\gamma_{0,\mathrm p}^k\cup\gamma_\setminussign^k}$ on $\partial F_\mathrm p^k$ and of $\boldsymbol w|_{\gamma_{0,\mathrm p}^k\cup\gamma_\mathrm n^k}$ on $\partial\Omega_\star^k$, respectively.
	Then by continuity of the solution on the data and from \cite[Appendix~A.6]{paper1defeaturing},
	\begin{align}
	\left\|\boldsymbol \nabla \boldsymbol u_{\boldsymbol w}\right\|_{0,\Omega^k} \lesssim \left( \|\boldsymbol w\|_{\boldsymbol H_{00}^{1/2}\left(\gamma_{0,\mathrm p}^k\cup\gamma_\setminussign^k\right)}^2 + \|\boldsymbol w\|_{\boldsymbol H_{00}^{1/2}\left(\gamma_{0,\mathrm p}^k\cup\gamma_\mathrm n^k\right)}^2 \right)^\frac{1}{2} \lesssim \|\boldsymbol w\|_{\boldsymbol H^{(k)}}. \label{eq:H(k)linelast}
	\end{align}
	So thanks to (\ref{eq:errOmegakstokes}) and (\ref{eq:H(k)linelast}), recalling that $\boldsymbol d^k\vert_\boundarypiece = \boldsymbol d_\boundarypiece$ on each $\boundarypiece\in \Sigma^k$ by definition, then
	\begin{align}
	\left\| \boldsymbol d^k\right\|_{\left(\boldsymbol H^{(k)}\right)^*} = \sup_{\substack{\boldsymbol w\in \boldsymbol H^{(k)}\\ \boldsymbol w\neq \boldsymbol 0}} \frac{\displaystyle \int_{\Gamma^k} \boldsymbol d^k \cdot \boldsymbol w \,\mathrm ds}{\|\boldsymbol w\|_{\boldsymbol H^{(k)}}}
	&\lesssim \sup_{\substack{\boldsymbol w\in \boldsymbol H^{(k)}\\ \boldsymbol w\neq \boldsymbol 0}} \frac{\displaystyle \sum_{\boundarypiece\in\Sigma^k} \int_{\boundarypiece} \boldsymbol d_\boundarypiece \cdot \boldsymbol u_{\boldsymbol w}\,\mathrm ds}{\left\|\boldsymbol \nabla \boldsymbol u_{\boldsymbol w}\right\|_{0,\Omega^k}} \nonumber \\
	& \leq \sup_{\substack{\boldsymbol v\in \boldsymbol H^1_{\boldsymbol 0,\partial \Omega^k\setminus\Gamma_N}\left(\Omega^k\right)\\ \boldsymbol v\neq \boldsymbol 0}} \frac{\displaystyle \sum_{\boundarypiece\in\Sigma^k} \int_{\boundarypiece} \boldsymbol d_\boundarypiece \cdot \boldsymbol v\,\mathrm ds}{\left\|\boldsymbol \nabla \boldsymbol v\right\|_{0,\Omega^k}} \lesssim \mu \|\boldsymbol \nabla \boldsymbol e\|_{0,\Omega^k}. \label{eq:decomperrdualOmegaklinelast}
	\end{align}
	
	Moreover, using Remark \ref{rmk:estmultitildestokes} if $n=3$, or Remark \ref{rmk:compatcondmultistokes} if $n=2$ and the data compatibility conditions~(\ref{eq:compatcondmultilinelast}) and~(\ref{eq:compatcondmultilinelast2}) are satisfied, then
	$$\mathscr{E}(\boldsymbol u_\mathrm d) \lesssim \mu^{-\frac{1}{2}}\left( \sum_{k=1}^{N_f} \sum_{\boundarypiece\in\Sigma^k} |\boundarypiece|^\frac{1}{n-1} \left\|\boldsymbol d_\boundarypiece\right\|_{0,\boundarypiece}^2 \right)^\frac{1}{2}.$$
	Therefore, using the triangle inequality and since $\left|\gamma_\mathrm n^k\right| \simeq \left|\gamma_\setminussign^k\right| \simeq \left| \gamma_{0,\mathrm p}^k \right| \simeq \left| \Gamma^k\right|$ for all $k=1,\ldots,N_f$, then
	\begin{align*}
	\sum_{\boundarypiece \in \Sigma^k} |\boundarypiece|^\frac{1}{n-1} \left\|\boldsymbol d_\boundarypiece\right\|_{0,\boundarypiece}^2 &\lesssim \sum_{\boundarypiece \in \Sigma^k} |\boundarypiece|^\frac{1}{n-1} \left( \left\|\boldsymbol\Pi_{m}(\boldsymbol d_\boundarypiece)\right\|_{0,\boundarypiece}^2 + \left\|\boldsymbol d_\boundarypiece - \boldsymbol\Pi_{m}(\boldsymbol d_\boundarypiece)\right\|_{0,\boundarypiece}^2 \right)\\
	&\lesssim \left|\Gamma^k\right|^\frac{1}{n-1} \left\|\boldsymbol\Pi_{m}\left(\boldsymbol d^k\right)\right\|_{0,\Gamma^k}^2 + \left|\Gamma^k\right|^\frac{1}{n-1} \left\|\boldsymbol d^k - \boldsymbol\Pi_{m}\left(\boldsymbol d^k\right)\right\|_{0,\Gamma^k}^2.
	\end{align*}
	Now, we use the definition of the broken norm in $\boldsymbol H^{(k)}$ to apply the inverse inequality of \cite[Appendix A.7]{paper1defeaturing}. Recalling Definition~\ref{def:oscmultistokes} and using again the triangle inequality, we thus obtain for all $k=1,\ldots,N_f$,
	\begin{align*}
	\sum_{\boundarypiece \in \Sigma^k} |\boundarypiece|^\frac{1}{n-1} \left\|\boldsymbol d_\boundarypiece\right\|_{0,\boundarypiece}^2 &\lesssim \left\| \boldsymbol\Pi_m\left(\boldsymbol d^k\right)\right\|_{\left(\boldsymbol H^{(k)}\right)^*}^2 + \mu\left( \text{osc}^k(\boldsymbol u_\mathrm d)\right)^2 \\
	&\lesssim \left\|\boldsymbol d^k \right\|^2_{\left(\boldsymbol H^{(k)}\right)^*} + \left\| \boldsymbol\Pi_m\left(\boldsymbol d^k\right) - \boldsymbol d^k \right\|^2_{\left(\boldsymbol H^{(k)}\right)^*} + \mu\left(\text{osc}^k(\boldsymbol u_\mathrm d)\right)^2.
	\end{align*}
	Finally, using (\ref{eq:decomperrdualOmegaklinelast}), applying \cite[Appendix A.4 and A.6]{paper1defeaturing}, we obtain
	\begin{align*}
	\mathscr E(\boldsymbol u_\mathrm d)^2 &\lesssim \mu^{-1}\left( \sum_{k=1}^{N_f} \left\| \boldsymbol d^k\right\|^2_{\left(\boldsymbol H^{(k)}\right)^*} + \sum_{k=1}^{N_f} \left\|\boldsymbol\Pi_m \left(\boldsymbol d^k\right) - \boldsymbol d^k\right\|^2_{\left(\boldsymbol H^{(k)}\right)^*} + \mu\,\mathrm{osc}(\boldsymbol u_\mathrm d)^2\right) \\
	& \lesssim \mu^{-1}\left(\mu^2 \sum_{k=1}^{N_f} \|\boldsymbol \nabla \boldsymbol e\|^2_{0,\Omega^k} + \sum_{k=1}^{N_f} \left\|\boldsymbol\Pi_m \left( \boldsymbol d^k\right) - \boldsymbol d^k \right\|^2_{\boldsymbol H_{00}^{-1/2}(\Gamma^k)} + \mu\,\mathrm{osc}(\boldsymbol u_\mathrm d)^2\right) \\
	& \lesssim \mu\|\boldsymbol \nabla \boldsymbol e\|^2_{0,\Omega} + \mathrm{osc}(\boldsymbol u_\mathrm d)^2 \lesssim \big( \mu^{\frac{1}{2}}\|\boldsymbol \nabla \boldsymbol e\|_{0,\Omega} + \mathrm{osc}(\boldsymbol u_\mathrm d)\big)^2. 
	\end{align*}
	To conclude, we use the coercivity of the bilinear form $\mathfrak a(\boldsymbol \cdot, \boldsymbol \cdot)$ in $\boldsymbol H^1_{\boldsymbol 0,\Gamma_D}(\Omega)$ to obtain $\mu^\frac{1}{2}\|\boldsymbol \nabla \boldsymbol e\|_{0,\Omega} \lesssim \vertiii{\boldsymbol e}_\Omega.$
\end{proof}

}

\end{document}